\documentclass[a4paper]{article} 
\usepackage{amsthm,amsmath,amssymb,mathrsfs,url,amsfonts} 

\theoremstyle{definition}
 \theoremstyle{definition}
 \newtheorem{theorem}{Theorem}[section]
 \newtheorem*{theorem*}{Theorem}
 \newtheorem{lemma}[theorem]{Lemma}
 \newtheorem{proposition}[theorem]{Proposition}
 \newtheorem{corollary}[theorem]{Corollary}
 \newtheorem{definition}{Definition}[section]
 \newtheorem{remark}{Remark}[section]

\newcommand\DN{\newcommand}
\numberwithin{equation}{section}
\numberwithin{theorem}{section}

\DN\lref[1]{Lemma~\ref{#1}}
\DN\tref[1]{Theorem~\ref{#1}}
\DN\pref[1]{Proposition~\ref{#1}}
\DN\sref[1]{Section~\ref{#1}}
\DN\dref[1]{Definition~\ref{#1}}
\DN\rref[1]{Remark~\ref{#1}} 
\DN\corref[1]{Corollary~\ref{#1}}
\DN\eref[1]{Example~\ref{#1}}

\newcounter{Const} 
\DN\Ct{\refstepcounter{Const}c_{\theConst}}
\DN\cref[1]{c_{\ref{#1}}}

\DN\tolaw{\stackrel{d}{\longrightarrow }}
\DN\ton{\stackrel{\n}{\longto}}

\DN\bs{\bigskip}\DN\ms{\medskip}
\DN\PD[2]{\frac{\partial #1 }{\partial #2}}
\DN\half{\frac{1}{2}}
\DN\map[3]{#1:#2 \to #3}
\DN\st{\,;\,}
\DN\ot{\otimes}\DN\ts{\times}

\DN\elaw{\stackrel{\mathrm{law}}{=}}
\DN\laweq{\stackrel{d}{=}}

\DN\R{\mathbb{R}}
\DN\Rd{\mathbb{R}^d}
\DN\N{\mathbb{N}}
\DN\Q{\mathbb{Q}}
\DN\Z{\mathbb{Z}}
\DN\C{\mathbb{C}}

\DN\nN{N}
\DN\rR{R}
\DN\sS{S}

\DN{\limi}[1]{\lim_{#1\to\infty}} 	
\DN{\limz}[1]{\lim_{#1\to0}}
\DN{\limdz}[1]{\lim_{#1\downarrow 0}}
\DN{\limsupi}[1]{\limsup_{#1\to\infty}}
\DN{\limsupz}[1]{\limsup_{#1\to0}}
\DN{\limsupdz}[1]{\limsup_{#1\downarrow0}}
\DN{\liminfi}[1]{\liminf_{#1\to\infty}}
\DN{\liminfz}[1]{\liminf_{#1\to0}}
\DN{\liminfdz}[1]{\liminf_{#1\downarrow0}}
\DN{\supnor}[1]{\| #1\|_{\infty}}
\DN{\sumii}[1]{\sum_{#1=1}^{\infty}}
\DN{\sumi}[1]{\sum_{#1=0}^{\infty}}

\DN\RA{\Rightarrow} \DN\LA{\Leftarrow}

\DN\PFEnd{\qed \smallskip}
\DN\PF{\begin{proof}} \DN\PFEND{\qed\end{proof}}

\DN\cREFa{\cref{;53} ( \nN )}

\DN\muetaN{\mu _{\mathrm{gin}}^{\etaN }}
\DN\KN{K_{ \mathrm{gin} }^{\etaN }}
\DN\KKN{\mathsf{K}_{ \mathrm{gin}}^{\etaN }}
\DN\etaN{\eta _{\nN }}
\DN\etaNi{\eta _{\nN ,i}}

	\DN\Pmg{\PPrNN }\DN\Emg{\EErNN } \DN\Pmgone{\Pmg }

\DN\Xti{X_t^i}
\DN\Mti{M_t^i}
\DN\si{s_i}
\DN\sj{s_j}
\DN\sk{s_k}
\DN\Xzi{X_0^i }

\DN\XoneX{\X _{\rN }^{\nN , [1]}}
\DN\rhoNone{\rho _{\rN }^{\nN , 1}}

\DN\lA{\langle} \DN\rA{\rangle}
\DN\QG{\As{QG} with $ \{ \ORnum \}_{ \Rm \in \N } $ for each $ \nuN $}
\DN\QGg{\As{QG} with $ \{ \ORnum \}_{ \Rm \in \N } $} 
\DN\QGO{\As{QG} with $ \{ \ORnum \}_{ \Rm \in \N } $ for each $ \nuN $ and that 
$ \{ \ORnum \}_{ \Rnu , m \in \N } $ satisfies \eqref{:63q}--\eqref{:63p}}

\DN\sumN{\sum_{ j\neq i }^{ \nN}}
\DN\V{V}
\DN\Vb{ \V _{ \beta }}
\DN\lE{\E }
	
\DN\E{ \mathcal{E}} 

\DN\pP{P} \DN\qQ{Q}

\DN\mF{\mathfrak}

\DN\Nm{ \nN , m }

	\DN\p{p}	
\DN\rsp{r,s,\p }
\DN\kkk{k}
\DN\nnn{n}
\DN\nnnN{ \nnn _{ \nN }}

\DN\dia{\diamondsuit}
\DN\HR{\mathcal{H}_{ \rR } }
\DN\HRR{\mathcal{H}_{ \rR -1 }}
\DN\HQd{\mathcal{I}_{ \qQ ,\rR }}
\DN\mM{M}
\DN\qN{ q _{ \nN }}
\DN\qM{ q _{ \mM }}
\DN\rN{ r _{ \nN }}

	\DN\chika{ \chi _{ \kR }} 

\DN\RN{\R ^{\nN }}
\DN\RNN{\R ^{\N }}
\DN\RI{\rR , \infty }
	\DN\Rm{\rR , m }
	\DN\Ql{\qQ , l }
	\DN\Inu{\infty , \nu }
	\DN\Rnu{ \rR , \3 }
	\DN\RRnu{ \rR +1 , \3 }
	\DN\Qnu{ \qQ , \3 } 

\DN\psiR{\psi _{ \Rnu }}
\DN\varphiQ{\varphi _{\Qnu }}
\DN\varphiR{\varphi _{\Rnu }}

\DN\nuN{\nu \in \N }
\DN\ZVbN{\mathcal{Z}_{\V ,\beta }^{ \nN }}
\DN\dRnu{ \dom _{ \Rnu }}

	\DN\kR{ \kappa , \rR }
	\DN\kRR{ \kappa , \rR , \rR ' }
	\DN\kRN{ \nN , \kappa , \rR , \rR ' }

\DN\akR{ a _{ \kR }}
\DN\akRR{ a _{ \kR + 1 }}
\DN\akkR{ a _{\kappa + 1 , \rR }}

	\DN\akrN{ a _{ \rN , \rN }}

\DN\0{ \lE _{ \rN }^{ \nN }} 

\DN\1{\dom _{\circ \3 }^{\mu }}

\DN\SrNSS{\SO _{\rN }^m \ts \sSS _{\rN }}
\DN\2{L ^2 ( \SrNSS , \murNNm ) }
 \DN\3{ \nu } 
\DN\4{( \E _{ \rN }^{ \nN } , \widehat{\dom }_{ \rN }^{ \nN })}

\DN\6{\eqref{:63q}--\eqref{:63p}}

	\DN\7{ \chi _N } 
 \DN\fqN{ \f _{ \qN }} \DN\fqM{ \f _{ \qM }}

\DN\XXX{Assume that $ \{ \fN \}_{\nN \in \N } $ with $ \fN \in \LmrNN $ weakly converges to $ \f \in \Lm $ in the sense of \dref{d:42}. }

\DN\g{g} \DN\f{f}

\DN\fH{\widehat{\f }}

\DN\fN{ \f _{ \nN } }
\DN\fNW{ \widehat{ \f }_{ \nN } }
\DN\fW{ \widehat{ \f } } 	\DN\fn{ \f _{ \nnn }} 

\DN\gN{ \g _{ \nN }}
\DN\gNnu{ \g _{ \nN , \nu }} 

\DN\y{\mathbf{y}} \DN\x{\mathbf{x}} \DN\z{ \mathbf{z}} 
\DN\xm{ \x _m} \DN\xn{ \x _n} \DN\ym{ \y _m} \DN\yn{ \y _n} \DN\xN{\x _{\nN }}
\DN\sN{\mathbf{s}_{\nN }}

\DN\sSS{ \mF{S}}
\DN\sss{ \mF{s}}
\DN\ww{ \mF{w}}
\DN\PP{ \mF{P}}	 \DN\EE{ \mF{E}}

\DN\PPxs{\widetilde{\underline{\PP }}_{ \rN , \xs }^{\nN , [m]}}
\DN\PPxsz{\widetilde{\underline{\PP }}_{ \rN , \sss }^{\nN , [0]}}
\DN\xs{( x , \sss )}
\DN\XX{ \mF{X}}

\DN\mrXX{\mr (\mathbf{X})}
\DN\mr{ \mF{m}_{\rR , T } }

\DN\da{d_{ \kR }(\sss )}
\DN\SSa{ \KK _{ \kR }}
\DN\SSaa{ \SSa ^+}
\DN\SN{ \sS ^{ \N }}

\DN\muR{ \mu _{ \rR }}
\DN\LmR{L^{2} (\muR )}
\DN\LmRC{L^2 (\CSR , \muR )}

\DN\Lm{L^{2}(\mu )}
\DN\LmN{L^{2}(\muN )}
\DN\LmRN{L^{2}(\muRN )}

\DN\LmrNN{L^{2}(\murNN )}

\DN\Lploc{L_{ \mathrm{loc}}^{p}}
\DN\Lqloc{L_{ \mathrm{loc}}^{q}}
\DN\Loneloc{L_{ \mathrm{loc}}^{1}}
\DN\Ltwoloc{L_{ \mathrm{loc}}^{2}}

\DN\SSS{ \sS \ts \sSS }

\DN\Sm{ \sS ^{ m }}
\DN\SR{ \sS _{ \rR }} 
\DN\SQ{\sS _{\qQ }}
\DN\SRR{ \sS _{ \rR -1}} 
\DN\SQQ{ \sS _{ \qQ -1}} 
\DN\SrN{ \sS _{ \rN }} 

\DN\TQl{\TQ ^l}
\DN\TQ{T _{\qQ }}
\DN\TQQ{T _{\qQ ' }}

\DN\SO{\overline{\sS }}
\DN\SOrN{\SO _{\rN }}
\DN\SOR{\SO _{\rR }}

\DN\SRm{ \SR ^{ m }}
\DN\SrNm{\sS _{ \rN }^m } 

\DN\ORm{ O _{\rR }^m }
\DN\ORnum{ O _{ \Rnu } ^{m}} 

\DN\SRbar{\SO _{\rR }}
\DN\Cf{\mathrm{Cf}}
\DN\CSRbar{\Cf (\SOR )}
\DN\CSrNbar{\Cf (\SO _{\rN } )}
\DN\CSR{ \Cf (\SR ) }
\DN\CSRm{ \Cf ^m (\SR ) } 

\DN\OO{\mathfrak{O}} 
\DN\OOR{\OO _{\rR }}
\DN\OORm{\OOR ^{ m }}

 \DN\OORnu{ \OO _{ \Rnu }}		
 \DN\OORRnu{ \OO _{ \RRnu }} 	

\DN\SSR{ \sSS _{ \rR }}
\DN\SSRr{ \sSS _{ \rR +1}}
\DN\SSq{ \sSS _{ \qQ }}
\DN\SSRl{ \SSR ^l} 
\DN\SSRm{ \SSR ^m} 
\DN\SSRrm{ \SSRr ^m} 
\DN\SSqm{ \SSq ^m} 
\DN\SSRn{ \SSR ^n} 

\DN\KK{\mathfrak{K}}
\DN\KKkR{ \KK _{ \kR } } 
\DN\KKk{\KK _{ \kappa }} 

\DN\KKrN{ \KK _{ \rN , \rN }}

\DN\OORnunu{ \OO _{ \Rnu +1 }}

\DN\OORnum { \OORnu ^{ m }} 
\DN\OOnu{ \OO _{ \nu }} 		
\DN\OOnuR{\OOnu ^{\rR }}

\DN\Orm{ \SRm } 
\DN\OOrm{ \SSRm }

\DN\MRnu{\mathfrak{M}_{\Rnu }}
\DN\NR{\mathfrak{N}_{\rR }}

\DN\SSs{ \sSS _{ \mathrm{s}}}
\DN\SSi{ \sSS _{ \mathrm{i}}}
\DN\SSsi{ \sSS _{ \mathrm{si}}}

\DN\LambdaR{ \Lambda _{ \rR }}
\DN\LambdaRm{ \LambdaR ^m }

\DN\WSs{ W (\SSs )}
\DN\WSS{ W (\sSS )}
\DN\WSsi{ W (\SSsi )}
\DN\WSsiNE{ W _{ \mathrm{NE}} (\SSsi ) }
\DN\WSsNE{ W _{ \mathrm{NE}} (\SSs ) }

\DN\pirN{\pi _{\SOrN }}
\DN\piR{ \pi_{ \rR } }
\DN\piQ{ \pi_{ \qQ } }
\DN\piRR{\pi _{\rR - 1}}
\DN\piRc{ \piR ^c}
\DN\piRd{ \piR ^{\dia }}

\DN\AQl{\AQ ^{ l_{ \qQ }}}
\DN\ARi{\mathbf{A}_{\rR }^{ i }}
\DN\ARl{\mathbf{A}_{\rR }^{ l }}
\DN\ARm{\mathbf{A}_{\rR }^m}
\DN\ARn{\mathbf{A}_{\rR }^n}
\DN\AR{\mathbf{A}_{\rR }} 
\DN\ARe{\mathbf{A}_{\rR ,\varepsilon }} 
\DN\AQ{\mathbf{A}_{\qQ ,\rR }} 
\DN\AQe{N (\Ql )} 
\DN\AQee{N (\qQ , l_{\qQ } )} 
\DN\BR{\mathbf{B}_{\rR }} 
\DN\BRm{\BR ^m}

\DN\As[1]{\textbf{(#1)}}

\DN\Rs{\rR , \sss } 
\DN\muk{\mu _{(k)}}
\DN\mukk{\mu _{(k+1)}}
\DN\mukRs{\mu _{(k) ,\Rs }}
\DN\mukR{\mu _{(k) ,\rR }}
\DN\mukkR{\mu _{(k+1), \rR }}

\DN\muRsone{ \muRs ^{[1]} } 
\DN\muRs{ \mu _{ \Rs }} 

\DN\muone{ \mu ^{[1]}} \DN\mum{ \mu ^{[m]}}

\DN\muN{ \mu ^{ \nN }}
\DN\muRN{ \muR ^{\nN }} 

\DN\murNN{\mu _{\rN }^{\nN }}
\DN\murNNm{\mu _{\rN }^{\nN , [ m ]} }
\DN\murNNone{ \mu _{ \rN }^{\nN , [1]} }

\DN\muNc{ \check{ \mu }^{ \nN }}

\DN\mul{ \mu \circ \lab ^{-1}}

\DN\mut{ \mu _{ \theta }}
\DN\mutx{ \mu _{ \theta ,x}}
\DN\muNtx{ \mu _{ \theta ,x}^{ \nN }}

\DN\muNt{ \mu _\theta ^{ \nN }}
\DN\muNtc{ \check{ \mu }_\theta ^{ \nN }}

\DN\muVN{ \mu _{V , \beta } ^{ \nN }}
\DN\nuVN{ \nu _{V , \beta } ^{ \nN }}
\DN\nuVNtwo{ \nu _{V , 2 } ^{ \nN }}

\DN\mVtN{ \mathbf{m}_{\V , \beta ,\theta }^{ \nN }}
\DN\muVtN{ \mu _{\V , \beta ,\theta }^{ \nN }}
\DN\muVNC{\check{ \mu }_{\V , \beta , \theta }^{ \nN }}

\DN\dom{ \mathcal{D}}
\DN\dmu{ \mathrm{d}^\mu}
\DN\dmuN{ \mathrm{d}^{ \nN }} 

\DN\di{ \dom _{ \circ}}
\DN\dimu{ \dom _{ \circ}^{ \mu }}
\DN\dik{ \di ^{k}}
\DN\domK{ \dom _{ \Inu }}
\DN\domKN{ \dom _{ \Inu }^{\nN }}
\DN\domI{ \dom _{ \infty }}

\DN\dN{ \dom ^{ \nN }}
\DN\dR{ \dom _{ \rR }}
\DN\dQ{ \dom _{ \qQ }}
\DN\dRk{ \dom _{ \Rnu }}
\DN\dRkn{ \dom _{ \Rnu ,\nnn }}

\DN\dimuN{ \dom _{ \circ}^{ \nN }}
\DN\dimuNm{ \dom _{ \circ}^{ \nN ,m}}

\DN\ld{ \underline{ \dom }}
\DN\ldR{ \ld _{ \rR }}
\DN\ldRn{ \ld _{ \rR , \nnn }}
\DN\ldRnN{ \ld _{ \rR , \nnn }^{ \nN }}

\DN\ldrNNm{ \ldW _{ \rN }^{ \nN , [ m ] }}

\DN\ldW{ \widetilde{ \ld }}
\DN\ldRnNW{ \ldW _{ \rR , \nnn }^{ \nN }}

\DN\ldRNW{ \ldW _{ \rR }^{ \nN }}
\DN\ldRnuNW{ \ldW _{ \Rnu }^{ \nN }}
\DN\ldrNNW{ \ldW _{ \rN }^{ \nN }}

\DN\ldRN{ \ldR ^{ \nN }}

\DN\ldRnu{ \ld _{ \Rnu }}
\DN\ldWnu{ \ld _{ \Inu }} 
\DN\ldRnuW{ \ld _{ \Rnu }} 	

\DN\ldRnun{ \ld _{ \Rnu , \nnn }}
\DN\ldRnunN{ \ldRnun ^{ \nN }}

\DN\enu{\varepsilon _{\nu }}
\DN\enunu{\varepsilon _{\nu +1 }}

\DN\EN{ \E ^{ \nN }} \DN\ER{ \E _{ \rR }}

\DN\lER{ \lE _{ \rR }}
\DN\lEInu{ \lE _{ \Inu }} 
\DN\Enu{ \E _{ \Inu }}
\DN\EnuN{ \Enu ^{ \nN }}
\DN\ERm{ \ER ^m} 

\DN\Edi{ \E _{ \infty }}

\DN\ERnu{ \E _{ \Rnu }}
	\DN\ERkK{ \E _{ \Rnu +1 }}
	\DN\ERkm{ \E _{ \Rnu }^{ m }}

\DN\lErN{ \lE _{ \rN }}
\DN\lErNone{ \lE _{\rN ,(1)}} 

\DN\lErNN{ \lE _{ \rN }^{ \nN }}
\DN\lErNNone{ \lErNone ^{ \nN }}

\DN\lERnu{ \lE _{ \Rnu } }
\DN\lERnuN{ \lERnu ^{ \nN } }
\DN\lErNnu{ \lE _{ \rN , \nu }}
\DN\lErNnuN{ \lErNnu ^{ \nN }}

\DN\ERN{ \ER ^{ \nN }}
\DN\ERNm{ \ER ^{ \nN ,m}}

\DN\lERN{ \lER ^{ \nN } }
\DN\lEp{ \lE _{ (\p )}}
\DN\lERNW{ \widetilde{ \lE }_{ \rR }^{ \nN }}

\DN\Eone{\E _{(1)}}
\DN\ERone{\E _{\rR , (1)}}

\DN\pN{ \p \in \N }

\DN\drN{ \dom _{ \rR }^{ \nN }}
\DN\domi{ \dom _{ \infty }}
\DN\ldi{ \ld_{ \infty }}
\DN\lcp[1]{(#1 )_{ \mathrm{reg}}}

\DN\ulab{ \mathfrak{u} }
\DN\lab{ \mathfrak{l} }
\DN\labj{\lab _j }
\DN\lpath{ \lab _{ \mathrm{path}}}
\DN\lkpath{ \lab _{k,\mathrm{path}}}
\DN\upath{ \ulab _{ \mathrm{path}}}
\DN\lpathN{ \lpath ^{ \nN }}
\DN\labN{ \lab ^{ \nN }} 
\DN\labNm{ \lab ^{ \nN ,m}} 
	\DN\lpathNi{ \lab _{\mathrm{path}, i } ^{ \nN }}
\DN\PPs{ \PP _{ \sss }}
\DN\EEs{ \EE _{ \sss }}
\DN\Pmu{ \PP _{ \mu }}
\DN\PPrNN{ \underline{\widetilde{\PP }}_{ \rN }^{ \nN }}
\DN\PPrNNxi{ \underline{\widetilde{\PP }}_{ \rN , \xi ^{\nN }d\muN }^{ \nN }}

\DN\EErNN{ \underline{\widetilde{\EE }}_{ \rN }^{ \nN }}

\DN\PmuX{ \PP _{ \mu }^{ \XX (t)}}
\DN\Pxt{ \PP _{(x,\sss ) }}
\DN\PPk{ \PP ^{k}}
\DN\Pmt{ \PP _{ \sss }}

\DN\rNone{ \rho ^{N,1}}
\DN\rNtwo{ \rho ^{N,2}}
\DN\rNm{ \rho ^{ \Nm }}
\DN\rNnk{ \rho ^{N,n+k}}

\DN\rhom{\rho ^m }

\DN\Capa{ \mathrm{Cap}}

\DN\B{\mathbf{B}} \DN\X{ \mathbf{X} }
\DN\XB{( \XXrNN , \PPrNN )} 

\DN\XN{ \X ^{ \nN }}
\DN\XXN{ \XX ^{ \nN }}
\DN\XrN{ \X _{ \rN }}
\DN\XrNN{\X _{ \rN }^{ \nN }}
\DN\XrNNm{\XrN ^{\nN ,[ m ]} }
\DN\XRN{ \X _{ \rR }}
\DN\XRNN{\X _{ \rR }^{ \nN }}

\DN\XXrNN{\XX _{ \rN }^{ \nN }}
\DN\XXrNNms{\XX ^{ m *}}

\DN\XXRN{\XX _{ \rR }^{ \nN }}
\DN\XXtt{ \X _{ \mathbf{t}}}
\DN\Xm{ \X ^{m}}
\DN\Xms{ \X ^{m*}}
\DN\XNm{ \X ^{ \Nm }}
\DN\XXtid{ \XX ^{i\dia } (t)}
\DN\XXuid{ \XX ^{i\dia } (u)}
\DN\XXNuid{ \XX _{\rN }^{N,i\dia } (u) }
\DN\XXNtid{ \XX _{\rN }^{N,i\dia } (t) }

\DN\sigmaRm{ \sigma _{\rR }^m }
\DN\sigmaRnm{ \sigma _{\rR ,n}^m }
\DN\sigmaRmn{ \sigma _{ \rR }^{ \nN , m}}

\DN\DDD{ \mathbb{D}}\DN\DDDa{ \DDD ^{a}}
\DN\DDDaR{ \DDDa _{ \rR }}
\DN\DDDR{ \DDDa _{ \rR }}
\DN\DDDRm{ \DDD _{ \rR } ^{a,m}} 
\DN\DDDam{ \DDD ^{a, m } }

\DN\Br{\mathcal{B}_{\rR } }
\DN\Brb{ \Br ^{b}}
\DN\Brnu{ \mathcal{B}_{ \Rnu }} 
\DN\Brnun{ \mathcal{B}_{ \Rnu }}
\DN\Brnunn{ \mathcal{B}_{ \Rnu + 1 }} 

\DN\muAVN{ \mu _{ \mathrm{Ai}, \beta ,V}^{ \nN }}
\DN\muA{ \mu _{ \mathrm{Ai},\beta } }
\DN\muAtwo{ \mu _{ \mathrm{Ai}, 2 }}

\DN\rhoAVNn{ \rho _{ \mathrm{Ai}, \beta , V}^{N,n}}
\DN\rhoAn{ \rho _{ \mathrm{Ai}, \beta }^n}
\DN\Ai{ \mathrm{Ai}}
\DN\KA{K_\mathrm{Ai}}
\DN\muG{ \mu _{ \mathrm{Gin}}}
\DN\muGN{ \mu _{ \mathrm{Gin}}^{ \nN }}
\DN\KG{K_{ \mathrm{Gin}}}
\DN\rhoGn{ \rho_{ \mathrm{Gin}}^n}
\DN\rhoGNn{ \rho_{ \mathrm{Gin}}^{N,n}}
\DN\muTrN{ \mu _{ \mathrm{Tr}^2}^{ \nN }}

\DN\rhoV{ \rho _V}
\DN\musin{ \mu _{ \sin ,2}}
\DN\musinbeta{ \mu _{ \sin , \beta }}
\DN\mVN{ \mathbf{m}_{\V , \beta }^{ \nN }}
\DN\PVN{ \pP _{\V , \beta }^{ \nN }}
\DN\Ksin{K_{ \sin }}
\DN\Ksinbeta{K_{ \sin ,\beta }}
\DN\rhosinm{ \rho _{ \sin ,2}^m}
\DN\rhosinmbeta{ \rho _{ \sin ,\beta }^m}
\DN\rhoVNtm{ \rho_{\V ,\beta ,\theta }^{ \Nm }}

\begin{document}

\title{\textsf{Dynamical Universality for Random Matrices}}

\author{Yosuke Kawamoto   \and      Hirofumi Osada \\
\\
\texttt{To appear in}\\
\text{Partial Differential Equations and Applications}}



\maketitle

\begin{abstract}
We establish an invariance principle corresponding to the universality of random matrices.
More precisely, we prove the dynamical universality of random matrices in the sense that, if the random point fields $ \muN $ of $ \nN $-particle systems describing the eigenvalues of random matrices or log-gases with general self-interaction potentials $ \V $ converge to some random point field $ \mu $, then the associated natural $ \muN $-reversible diffusions represented by solutions of stochastic differential equations (SDEs) converge to some $ \mu $-reversible diffusion given by the solution of an infinite-dimensional SDE (ISDE). 
Our results are general theorems that can be applied to various random point fields related to random matrices such as sine, Airy, Bessel, and Ginibre random point fields. 
In general, the representations of finite-dimensional SDEs describing $ \nN $-particle systems are very complicated. 
Nevertheless, the limit ISDE has a simple and universal representation that depends on a class of random matrices appearing in the bulk, and at the soft- and at hard-edge positions.
Thus, we prove that ISDEs such as the infinite-dimensional Dyson model and the Airy, Bessel, and Ginibre interacting Brownian motions are universal dynamical objects. 

\\
 \textbf{Subclass: }{MSC 60B20 \and MSC 60J40 \and MSC 60H10}
\end{abstract}
\noindent 
\footnote{This work was supported by JSPS KAKENHI Grant Numbers JP16H06338, \\\quad \quad JP20K20885, JP21H04432, and JP21K13812. }
\footnote{Y. Kawamoto, 
             Fukuoka Dental College, Fukuoka 814-0193, Japan \\\quad \quad 
              \texttt{kawamoto@college.fdcnet.ac.jp}       } 
\footnote{
           H. Osada, 
              Faculty of Mathematics, Kyushu University, Fukuoka 819-0395, Japan\\\quad \quad 
                            \texttt{osada@math.kyushu-u.ac.jp}     
}
 { \small \tableofcontents }

\section{Introduction}\label{s:1} 
The concept of universality in strongly correlated systems was envisioned by Wigner, who conjectured that the eigenvalue distribution of a large random matrix is universal, that is, the eigenvalue distribution depends only on the symmetry classes of matrices, rather than on the distributions of matrix components.
The universality of random matrices is a central concept in random matrix theory, and has been studied intensively over the past two decades 
(see, e.g.\,\cite{bey.14-duke,dkmvz,dkmvz.2,d-g,d-g2,J01,LL16,L09,lsy,TV.15}). 
However, its dynamical counterpart has been much less studied. In the present paper, we establish the concept of dynamical universality for random matrices in a general framework. 

Let us recall some universality results for the sine$_\beta $ random point field derived by Deift {\em et al.} \cite{dkmvz,dkmvz.2} following Deift and Gioev \cite{d-g2}. 
We consider the ensembles $ \mathcal{M}^{ \nN } $ with the distribution 
\begin{align} &\notag 
 \mathcal{P} _{\beta }^{ \nN } ( d \mathcal{M}^{ \nN } ) = \frac{1}{\ZVbN } 
 e^ {- \nN \mathrm{tr} \V ( \mathcal{M}^{ \nN } )} d \mathcal{M}^{ \nN } 
\end{align}
for $ \beta = 1,2$, and $ 4 $. 
Here, the ensembles $ \mathcal{M}^{ \nN } $ consist of $ \nN \ts \nN $ real symmetric matrices, 
$ \nN \ts \nN $ Hermitian matrices, and 
$ 2\nN \ts 2\nN $ Hermitian self-dual matrices for $ \beta = 1,2$, and $ 4 $, respectively. 
For $ \beta = 2 $, the potential $ \V $ is a real analytic function satisfying 
\begin{align}\label{:11Q}
& \limi{|x|} \frac{ \V (x)}{ \log |x| } = \infty 
.\end{align}
For $ \beta = 1,4$, the potential $ \V $ is a real polynomial such that 
\begin{align}\label{:11q}&
\V (t) = \sum_{n=0}^{2m} v _{n} t^n , \quad v _{2m} > 0
.\end{align}
The density of the distribution of eigenvalues $ \xN $ of $ \mathcal{M}^{ \nN } $ is given by 
\begin{align}\label{:11s}
\PVN (\xN ) = \frac{1}{\ZVbN } 
\prod_{i<j}^{ \nN } | x_i - x_j |^{\beta } \prod_{ k =1}^{ \nN } e^{- \nN \Vb (x_k)}
,\end{align}
where $ \xN =(x_1,\ldots,x_{ \nN })\in \RN $, 
$\ZVbN $ is the normalizing constant, and we set 
\begin{align}\label{:11u}&
 \Vb (x) = 
 \begin{cases}
\V (x) , & \beta = 1,2 , \\
 2 \V (x) & \beta = 4 
.\end{cases}
\end{align}

Note that, according to the logarithmic interaction potential, particles repel each other. 
The logarithmic interaction potential has a strong long-range effect that causes special phenomena to occur in particle systems. 
One example is the convergence of the empirical distribution to a deterministic distribution that has non-degenerate density, typically known as Wigner's semi-circle law in the limit. 

Let $ \nuVN $ be the random point field such that its labeled density is given by $\PVN $. 
The behavior of the system as the number of particles $N$ tends to infinity has been extensively studied.
We set $ \mF{x} = \sum_{ i }\delta_{x_i}$, where $ \delta_a $ denotes the delta measure at $ a $.
Then, there exists a probability density function $ \varrho_{V}$ on $ \R$ such that 
\begin{align}\label{:11b}
\limi{N} \int \frac{1}{N}\mF{x} ((-\infty, {s}]) \nuVN ( d \mF{x})
=\int_{-\infty}^s\varrho_{V}(x)\,dx
.\end{align}
%
If $V(x)=x^2$, then $ \nuVNtwo $ gives the eigenvalue distribution of the Gaussian unitary ensemble (GUE), which is the Hermitian random matrix whose entries follow an independent and identically distributed Gaussian distribution. 
The probability measure $ \varrho_{V} d x $ is simply the Wigner semicircle law, 
which is given by 
$ \varrho_{V}( x ) = 
\frac{1}{ \pi }\sqrt{2- x ^2}\mathbf{1}_{ \{| x |<\sqrt{2}\}}$ 
 (see, e.g., \cite{AGZ,mehta}).

The convergence in \eqref{:11b} is in the macroscopic regime.
Next, we consider the microscopic scaling limit.
More precisely, we consider a local fluctuation of \eqref{:11s} and obtain a random point field with infinitely many particles as the limit.
Here, we take the bulk scaling limit.
For a constant $ \theta \in\R $ satisfying 
\begin{align}\label{:11c}
\varrho_{V}(\theta )>0
,\end{align}
we set the bulk scaling at $ \theta $ as 
\begin{align}\label{:11d}
 x \mapsto \frac{s}{N\varrho_{V}(\theta )} +\theta 
.\end{align}
Let $ \mVtN $ be the rescaled density function of $ \PVN $ in \eqref{:11s} 
under the scaling defined by \eqref{:11d}. Then, $ \mVtN $ is given by 
\begin{align}\label{:11e}
\mVtN (\mathbf{s }_{ \nN })=\frac{1}{\mathcal{Z} _{V,\beta ,\theta }^{ \nN }} 
\prod_{i<j}^{ \nN } | \si - \sj |^{ \beta }
\prod_{ k =1}^{ \nN } 
\exp \Big(-N V_{ \beta } \Big(\frac{s_{ k }}{N\varrho_{V}(\theta )} + 
\theta \Big) \Big) 
.\end{align}
We define $ \muVtN $ as the random point field whose labeled density is given by $ \mVtN $.

Let $ \musinbeta $ be the sine$ _\beta $ random point field. 
Let $ \rho _{\mathrm{sin},\beta }^m $ be the $ m $-point correlation function of $ \musinbeta $. 
If $ \beta = 2 $, then $ \musinbeta $ is the determinantal random point field whose kernel is given by 
\begin{align*}
\Ksin (x,y)=\frac{ \sin \pi (x-y)}{ \pi (x-y)}
.\end{align*}
Then, by definition, 
\begin{align} &\notag 
\rhosinm (x_1,\ldots,x_m )=\det [\Ksin (x_i,x_j)]_{1\le i, j\le m }
.\end{align}
Similar formulae are also known for $ \beta = 1,4 $ (see \cite{mehta}). 

Bulk universality for log-gases asserts that, for a suitable and wide class of $ \V $ and for any $ \theta $ satisfying \eqref{:11c}, 
\begin{align} &\notag 
\limi{N}\muVtN = \musinbeta 
\quad\text{ weakly,}
\end{align}
or, more strongly, for each $ m \in \N $ the $ m $-point correlation function $ \rhoVNtm $ of 
 $ \muVtN $ satisfies 
\begin{align} \label{:11j}
\limi{N} \rhoVNtm = \rhosinmbeta 
\end{align}
uniformly on each compact set.
Note that the limit $ \musinbeta $ is independent of $V$ and $ \theta $, and, in this sense, 
the sine$ _{\beta }$ random point field can be thought of as a universal object. We call such universality \emph{static} or \emph{geometric} because it involves no time evolution.

Consideration of the classical invariance principle yields a natural question: what is the dynamical counterpart of geometric universality?

We consider an $ \nN $-dimensional stochastic differential equation (SDE) of 
$ \X ^{ \nN }= (X^{ \nN , i})_{i=1}^{ \nN } $ 
corresponding to $ \muVtN $ such that, for $1\le i\le N$, 
\begin{align}\label{:12a}
X^{ \nN , i } (t) - X^{ \nN , i } (0) = 
B^i(t) & + \frac{\beta }{2} 
\int_0^t 
 \sumN \frac{1}{X^{ \nN , i } (u) - X^{ \nN , j } (u)} du
\\ \notag &
 - \frac{1}{2} 
 \int_0^t 
\frac{1}{\varrho_{V}(\theta )} 
\Vb ' \Big(\frac{X^{ \nN , i } (u) }{N\varrho_{V}(\theta )} + \theta \Big) du 
.\end{align}
Here, $ B^i $, $i=1,\ldots,N $, are independent standard Brownian motions with 
$ B^i(0) = 0 $. 
We derive \eqref{:12a} from $ \muVtN $ as follows. 
Let $ \muVNC $ be the distribution of the labeled particle system of $ \muVtN $. 
Then, 
\begin{align}&\label{:12b}
 \muVNC (d\xN )= \mVtN(\xN ) d\xN 
,\end{align}
and consider the Dirichlet form on $ L^2(\RN , \muVNC ) $ such that 
\begin{align} & \notag 
\mathcal{E}^{\muVtN } (f,g) = 
\int_{ \RN } \frac{1}{2} \sum_{i=1}^{N } \PD{f}{x_i} \PD{g}{x_i} 
 \mVtN(\xN ) d\xN 
.\end{align}
Integrating by parts and using \eqref{:12b}, we have that
\begin{align}\label{:12d}
\mathcal{E}^{\muVtN } (f,g) = & 
- \int_{ \RN } \half 
\sum_{i=1}^{ \nN } \Big\{ 
\frac{ \partial ^2 f}{(\partial x_i)^2} + \PD{ \log \mVtN }{x_i} \PD{f}{x_i} \Big\} g \, 
 \mVtN(\xN ) d\xN 
\\ \notag =& 
- \int_{ \RN } ( A^{N } f ) g \, \mVtN(\xN ) d\xN 
.\end{align}
Here, from \eqref{:11e} and \eqref{:12d}, we see that $ A^{N }$ is given by 
\begin{align} &\notag 
A^{N } = 
\half \Delta +
 \frac{\beta }{2} 
 \sum_{i=1}^{N } \sumN \frac{1}{x_i-x_j} \PD{}{x_i} 
 - \frac{ 1 }{2} 
\sum_{ k =1}^{N } 
\frac{1}{\varrho_{V}(\theta )} \Vb ' \Big(\frac{x_k }{N\varrho_{V}(\theta )} + 
\theta \Big) 
\PD{}{x_k } 
.\end{align}
Because $ \mVtN(\xN ) $ is bounded and continuous, we easily find that 
the bilinear form $ (\mathcal{E}^{\muVtN } , C_0^{ \infty}(\RN ) ) $ 
is closable on $ L^2(\RN , \muVNC ) $. 
Then, there exists a unique $ L^2$-Markovian semi-group $ \{T_t^{ \nN }\} $ 
on $ L^2(\RN , \muVNC ) $ associated with the closure of 
$ (\mathcal{E}^{\muVtN } , C_0^{ \infty}(\RN ) ) $. 
The semi-group $ \{ T_t^{ \nN }\} $ is given by the solution $ \X ^{ \nN } $ of \eqref{:12a}:
\begin{align*}&
T_t^{ \nN }f (\xN ) = E _{\xN } [f (\X ^{ \nN } (t))]
.\end{align*}
Here, $ E _{\xN }$ is the expectation with respect to the solution $ \X ^{ \nN }$ 
starting at $ \xN $. 
By construction, $ \X ^{ \nN } $ is reversible with respect to $ \muVNC $. 
We denote the unlabeled dynamics of $ \X ^{ \nN } $ as $ \XX ^{ \nN } $. 
Then, $ \XX ^{ \nN }(t) = \sum_{i=1}^{ \nN }\delta_{X^{ \nN , i }(t)}$ by definition. 
It is clear that $ \XX ^{ \nN }$ is reversible with respect to $ \muVtN $. 

If the initial conditions converge, it is interesting to determine whether the stochastic process 
$ \X ^{ \nN }$ converges or not, and to exploit the infinite-dimensional SDE (ISDE) satisfied by the limit stochastic process $ \X $. 

We set $ \X ^{ \Nm }=(X^{ \nN , i })_{i=1}^m$ and 
$ \X ^m = (X^i)_{i=1}^m$ for each $ m \in \N $. 
We shall prove that, loosely speaking, 
if the initial distributions of $ \X ^{ \Nm } $ converge to that of $ \X ^m $ in distribution, 
then 
\begin{align} &\notag 
\limi{N} \X ^{ \Nm } = \X ^m \quad \text{ weakly in $ C([0,\infty);\R^m) $}
.\end{align}
Here, $ \XN $ is a solution to \eqref{:16c} for $ \muN = \muVtN $. 
The limit ISDE is Dyson's model in infinite dimensions \cite{spohn.dyson}, 
which is an ISDE of $ \X = (X^i)_{i\in\N }$ such that, for $ i\in\N $, 
\begin{align}\label{:13b}
 X^{i} (t) - X^{i}(0) = 
 B^i (t) + \frac{\beta }{2} \int_0^t \limi{ \rR } 
 \sum_{ j\neq i \atop |X^i (u) - X^j (u) | < \rR } 
 \frac{1}{X^{i} (u) -X^{j}(u)} du 
.\end{align}
The associated unlabeled dynamics 
$ \XX (t) =\sum_{i=1}^{ \infty} \delta_{X^i (t)}$ are 
reversible with respect to the sine$ _{ \beta }$ random point field \cite{o.isde}. 
From the static universality defined in \eqref{:11j}, we expect the limit of \eqref{:12a} as $N\to\infty $ to be given by \eqref{:13b}.
In particular, the limit does not depend on $V$ and $ \theta $.
In other words, ISDE \eqref{:13b} is expected to be a dynamical universal object, 
which is a consequence of the present paper. 

The simplest case is $ \V (t) = t^2$, $ \beta = 2 $, and $ | \theta | < \sqrt{2}$. 
Because $ \varrho_{V}( \theta )=\frac{1}{ \pi }\sqrt{2-\theta^2} $, 
 \eqref{:12a} becomes 
\begin{align}\label{:13c}
X^{ \nN , i } (t) - X^{ \nN , i } (0) = 
B^i(t) & + 
\int_0^t \sum_{ j\neq i }^{ \nN} \frac{1}{X^{ \nN , i } (u) - X^{ \nN , j } (u) } du 
\\ \notag & - 
 \frac{\pi ^2}{N( 2 - \theta ^2 )} \int _0 ^t X^{ \nN , i } (u) du 
 - \frac{\pi \theta }{\sqrt{ 2 - \theta ^2}} t 
.\end{align}
In \cite{k-o.sdeg}, the authors proved that the solution of \eqref{:13c} converges to that of 
\eqref{:13b} in distribution in $ C([0,\infty);\R^m) $ for each $ m $ if 
the dynamics start from reversible measures. 
Note that the term $ - \frac{\pi \theta }{\sqrt{ 2 - \theta ^2} } t $ 
 in \eqref{:13c} disappears in \eqref{:13b}. 

Suppose that $ \V $ is a real analytic function satisfying \eqref{:11Q} and $ \beta = 2 $. 
Set 
\begin{align}\label{:13e}&
\cref{;13} := 
\half 
\frac{1}{\varrho_{V}(\theta )} \V ' (\theta )
.\end{align}
Then, as $ \nN \to \infty $, the remaining part in the last term of \eqref{:12a} is 
$ - \Ct t \label{;13} $, which vanishes in \eqref{:13b}. 
Thus, we see that the convergence has an SDE gap. This phenomenon is a result of the long-range interaction of the logarithmic potential, and never happens for translation-invariant random point fields with potentials of Ruelle's class.

Our result is an invariance principle in the sense that weak convergence occurs in the path space $ C([0,\infty); \R ^{\N } ) $. 
As an application of the invariance principle, we see 
from \eqref{:12a}, \eqref{:13e}, and \eqref{:16c} that, for $ a \ge 0 $ and $ i \le \nnnN $, 
\begin{align}\label{:13f} 
\limi{ \nN } \pP \Big( \max_{0 \le t \le T } 
 \Big\{
 X^{ \nN , i } (t) & - X^{ \nN , i }(0) - 
 \int_0^t \sum_{ j\neq i }^{ \nnnN} \frac{1}{X^{ \nN , i } (u) - X^{ \nN , j } (u) } du 
 \\ \notag &
 + \cref{;13} t 
 - 
 \int_{ \sS \backslash \sS _{ \rN } } \frac{1}{ X^{ \nN , i } (t) - y } \rNone ( y ) dy 
 \Big\} \ge a 
\Big)
\\ \notag &
= 2 \int_{ x \ge a }\frac{1}{\sqrt{2\pi T }} e^{- |x|^2 /2T} dx 
.\end{align}
Here, $ \rNone $ is the one-point correlation function of 
$ \muN := \mu _{\V , 2 ,\theta }^{ \nN } $ with respect to the Lebesgue measure. 
The process $ (X^{\nN , i})_{i=1}^{ \nnnN }$ is a solution of SDE \eqref{:16c} 
on $ \{ |x| < \rN \} ^{ \nnnN }$ for $ \muN $ with $ \V $ as above. 
The number of particles in $ \sS _{\rN } $ is denoted by $ \nnnN $, and so 
$ X^{\nN , i }(0) \in \sS _{\rN }$ for $ i=1,\ldots,\nnnN $. 
The radius of the domain $ \sS _{\rN }$ satisfies $ \limi{\nN } \rN = \infty $. 
Furthermore, \eqref{:16c} is given by \eqref{:12a} with the reflecting boundary condition and 
the free potential in \eqref{:16e} for $ \Psi ^{ \nN } (x) = - 2 \log |x| $. 
We prove that \eqref{:13f} holds in \sref{s:91}. 

The next example is the Ginibre random point field. 
The Ginibre random point field $ \mu _{ \mathrm{gin}}$ is a determinantal random point field 
on $ \R ^2 $ that has the kernel $ K_{ \mathrm{gin}} $ with respect to the Lebesgue measure such that
\begin{align}&\label{:14a}
K_{ \mathrm{gin}} (x,y) = \frac{1}{ \pi} 
\exp \Big(-\frac{|x|^2+|y|^2}{2} + x \bar{y}\Big)
.\end{align}
Here, we naturally regard $ \R ^2 $ as $ \C $ by $ (x,y)\mapsto x + \sqrt{-1}y$. 
Then, by definition, the $ m $-point correlation function $ \rho _{ \mathrm{gin}}^m $ 
with respect to the Lebesgue measure is given by 
\begin{align} &\notag 
\rho _{ \mathrm{gin}}^m (x_1,\ldots,x_{ m }) = 
\det [K_{ \mathrm{gin}}(x_i,x_j) ]_{i.j=1}^m 
.\end{align}
The Ginibre random point field $ \mu _{ \mathrm{gin}} $ 
is a limit of the random point fields $ \mu _{ \mathrm{gin}}^{ \nN } $ 
arising from the ensemble of the complex non-Hermitian Gaussian random matrices 
$ M_{ \mathrm{gin}}^{ \nN } = ( M_{ \mathrm{gin}}^{ \nN }({i,j}) )_{i.j=1}^{ \nN } $ 
in which the $ 2\nN ^2$ parameters are independent Gaussian random variables 
with mean zero and variance $ 1/2$. 
The labeled density $ \mathbf{m}_{\mathrm{gin}}^{ \nN } $ 
of the distribution of the eigenvalues of $ M_{ \mathrm{gin}}^{ \nN }$ is given by 
\begin{align}\label{:14c}&
 \mathbf{m}_{\mathrm{gin}}^{ \nN }( x_1 ,\ldots, x_{ \nN }) = 
 \frac{1}{\mathcal{Z}} 
 \{ \prod_{ i < j }^{ \nN } | x_i - x_j | ^2 \} 
 e^{ - \sum_{ k = 1}^{ \nN } | x_k |^2 }
.\end{align}
From \eqref{:14c}, similar to the case for \eqref{:12a}--\eqref{:12d}, we derive the SDE describing the labeled $ \nN $-particle dynamics $ \mathbf{X}^{ \nN } = (X^{ \nN , i})_{i=1,\ldots,\nN}$ 
as follows. For $ i = 1,\ldots,\nN $, 
\begin{align}\label{:14f}&
X^{ \nN , i } (t) - X^{ \nN , i } (0) = 
 B^i (t) - \int _0^t X^{ \nN , i } (u) du + 
 \int_0^t \sumN 
 \frac{X^{ \nN , i } (u) - X^{ \nN , j } (u)} { | X^{ \nN , i } (u) - X^{ \nN , j } (u) | ^2 } du 
.\end{align}
It has been shown \cite{k-o.fpa} that, under suitable assumptions regarding the initial distributions, the dynamics of the first $ m $ particles 
$ \mathbf{X}^{ \nN , m } = ( X^{ \nN , 1},\ldots,X^{ \nN , m })$ 
converge weakly in $ C([0,\infty); (\R ^2)^m)$ 
to those of the unique strong solutions $ \mathbf{X} = (X^i)_{ i=1}^{\infty}$ of the ISDE such that 
\begin{align}\label{:14h}& 
 X^i (t) - X^i (0) = B^i (t) - \int _0^t X^i (u) du + 
 \int_0^t \limi{\rR} 
 \sum_{ j \ne i \atop | X^j (u) | < \rR } 
 \frac{X^i (u) - X^j (u) } { | X^i (u) - X^j (u) | ^2 } du 
.\end{align}
It has also been shown \cite{o.isde,o-t.tail} that the solution $ \X $ of \eqref{:14h} satisfies the second ISDE
\begin{align} \label{:14i} &
X^i (t) - X^i (0) = B^i (t) + \int_0^t \limi{\rR} 
 \sum_{ j \ne i \atop |X^i (u) - X^j (u) | < \rR } 
 \frac{X^i (u) - X^j (u) } { | X^i (u) - X^j (u) | ^2 } du 
.\end{align}
Thus, both ISDEs have the same unique strong solution $ \X $, and the solution 
$ \X $ satisfies two different ISDEs \cite{o.isde,o-t.tail}. 
Such multiple representations of ISDEs are the result of the long-range nature of the logarithmic interaction potential in the drift terms of \eqref{:14h} and \eqref{:14i}. 

The unlabeled dynamics $ \XX (t) = \sum_{ i= 1 }^{\infty} \delta _{ X^i (t)}$ describe $ \mu _{ \mathrm{gin}} $-reversible diffusion \cite{o.rm}. 
The ISDEs are obtained from the general theory for random point fields and the associated ISDEs in \cite{o.isde,o-t.tail}. 
Note that, even in such a typical case, the convergence of finite particle systems is a sensitive problem, as we see from \eqref{:14f}--\eqref{:14i}. 

Akemann-Cikovic-Venker \cite{acv} proved the following universality of the Ginibre random point field. 
Let $ \mathcal{J}(\nN ) $ be the space of the normal matrices of order $ \nN $. 
For constants $ \gamma \ge 0$, $K_p\in\R $, and $ \omega \in[0,1) $, 
consider the probability measure on $ \mathcal{J}(\nN ) $ whose density is given by 
\begin{align} &\notag 
\sigma (J) = \frac{1}{ \mathcal{Z} }
\exp \Big\{
\frac{ \nN }{1-\omega ^2} \mathrm{Tr}(J J^* - \frac{ \omega }{2} (J^2+J^{*2})) - 
\gamma (\mathrm{Tr}JJ^*-\nN K_p)^2 
\Big\}.
\end{align}
Then, the joint density of the eigenvalues is proportional to 
\begin{align} \label{:14l} &
\{ \prod_{ i < j }^{ \nN }|z_i -z_j|^2 \}
 \\& \notag 
\ts 
\exp 
\Big\{ -\frac{ \nN }{1-\omega ^2} \Big( \sum_{i=1}^{ \nN }|z_i|^2-\frac{ \omega }{2} 
\sum_{i=1}^{ \nN }(z_i^2 +\bar{z_i}^2 ) \Big) 
 -\gamma 
 \Big(\sum_{i=1}^{ \nN }|z_i|^2- \nN K_p \Big)^2 
\Big\} 
\end{align}
Let $ \rNm _{ \mathrm{gin}} $ be the $ m $-point correlation function of the eigenvalue density corresponding to \eqref{:14l}. 
For positive constants $ \cref{;8a}$ and $ \cref{;8b}$, we set 
\begin{align} &\notag 
E =\{ z \in \C ; \cref{;8a} (\Re z)^2 + \cref{;8b} (\Im z)^2 < 1 \} 
.\end{align}
We quote a universality result for the Ginibre random point field from \cite{acv}. 
\begin{proposition}[{\cite[Theorem 1]{acv}}]	\label{l:11} 
There exist positive constants $ \Ct \label{;8a} $, $ \Ct \label{;8b} $, $ \Ct \label{;8c} $ 
depending on $ K_p $, $ \gamma $, $ \omega $ such that, for any 
$ \zeta \in \C \setminus \partial E $, 
\begin{align}& \notag 
\limi{ \nN }\frac{1}{ \nN } \rNone _{ \mathrm{gin}} (\zeta ) = 
\frac{ \cref{;8c}}{ \pi } 1_{E} (\zeta ) 
.\end{align}
Furthermore, for each $ m \in \nN $ and $ \zeta \in E $, 
\begin{align}& \notag 
\frac{1}{(\cref{;8c}\nN )^m } \rNm_{ \mathrm{gin}} 
\Big(\zeta + 
\frac{z_1}{ \sqrt{ \cref{;8c}\nN }},\ldots, \zeta + 
\frac{z_{ m }}{ \sqrt{ \cref{;8c}\nN }}\Big) 
= 
 \rhom _{ \mathrm{gin}} (z_1,\ldots,z_{ m }) + 
 \mathcal{O} \Big(\frac{1}{ \sqrt{ \nN }}\Big) 
.\end{align}
Here, the error term $ \mathcal{O} \Big(\frac{1}{ \sqrt{ \nN }}\Big) $ can be taken uniformly 
on each compact set in $ \mathbb{C}^m$. 
\end{proposition}

We shall investigate a dynamical counterpart of this result. 

From the same calculation as for \eqref{:12b}--\eqref{:12d}, 
we obtain from \eqref{:14l} the SDE of the $ \nN $ particles 
$ \mathbf{X}^{ \nN } = (X^{\nN , i } )_{i=1}^{ \nN }$ 
such that, for $ i=1,\ldots,\nN $, 
\begin{align}& \label{:15g}
X^{\nN , i } (t) - X^{\nN , i } (0) = B^i (t) + 
 \int_0^t 
\sum_{j\neq i}^{ \nN } 
\frac{ X^{\nN , i } (u)-X^{\nN , j }(u) }
{|X^{\nN , i } (u)-X^{\nN , j } (u) |^2} du 
\\ \notag &
- 
\int_0^t 
 \frac{ \nN }{1-\omega ^2} \frac{1}{ \sqrt{ \cref{;8c} \nN }}
 \Big(\zeta +\frac{X^{\nN , i } (u)}{ \sqrt{ \cref{;8c} \nN }}\Big) 
du 
\\ \notag &
+
\int_0^t 
 \frac{ \nN }{1-\omega ^2}\frac{\omega }{2}\frac{1}{ \sqrt{ \cref{;8c} \nN }}
\Big\{ 
 \Big(\zeta +\frac{X^{\nN , i } (u)}{ \sqrt{ \cref{;8c} \nN }}\Big) 
 + 
 \Big(\zeta + 
\frac{X^{\nN , i } (u)}{ \sqrt{ \cref{;8c} \nN }}\Big)^{ \dagger }\Big\} 
du 
\\ \notag & 
- \int_0^t 
\frac{2 \gamma}{ \sqrt{ \cref{;8c} \nN }} 
\Big(\zeta +\frac{X^{\nN , i } (u)}{ \sqrt{ \cref{;8c} \nN }}\Big)
\Big\{ 
\sum_{ j =1}^{ \nN }\Big|\zeta +\frac{X^{\nN , j }(u)}{ \sqrt{ \cref{;8c} \nN }}
\Big|^2- \nN K_{p} \Big\} du 
.\end{align}
Here, we set $ (x,y)^{ \dagger } = (x,-y) \in \R ^2$. 
The limit ISDE corresponding to $ \muG $ is 
\begin{align*}
X^{i} (t) - X^i (0) = B^i (t) + 
\int_0^t \limi{ \rR } \sum_{ j \ne i \atop |X^i(u)-X^j(u) |< \rR } 
\frac{X^i(u)-X^j(u) }{|X^i(u)-X^j(u) |^2}du 
\quad (i \in \N ).
\end{align*}
Thus, although the representation of the $ \nN $-particle SDE is quite complicated, the limit ISDE is very simple and universal.

To prove such dynamical finite particle approximations, the authors have previously established a general theory \cite{k-o.fpa}.
The framework in \cite{k-o.fpa} does not depend on the dimension of the underlying space, the inverse temperature, or the integrable structures; thus, the theory can be applied to many examples.
A key point in our previous paper is the control of drift terms in finite-dimensional SDEs, which provide a sensitive estimate for the long-range interaction potential. 
Actually, we have proved the dynamical bulk scaling limit by completing such an estimate \cite{k-o.sdeg}. %
However, when the potentials become some general $V(x) $, the calculations are more difficult.
In particular, when we consider an ISDE related to the Airy random point field, which arises from the soft-edge scaling limit of eigenvalue distributions of random matrices, this presents a more complicated problem. 
Indeed, the drift term in the corresponding finite-dimensional SDE 
includes a divergent term \cite{o-t.airy}. 

To overcome this difficulty, we construct a new method in the present paper. 
This approach uses the convergence concept of Dirichlet forms associated with finite or infinite particle systems. 
Let $ \rho ^{\nN , m }$ and $ \rho ^m $ be the $ m $-point correlation functions of $ \muN $ and $ \mu $, respectively. 
Let $ \sigma _{\rR }^m $ be the $ m $-point density function of $ \mu $ on $ \SR = \{ |s| < \rR \} $. 
As well as the existence of infinite particle dynamics, we assume two main conditions for convergence: 
\begin{itemize}
 \item[(I)] 
We assume that $ \limi{\nN } \rho ^{\nN , m } = \rho ^m $ uniformly on each compact set, and 
 the capacity of the zero points of $ \sigma _{\rR }^m $ vanishes. 

\item[(II)] 
We assume the uniqueness in law of weak solutions to the limit ISDE. 
\end{itemize}
\noindent 
Condition \thetag{I} is related to the static property of random point fields, while 
Condition \thetag{II} concerns the dynamical property of infinite-particle systems in the limit. 
These two ingredients are sufficient for dynamical universality. 

Only condition \thetag{I} is related to finite-particle systems. 
This is purely because of the static property of such systems---we do not require any assumption regarding estimates related to the dynamics of finite-particle systems, such as the estimates of drift terms in \cite{k-o.fpa}.
Additionally, the method derived in the present paper is independent of the dimension of an underlying space, the inverse temperature, or the integrable structures, as in \cite{k-o.fpa}.

Consequently, if the limit ISDE has a unique weak solution, the strong convergence of random point fields automatically implies dynamical convergence. 
It has been proved that several ISDEs with coefficients given by the logarithmic interaction potential have unique solutions \cite{o-t.tail,tsai.14}.
Therefore, the static strong universality of random matrices can be strengthened to dynamical universality, not only for Dyson's Brownian motion and the Ginibre random point field, but also for ISDEs related to the Airy random point field, the Bessel random point field, and so on.

\smallskip 

We now explain some of our main results. 
There exists a natural correspondence among a random point field $ \mu $ on $ \sS $, an unlabeled diffusion $ \XX $ in the configuration space $ \sSS $ over $ \sS$, 
and an ISDE $ \X $ on $ \SN $ (see \cite{o.isde,o.rm}). 
For a given random point field $ \mu $, the associated unlabeled diffusion $ \XX $ is given by distorted Brownian motion, which is a Dirichlet form whose energy and time change measures are common (see \sref{s:2}). 
The correspondence between $ \XX $ and $ \X $ is given by 
$ \XX (t) = \sum_{i=1}^{\infty} \delta_{ X^i (t)}$ and $ \X = (X^i)_{i\in \N }$. 
Furthermore, the labeled dynamics $ \X $ have a representation as a solution of an ISDE. 
The ISDE is described in terms of the logarithmic derivative $ \dmu $ of $ \mu $ as follows: 
\begin{align}\label{:16a}&
X^i(t) -X^i(0) = 
B^i(t) + \half \int_0^t \dmu ( X^i (u), \sum_{ j \ne i}^{\infty} \delta_{X^j(u)}) du 
.\end{align}
See \dref{d:21} for the definition of the logarithmic derivative $ \dmu $. 

For a random point field $ \muN $, we consider a window $ \SOrN = \{ | x | \le \rN \} $, $ 0 < \rN < \infty $, such that 
\begin{align}\label{:16b}&
\limi{ \nN } \rN = \infty 
\end{align}
and introduce the diffusion $ ( X_{ \rN }^{ \nN , i})_{i=1}^{\nnnN }$ 
with the state space $ (\SO _{\rN })^{ \nnnN } $ such that, for $ 1\le i \le \nnnN $, 
where $ \nnnN $ is the number of particles in $ \SOrN $, 
\begin{align}\label{:16c}
 X_{ \rN }^{ \nN , i } (t) - X_{ \rN }^{ \nN , i } (0) = B^i (t) 
 + & 
 \half \int_0^t \dmuN 
 ( X_{ \rN }^{ \nN , i } (u), \sum_{ j \ne i}^{\nnnN } 
 \delta_{ X_{ \rN }^{ \nN , j } (u)}) du 
 \\ \notag + & 
 \half \int_0^t U ^{\nN }( X_{ \rN }^{ \nN , i } (u)) du + 
 \half \mathrm{BC}^{ \nN }( X_{ \rN }^{ \nN , i } (t) ) 
.\end{align}
Here, $ \mathrm{BC}^{ \nN }( X_{ \rN }^{ \nN , i } (t) ) $ 
comes from a boundary condition of $ \partial \SrN $ 
and $ U ^{\nN }$ is a free potential caused by the average of the outside particles. 
Furthermore, $ \dmuN $ is the logarithmic derivative of $ \muN $. 
If $ \dmuN $ is given by an interaction potential $ \Psi ^{ \nN } $ 
with inverse temperature $ \beta $, then $ U ^{\nN }$ becomes 
\begin{align}\label{:16e}&
 U ^{\nN } (x) = \frac{\beta }{2}\int_{\SrN ^c} \nabla \Psi ^{ \nN} (x,y) \rNone (y) dy
,\end{align}
where $ \rNone $ is the one-point correlation function of $ \muN $ 
(see \eqref{:36b}--\eqref{:36d}). 
By construction, the associated unlabeled diffusion 
$ \XXN (t) = \sum_{i=1}^{\nnnN } 
\delta _{ X_{ \rN }^{ \nN , i } (t)}$ is $ \muN \circ \pi _{ \rN }^{-1} $-reversible. 
Here, $ \pi _{ \rN } ( \sss ) = \sss (\cdot \cap \SrN )$. 
We label the particles here in such a way that 
$ |X_{\rN }^{\nN ,i } (0)| \le |X_{\rN }^{\nN , i + 1 } (0)| $ for all $ i $. 
Note that the number $ \nnnN $ is, therefore, random. 

The last term in \eqref{:16c} will vanish as $ \nN \to \infty $ because of \eqref{:16b}. 
The relation between $ \dmuN $ and $ \dmu $ is not transparent, as we saw in \eqref{:12a}, \eqref{:13c}, and \eqref{:15g}. 
Nevertheless, in \corref{l:35B}, we shall prove that there exists a sequence $ \{ \rN \}$ 
satisfying \eqref{:16b} such that $ ( X_{ \rN }^{ \nN , i})_{i=1}^{ m }$ converge weakly in 
$ C([0,\infty); \sS ^m)$ to the first $ m $ components of $ \X $. 
In this sense, the geometric universality in condition \thetag{I} to a random point field $ \mu $ 
yields the dynamical universality to the solution $ \X $ of ISDE \eqref{:16a}.

The idea behind the proof of the dynamical universality is as follows.
One of the main tools for the proof is the generalized Mosco convergence, in the sense of Kuwae--Shioya \cite{k-s}, of Dirichlet forms (cf.\,\cite{kol.2006}). This convergence is equivalent to the strong convergence of semi-groups corresponding to Dirichlet forms. 
The Mosco convergence concept consists of two convergence relations related to Dirichlet forms (see \dref{d:44}). 

For a random point field $ \mu $, there are two canonical Dirichlet forms, called the upper and lower Dirichlet forms (see \sref{s:2}). 
Accordingly, two natural schemes of finite-volume Dirichlet forms exist, and each scheme converges to the limit Dirichlet form. 
We shall prove that these two schemes of Dirichlet forms realize the two convergence relations in the Mosco convergence definition, respectively. We use condition \thetag{I} at this stage.

In addition, the two canonical Dirichlet forms in the limit are the same under the uniqueness in law of weak solutions to the ISDE. 
Hence, we conclude the Mosco convergence from conditions \thetag{I} and \thetag{II}. 

In \cite{k-o-t.udf}, we proved the uniqueness of Dirichlet forms applicable to the current situation. 
This uniqueness theorem is robust and can be applied to random point fields from random matrix theory despite the long-range interaction of these random point fields. 

Note that, in general, dynamical convergence fails under only the weak convergence of measures, even in one-dimensional diffusion.
 A typical example is a homogenization problem. 
Hence, we must assume a stronger convergence of the random point fields $ \muN $, such as that specified by condition (I). Thus, the strong convergence in condition (I) is a valid assumption. 
We note that examples of the universality of random matrices satisfying the uniform convergence of correlation functions on each compact set 
can be found in \cite{dkmvz,dkmvz.2,d-g,d-g2,LL16,L09,shc.2011}. 

In our argument, it is critical to take solutions of ISDEs as the limiting point of the stochastic dynamics. 
There are different constructions of infinite-dimensional stochastic dynamics arising from random matrices for $ d= 1 $ and $ \beta = 2 $. In \cite{k-t.07}, the dynamics were constructed by using spatial--temporal correlation functions. 
Using \cite{o-t.core,o-t.sm,o-t.tail}, we find that this construction defines the same stochastic dynamics as given by solutions of ISDEs. 
In \cite{c-h.2014}, the dynamics were constructed using the Brownian Gibbs property. 
It is plausible, but has not yet been proved, that this construction also defines the same dynamics given by the solutions of ISDEs. 

The remainder of this paper is organized as follows.
In \sref{s:2}, we set up Dirichlet forms and recall the relations among random point fields, unlabeled diffusions, and an SDE representation. 
In particular, two types of unlabeled diffusions are presented. 
In \sref{s:3}, we state the main theorems (\tref{l:31}--\tref{l:39}). 
In \sref{s:4}, we recall the concept of Mosco convergence in the sense of Kuwae-Shioya, before proving \tref{l:31} in \sref{s:5}. 
In \sref{s:6}, we prepare some results related to cut-off Dirichlet forms for the proof of \tref{l:32}--\tref{l:39}. 
In \sref{s:7}, we prove \tref{l:32}--\tref{l:34}, and then in \sref{s:8}, we prove \tref{l:35}--\tref{l:39}. 
In \sref{s:X}, we present a sufficient condition for \As{C3}. 
Finally, \sref{s:9} presents some examples of dynamical universality arising from random matrix theory.

\section{Preliminaries} \label{s:2}
\subsection{Two spatial approximations of Dirichlet forms}\label{s:21}
In this section, we prepare Dirichlet forms and the associated dynamics, 
following \cite{o.dfa,o.isde,k-o-t.udf}. 
Let $ \sS $ be a connected open set in $ \Rd $ such that $ \overline{ \sS }^{ \mathrm{int}} = \sS $. 
We take $ \sS $ as the underlying space, and denote the configuration space over $ \sS $ as $ \sSS $ . 
By definition, $ \sSS $ is the set of Radon measures consisting of sums of point measures: 
\begin{align*}
 \sSS =\{ \sss = \sum_{i} \delta _{\si } \st \si \in \sS , \sss(K) < \infty 
\text{ for any compact set $K$ in $ \sS $}\}
.\end{align*}
Here, we regard the zero measure as an element of $ \sSS $. 
The set $ \sSS $ is equipped with the vague topology, under which $ \sSS $ is a Polish space. 
We set $ \SR = \{ |s| < \rR \} $ and 
\begin{align} &\notag 
\SSRm =\{ \sss \in \sSS \st \sss (\SR )=m\}
.\end{align}

For a set $A\subset\sS $, let $ \map{ \pi_{A}}{ \sSS }{ \sSS }$ be the projection map given by 
$ \pi_{A}(\sss )=\sss (\cdot \cap A) $. We often write $ \piR = \pi_{ \SR }$. 
A function $ f $ on $ \sSS $ is said to be local if $ f $ is $ \sigma [\pi_{K}]$-measurable 
for some compact set $ K $ in $ \sS $. 
For such a local function $ f $ on $ \sSS $, $ f $ is said to be smooth if $ \check{f}=\check{f}_{O }$ is smooth for a relatively compact open set $ O \subset \sS $ such that $ K \subset O $. 
Here, $ \check{f}_{O } $ is a function defined on $ \cup_{ \kkk =0}^{ \infty} O^\kkk $ such that, 
 for each $ k $, $ \check{f}_{O } (x_1,\ldots x_{ \kkk }) $ restricted on $ O^\kkk $ is symmetric in 
$ x_i $ ($ i =1,\ldots, \kkk $) and $ \check{f}_{O }(x_1,\ldots,x_{ \kkk }) = f (\mF{x}) $, 
where $ \sum _{i=1}^{ \kkk } \delta _{x_i} = \pi _{ O } (\mF{x} ) $. 
The case $ \kkk =0$ corresponds to a constant function.
Because $ \mF{x}$ is a configuration and $ O $ is relatively compact, the cardinality of the particles of 
$ \mF{x}$ is finite in $ O $. 
Note that $ \check{f}_{O } $ has the consistency property such that 
\begin{align} &\notag 
\check{f}_{O }(x_1,\ldots,x_{ \kkk }) = \check{f}_{O' } (x_1,\ldots,x_{ \kkk }) 
\quad \text{ for all }(x_1,\ldots,x_{ \kkk }) 
\in O^\kkk \cap O^{'\kkk }
.\end{align}
Thus, we see that $ f (\mF{x}) = \check{f} (x_1,\ldots x_{ \kkk }) $ is well-defined.

Next, we introduce carr\'{e} du champ operators on $ \sSS $. 
Let $ a = (a_{pq})_{p,q=1}^d $ be an $ \R ^{d^2}$-valued function defined on $ \SSS $ 
such that $ a_{pq}=a_{qp}$ and $ a $ is elliptic and bounded: there exists a constant $ \Ct \label{;21B}$ such that, for all $ \xs \in \SSS $, 
\begin{align}\label{:21c}&
\cref{;21B}^{-1} |\xi |^2 \le (a \xi , \xi )_{ \Rd } \le \cref{;21B} |\xi |^2 
\quad \text{ for all } \xi \in \Rd 
.\end{align}
Because $ a (x ,\sss ) = a (x, \sum_{i} \delta_{\si } ) $ is symmetric in $ (\si )_i $, 
we can construct a function, denoted by the same symbol $ a (x, s_1,s_2,\ldots ) $ such that 
$ a $ is symmetric in $ ( s_1,s_2,\ldots ) $ for each $ x \in \sS $ and 
\begin{align}\label{:21d}&
 a (x, s_1,s_2,\ldots ) = a (x, \sum_{i} \delta_{\si } ) 
.\end{align} 

Let $ \di $ be the set of all local, smooth functions on $ \sSS $. 
For $ f,g \in \di $, we set
\begin{align}& \label{:21i}
 \DDDa [f,g] (\sss ) = \half \sum _{\si \in \sS } 
(a (\si ,\sss _{i\dia }) 
 \nabla_{\si } \check{f} (\mathbf{s}) , \nabla _{\si } \check{g} (\mathbf{s}) )_{ \Rd }
,\\\label{:21j}&
 \DDDaR [f,g] (\sss )= \half \sum _{\si \in \SR } 
(a (\si ,\sss _{i\dia })
 \nabla_{\si } \check{f} (\mathbf{s}) , \nabla _{\si } \check{g}(\mathbf{s}) )_{ \Rd }
.\end{align}
Here, for $ \mathbf{s}=(\si )_i $, we set 
$ \sss =\sum_i \delta_{\si }$ and $ \sss _{i\dia }= \sum_{j\not=i} \delta_{\sj }$. 
The right-hand sides of \eqref{:21i} and \eqref{:21j} are symmetric functions in 
$ \mathbf{s}=(\si )_i $. Hence, we regard them as functions in $ \sss $. Let 
\begin{align} \label{:21k} &
 \DDDRm [f,g](\sss ) = 
 1_{ \SSRm } (\sss ) \, \DDDaR [f,g] (\sss ) 
.\end{align}
Then, by construction, 
\begin{align} &\notag 
\DDDR = \sum_{m=1}^{ \infty} \DDDRm 
,\quad 
 \limi{ \rR } \DDDR [f,g] = \DDDa [f,g] \quad \text{ for $ f , \g \in \di $}
.\end{align}

A probability measure $ \mu $ on $ \sSS $ is called a random point field. 
For a random point field $ \mu $, we set $ \Lm := L^2(\sSS ,\mu ) $ and 
$ \dimu =\{f\in\di \cap \Lm \st \E (f,f)<\infty \} $. 
For each $ \Rm \in \N $, we define the bilinear forms on $ \Lm $ such that 
\begin{align} &\notag 
\E (f,g)=\int _\sSS \DDDa [f,g] d\mu 
, \quad 
\ER (f,g)=\int _\sSS \DDDR [f,g] d\mu 
\\ &\notag 
 \ERm (f,g) = \int _\sSS \DDDRm [f,g] d\mu 
.\end{align}
We assume the following. 

\smallskip 
\noindent
\As{A1}
$ (\ERm ,\dimu ) $ is closable on $ \Lm $ for each $ \Rm \in \N $. 

\smallskip 

Set $ \Br = \{f\st f\text{ is } \sigma [\piR ]\text{-measurable}\}$. 
If $ f \in \Br $, then $ f (\sss ) $ is independent of $ \piR ^c (\sss ) $. 
Hence, for $ f \in \di \cap \Br $, we have that
\begin{align}&\label{:21x}
 \DDDa [f,f] (\sss ) = \DDDR [f,f] (\sss ) = \DDDRm [f,f] (\sss ) 
\quad \text{ for all } \sss \in \SSRm 
.\end{align}
From \eqref{:21x}, we obtain that, for $ f \in \di \cap \Br $,
\begin{align} &\notag 
 \E (f,f) = \ER (f,f) = \sum_{m = 1 }^{\infty} \ERm (f,f) 
.\end{align}
This obvious identity is one of the key points of the argument in \cite{o.dfa}. 
In the following, we quote a sequence of results from \cite{o.dfa}. 

\begin{lemma}[{\cite[Lemma 2.2]{o.dfa}}] \label{l:21} 
Assume that \As{A1} is satisfied. Then, the following hold: 

\noindent\thetag{1} 
$ (\ER , \dimu ) $ is closable on $ \Lm $. 

\noindent\thetag{2} 
$ (\E , \dimu \cap\Br ) $ is closable on $ \Lm $. 
\end{lemma}
\PF 
Claim \thetag{1} follows from Lemma 2.2 \thetag{1} in \cite{o.dfa}. 
Claim \thetag{2} follows from Claim \thetag{1} and 
$ \dimu \cap\Br \subset \dimu $. 
\PFEND

We write $ (\E ^{1} ,\dom ^{1} )\le (\E ^{2}, \dom ^{2}) $ if 
$ \dom ^{1} \supset \dom ^{2} $ and $ \E ^{1}(f,f)\le \E ^{2}(f,f) $ 
for any $ f\in\dom ^{2} $. 
For a sequence $ \{(\E ^n ,\dom ^n )\}_{n\in\N}$ of positive definite, symmetric bilinear forms on $ \Lm $, we say that $ \{(\E ^n ,\dom ^n )\}$ is increasing if $ (\E ^n ,\dom ^n )\le (\E ^{n+1}, \dom ^{n+1}) $ for any $n\in\N$, and decreasing if $ (\E ^n ,\dom ^n )\ge (\E ^{n+1}, \dom ^{n+1}) $ 
for any $n\in\N$. 

Taking \lref{l:21} into account, we denote the closures of $ (\ER , \dimu ) $ and $ (\E , \dimu \cap\Br ) $ 
on $ \Lm $ as $ (\lER , \ldR ) $ and $ (\ER , \dR ) $, respectively. 
Note that $ (\lER , \ldR ) $ is an extension of $ (\ER , \dR ) $ in the sense that 
$ \dR \subset \ldR $ and the value $ \ER (f,f) $ for $ f \in \dR $ 
of $ (\ER , \dR ) $ coincides with that of $ (\lER , \ldR ) $. 
\begin{lemma}[{\cite[Lemma 2.2]{o.dfa}}] \label{l:22} 
Assume that \As{A1} is satisfied. Then, the following hold:

\noindent
\thetag{1}$ \{(\lER , \ldR )\}_{ \rR \in \N }$ is increasing. 

\noindent
\thetag{2}$ \{ (\ER , \dR ) \}_{ \rR \in \N }$ is decreasing.
\end{lemma}

By definition, the largest closable part $ ((\tilde{ \E })_{ \mathrm{reg}}, (\tilde{ \dom })_{ \mathrm{reg}}) $ of a given positive symmetric form $ (\tilde{ \E }, \tilde{ \dom }) $ 
with a dense domain is a closable form such that 
$ ((\tilde{ \E })_{ \mathrm{reg}}, (\tilde{ \dom })_{ \mathrm{reg}}) $
 is the largest element of closable forms dominated by $ (\tilde{ \E }, \tilde{ \dom }) $. 
Such a form exists and is unique \cite[Theorem S.15]{rs-1}. 

Let $ \dom _\infty =\bigcup_{ \rR \in \N } \dR $. 
Let $ (\E _\infty, \dom _\infty) $ be the symmetric form such that
\begin{align*}
\E _\infty(f,f) =\limi{ \rR }\ER (f,f) 
.\end{align*}
From \lref{l:21} and \lref{l:22}, we obtain the following. 
\begin{lemma}[\cite{o.dfa}] \label{l:23} Assume \As{A1}. Then, the following hold: \\
\thetag{1} $ (\E ,\dimu ) $ is closable on $ \Lm $. \\
\thetag{2} 
The closure of $ ((\E _\infty)_{ \mathrm{reg}}, (\dom _\infty )_{ \mathrm{reg}}) $ on $ \Lm $ 
coincides with that of $ (\E ,\dimu ) $ on $ \Lm $. 
\end{lemma}

From \lref{l:23} \thetag{1}, we denote the closure of $ (\E ,\dimu ) $ on $ \Lm $ as $ (\E ,\dom ) $. 

By \lref{l:22} \thetag{2}, $ \{(\lER , \ldR )\}_{ \rR \in \N }$ is increasing. 
Let $ (\lE ,\ld ) $ be the closed symmetric form 
given by the increasing limit of 
$ \{(\lER , \ldR )\}_{ \rR \in \N }$ 
as $ \rR \to\infty $. Then, by definition, 
\begin{align}\label{:23b}&
\lE (f,f) = \limi{ \rR } \lER (f,f) ,\quad 
\ld = \{ f \in \bigcap_{ \rR =1}^{ \infty}\ldR \st \limi{ \rR } \lER (f,f) < \infty \} 
.\end{align}

We say a sequence of positive closed bilinear forms $ \{ (\E ^N , \dom ^N ) \} $ on $ \Lm $ 
converges to a positive closed bilinear form $ (\E , \dom ) $ on $ \Lm $ 
in the strong resolvent sense if the sequence of their resolvents $ \{ R_{ \alpha }^N \}$ converges to 
the resolvent $ R_{ \alpha }$ of $ (\E , \dom ) $ on $ \Lm $ strongly in $ \Lm $ 
for each $ \alpha > 0 $.

Summarizing the above, we obtain the following lemma.
\begin{lemma} \label{l:24} 
Assume that \As{A1} is satisfied. 
Then, the following hold: \\
\thetag{1} 
$ (\E ,\dom ) $ is the strong resolvent limit of $ \{(\ER , \dR )\}_{ \rR \in \N }$ as $ \rR \to\infty $.
\\\thetag{2} 
$ (\lE ,\ld ) $ is the strong resolvent limit of 
$ \{(\lER , \ldR ) \}_{ \rR \in \N }$ 
as $ \rR \to\infty $.
\\\thetag{3} 
$ (\lE ,\ld ) \le (\E ,\dom ) $.
\end{lemma}

\PF 
The first two statements follow from Theorem 3 in \cite{o.dfa}. 
The third follows from Remark \thetag{3} in \cite[120 p]{o.dfa}. 
\PFEND

Taking \lref{l:24} \thetag{3} into account, we call 
$ (\lE ,\ld ) $ and $ (\E ,\dom ) $ the lower and upper Dirichlet forms, respectively. 
We make an assumption. 

\smallskip 
\noindent \As{A2} $ (\lE ,\ld ) = (\E ,\dom ) $. 
\smallskip 

We shall present a sufficient condition for \As{A2} in \lref{l:29}.

\subsection{A diffusion associated with the upper Dirichlet form}\label{s:22}
In \sref{s:22}, we introduce the $ \sSS $-valued diffusion $ \XX $ 
given by the Dirichlet form $ (\E ,\dom ) $ on $ \Lm $, and the associated $ \SN $-valued 
labeled process $ \X $. 

We denote the density function of $ \mu $ on 
$ \SSRm $ with respect to the Lebesgue measure on $ \SRm $ as $ \sigmaRm $, that is, 
$ \sigmaRm $ is the symmetric function such that 
\begin{align}&\notag 
\frac{1}{m!}\int _{ \SRm } f_{ \rR }^m ( \xm ) 
\sigmaRm ( \xm )d\xm = \int_{ \SSRm }f (\mF{x}) d\mu 
\end{align}
for any bounded $ \sigma [\piR ]$-measurable functions $ f $, where $ f_{ \rR }^m $ 
is a symmetric function on $ \SRm $ such that $ f_{ \rR }^m (\xm ) = f (\mF{x}) $ 
for $ \xm = (x_1,\ldots,x_m) $ and $ \mF{x} = \sum_{i=1}^m \delta_{x_i} \in \SSRm $. 
We set the following condition:

\ms
\noindent
\As{A3} 
$ \mu $ has a density function $ \sigmaRm $ for each $ \Rm \in \N $, and $ \mu $ satisfies 
\begin{align} &\notag 
\sum_{m=1}^{ \infty}m \mu (\SSRm )<\infty \text{ for each } \rR \in\N 
.\end{align}

We recall the concepts of quasi-regularity and locality of Dirichlet forms \cite{mr,fot.2}. 
Quasi-regularity and locality guarantee the existence of the associated diffusion. 
Here, a diffusion process is a strong Markov process with continuous sample paths 
starting at each point in the state space. 

We say that a diffusion $ (\XX , \PP ) = (\{ \XX (t)\} , \{ \PPs \}) $ is associated with 
the Dirichlet form $ (\E ,\dom ) $ on $ \Lm $ if 
$T_t f (\sss ) = \EEs [f (\XX (t) )]$ for any $ f \in \Lm$, where $ \{ T_t \} $ is the $ \Lm $-semi-group 
given by the Dirichlet form $ (\E ,\dom ) $ on $ \Lm $ and $ \EEs $ 
is the expectation with respect to $ \PPs $. By definition, $ \PPs (\XX (0) = \sss ) =1$. 
We say that $ \XX (t) $ is $ \mu $-reversible if $ T_t $ is $ \mu $-symmetric and 
$ \mu $ is an invariant probability measure of $ T_t $. 
The unlabeled diffusion $ (\XX , \PP ) $ associated with $ (\E ,\dom ) $ on $ \Lm $ is constructed in \cite{o.dfa,k-o-t.udf}. 
\begin{lemma}[{\cite[Theorem 1, Corollary 1]{o.dfa}, \cite[Lemma 2.5]{k-o-t.udf}}] \label{l:25}
Assume that \As{A1} and \As{A3} hold. 
Then, $ (\E , \dom ) $ is a local, quasi-regular Dirichlet form on $ \Lm $.
In particular, there exists an $ \sSS $-valued, $ \mu $-reversible diffusion 
$ (\XX , \PP ) $ associated with $ (\E ,\dom ) $ on $ \Lm $. 
\end{lemma}

Let $ \ulab $ be the map defined on $ \{ \cup_{m=0}^\infty \Sm \}\cup\sS ^\N $ such that 
$ \ulab ( \mathbf{s} )=\sum_{i}\delta _{\si }$ for $ \mathbf{s} =(\si )_{i }$. 
Here, $ \sS ^0 = \{\emptyset \}$ and $ \ulab (\emptyset ) = \mathfrak{0}$, 
where $ \mathfrak{0}$ is the zero measure. We call $ \ulab $ the unlabeling map. 

Let $ \SSsi $ be the subset of $ \sSS $ consisting of the single and infinite configurations. 
By definition, 
$ \SSsi = \SSs \cap \SSi $, where $ \SSs $ and $ \SSi $ are given by 
\begin{align} \label{:26g} &
\SSs = \{ \sss \in \sSS \, ;\, \sss ( \{ x \} ) \le 1 \text{ for all } x \in \sS \} ,\quad 
 \SSi = \{ \sss \in \sSS \, ;\, \sss (\sS ) = \infty \} 
.\end{align} 
We denote the set consisting of $ \mathfrak{A}$-valued continuous paths on $ [0,\infty) $ as $ W (\mathfrak{A})$. 
Each $ \ww \in \WSs $ can be written as $ \ww (t) = \sum_i \delta_{w^i (t)}$, 
where $ w^i $ is an $ \sS $-valued continuous path defined on an interval $ I_i $ of the form 
$ [0,b_i) $ or $ (a_i,b_i) $, where $ 0 \le a_i < b_i \le \infty $. 
Note that intervals of this form are unique up to labeling. 
Additionally, note that if $ \partial \sS = \emptyset $, then 
$ \lim_{t\downarrow a_i} |w^i (t)| =\infty $ and 
$ \lim_{t\uparrow b_i} |w^i (t)| =\infty $ for $ b_i < \infty $ for all $ i $. 
For $ \mathfrak{A} \subset \SSs $, we set 
\begin{align} \label{:26i}&
 W _{ \mathrm{NE}} (\mathfrak{A} ) = 
 \{ \ww \in W (\mathfrak{A} ) \, ;\, I_i = [0,\infty ) \text{ for all $ i $}\} 
.\end{align}
We say that the tagged path $ w^i $ of $ \ww $ {\it does not explode} if $ b_i = \infty $, 
and {\it does not enter} if $ I_i = [0,b_i) $. 
By definition, $ \WSsiNE $ is the set consisting of non-exploding and non-entering paths of 
infinitely many non-colliding particles. 

The measurable map 
$ \map{ \lab }{ \SSs \backslash \{ \mathfrak{0} \} } 
{ \{ \cup_{m=1}^\infty \Sm \} \cup \sS ^\N }$ 
is called a label if 
$ \ulab \circ \lab (\sss ) = \sss $ for all $ \sss \in \SSs \backslash \{ \mathfrak{0} \} $. 
For a label $ \lab $, we set 
\begin{align}& \label{:26l}
 \map{ \lpath }{ W _{ \mathrm{NE}} ( \SSs \backslash \{ \mathfrak{0} \} ) }
 {C([0,\infty); \{ \cup_{m=1}^\infty \Sm \} \cup \sS ^{ \N } )} 
\end{align}
by 
\begin{align}& \label{:26m}
\text{$ \lpath (\ww )(0) = \lab (\ww (0)) $, \quad $ \ulab ( \lpath (\ww ) (t)) = \ww (t) $.}
\end{align}
The map $ \lpath $ is uniquely determined by $ \lab $. 


Let $ \Pmu $ be the distribution of the unlabeled diffusion with $ \XX (0) \elaw \mu $, where 
 $ X \elaw \mu $ denotes that the distribution of a random variable $ X $ equals $ \mu $. 
 We assume that the stochastic process $ \XX $ is defined on $ \WSS $ such that 
$ \XX (\ww )(t) = \ww (t) $ for $ \ww =\{ \ww (t) \} \in \WSS $. 
We make the following assumption.

\ms
\noindent 
\As{A4} $ \Pmu ( \WSsiNE ) = 1 $.

\begin{lemma} \label{l:26}
Assume that \As{A1}, \As{A3}, and \As{A4} hold. 
Let $ (\XX , \PP ) $ be the diffusion in \lref{l:25}. 
Then, there exists a unique labeled process $ \X \in C([0,\infty);S ^{ \N } ) $ such that 
$ \X (0) = \lab (\sss ) $ for $ \mu $-a.s.\,$ \sss $ 
and $ \ulab (\X (t)) = \XX (t) $. 
\end{lemma}
\PF
Because of \As{A4}, $ \X := \lpath (\XX ) $ 
is well-defined for $ \Pmu $-a.s. 
Let $ \Pmu (\cdot | \XX (0)=\sss ) $ be the regular conditional probability. 
Then, $ \X = \lpath (\XX ) $ under $ \Pmu (\cdot | \XX (0)=\sss ) $ satisfies the claim. 
\PFEND
\begin{remark}\label{r:26}
Recall that $ \sS $ is an open set. If $ \partial \sS \ne \emptyset $, 
then assumption \As{A4} implies that none of the tagged particles $ X^i $ of 
$ \X = ( X^i ) _{ i \in \N } $ hits the boundary $ \partial \sS $. 
\end{remark}

\subsection{ISDE describing two limit stochastic dynamics, and 
identity of the upper and lower Dirichlet forms}\label{s:23}

Once we have established the labeled dynamics, the next task is to describe the dynamics more explicitly. 
We shall present an ISDE representation of the limit labeled dynamics $ \X $ in \lref{l:27}. 
For this, we recall the concept of the logarithmic derivative of $ \mu $ in \dref{d:21}.

A symmetric and locally integrable function $ \map{ \rho ^n }{ \sS ^n}{[0,\infty ) } $ is called the $ \nnn $-point correlation function of $ \mu $ with respect to the Lebesgue measure if $ \rho ^n $ satisfies 
\begin{align}&\notag 
\int_{A_1^{k_1}\ts \cdots \ts A_m^{k_m}} 
\rho ^n ( \xn ) d\xn 
 = \int _{ \sSS } \prod _{i = 1}^{m} 
\frac{ \sss (A_i) ! } {(\sss (A_i) - k_i )!} d\mu 
 \end{align}
for any sequence of disjoint bounded measurable sets $ A_1,\ldots,A_m \in \mathcal{B}(\sS ) $ and a sequence of natural numbers $ k_1,\ldots,k_m $ satisfying $ k_1+\cdots + k_m =\nnn $. 

Let $ \tilde{ \mu }^{[1]}$ be the measure on $ ( \SSS , \mathcal{B}( \SSS ) ) $ determined by
\begin{align} &\notag 
\tilde{ \mu }^{[1]}(A\times B )=
\int_{ B } \sss (A)\mu(d\sss ), 
\quad A\in\mathcal{B}(S), \; B \in \mathcal{B}(\sSS ) 
.\end{align}
If $ \mu$ has a one-point correlation function $ \rho ^1$, then 
there exists a regular conditional probability $ \tilde{ \mu}_x$ of $ \mu $ satisfying
\begin{align}\notag 
&\int_{A}\tilde{ \mu}_x ( B )\rho ^1(x)dx = 
\tilde{ \mu }^{[1]} (A\times B ), 
\quad A\in\mathcal{B}(S), \; B \in \mathcal{B}(\sSS ).
\end{align}
The measure $ \tilde{ \mu}_x$ is called the Palm measure of $ \mu$ \cite{kal}. 
In this paper, we use the probability measure 
$ \mu_x(\cdot) := \tilde{ \mu}_x(\cdot -\delta_x) $ 
instead of $ \tilde{ \mu}_x$. 
We call $ \mu_x $ the reduced Palm measure of $ \mu $ conditioned at $ x $. 
Informally, $ \mu_x $ is given by 
\begin{align} & \notag 
 \mu _{x} = \mu (\cdot - \delta_x | \, \sss (\{ x \} ) \ge 1 )
.\end{align} 
We consider the Radon measure $ \muone $ on $ \sS \times \sSS $ such that 
\begin{align} &\notag 
 \muone (dx d\sss ) = \rho ^1 (x) \mu _{x} (d\sss ) dx 
.\end{align}
We take $ \muone $ instead of $ \tilde{ \mu }^{[1]}$, and call this the reduced one-Campbell measure of $ \mu $.

We write $ f \in L_{ \mathrm{loc}}^p(\muone ) $ if 
$ f \in L^p(\SR \ts \sSS , \muone ) $ for all $ \rR \in\N $. 
We set 
\begin{align}&\notag 
C_{0}^{ \infty}(\sS )\ot \di = 
\{ \sum_{i=1}^{ \nN }f_i (x) g_i (\mF{y})\, ;\, f_i 
\in C_{0}^{ \infty}(\sS ),\, g_i \in \di ,\, \nN \in \N \} 
.\end{align}

Let $ \Brb $, $ \rR \in \N $, be a set consisting of $ \sigma [\piR ]$-measurable, bounded functions. Let $ \mathcal{B}_{\infty }^b $ be the set of bounded $ \mathcal{B}( \sSS ) $-measurable functions. 
We can naturally regard $ \Brb $ and $ \di \cap \Brb $ as function spaces on $ \SSR $. 
For $ 0 < \rR \le \infty $, we set 
\begin{align}\label{:27v}&
\muR = \mu \circ \piR ^{-1}
.\end{align}
We denote the reduced one-Campbell measure of $ \muR $ as $ \muR ^{[1]} $. 
If $ \rR = \infty $, then $ \muR = \mu $ and $ \muR ^{[1]} = \mu ^{[1]} $. 
We recall the concept of the logarithmic derivative $ \dmu $ of $ \mu $ in \cite{o.isde}. 
%
\begin{definition} \label{d:21}
Let $ 0 < \rR \le \infty $. 
An $ \Rd $-valued function $ \dmu _{ \rR } \in L_{ \mathrm{loc}}^1(\muR ^{[1]})^d $ 
is called a {\em logarithmic derivative} $ \dmu _{ \rR }$ of $ \mu $ on $ \SR $ if, for all 
$ h \in C_{0}^{ \infty}(\SR )\ot \{ \di \cap \Brb \}$, 
\begin{align}&\notag 
 \int _{ \SR \ts \sSS } \dmu _{ \rR }(x,\mF{y}) h (x,\mF{y}) \muR ^{[1]} (dx d\mF{y}) = 
- \int _{ \SR \ts \sSS } \nabla _x h (x,\mF{y}) \muR ^{[1]} (dx d\mF{y}) 
.\end{align}
We write $ \dmu = \dmu _{ \rR }$ for $ \rR = \infty $. 
\end{definition}
 
Let $ a $ be as in \eqref{:21d}. 
We set $ \nabla a (x,\sss ) = 
( \sum_{j = 1}^d \PD{a _{ij}}{x_j}(x,\sss ))_{i=1}^d$, where $ x=(x_1,\ldots,x_d) $. 
We make the following assumption. 

\smallskip 
\noindent \As{A5} 
The logarithmic derivative $ \dmu _{ \rR } $ of $ \mu $ on $ \SR $ exists 
for each $ \rR \in \N \cup \{ \infty \} $ and $ \nabla a \in L_{ \mathrm{loc}}^1(\muR ^{[1]} )^d $. 

\ms

Let $ \sigma (x,\sss ) $ be a matrix-valued function such that $ \sigma ^{t}\sigma = a $. 
We consider the ISDE of $ \X = (X^i)_{i\in\N }$ on $ \SN $ defined by 
\begin{align}\label{:27a}
 d X^{i} (t) = \sigma (X^i(t), \XXtid ) dB^i (t) + 
\half \{ \nabla a + a \dmu \} (X^i(t), \XXtid ) dt 
.\end{align}
Here, $ i\in\N$ and $ \XXtid $ denotes $ \sum_{j\neq i}^{ \infty} \delta _{X^j (t)}$. 
We say that a continuous process $ \X =(X^i)_{i\in\N }$ 
defined on a filtered space $ (\Omega , \mathcal{F}, \pP ,\{ \mathcal{F}_t \} ) $ 
is a weak solution of \eqref{:27a} if there exists a Brownian motion 
$ \mathbf{B}=(B^i)_{ i \in \N } $ 
on the same space such that $ (\mathbf{X},\mathbf{B}) $ satisfies the following for all $ t $:
\begin{align}& \notag 
X^i(t) - X^i (0) = \int_0^t \sigma (X^i (u), \XXuid ) dB^i (u) +
 \int_0^t \half \{ \nabla a + a \dmu \} (X^i (u), \XXuid ) du 
.\end{align}
Here, $ i \in \N $ and we implicitly assume that $ \partial \sS = \emptyset $, 
or no tagged particles hit the boundary. 
Otherwise, a boundary term generally appears in \eqref{:27a}. 
As we saw in \rref{r:26}, the assumption stated above follows from \As{A4}. 

\begin{lemma}[\cite{o.isde}]\label{l:27}
Assume that \As{A1} and \As{A3}--\As{A5} hold. 
Then, the labeled process $ \X = \lpath (\XX ) $ under $ \PPs $ 
given by \lref{l:26} is a weak solution of \eqref{:27a} 
such that $ \X (0) = \lab (\sss ) $ for $ \mu $-a.s.\,$ \sss $. 
\end{lemma}
\PF
If $ \sS = \Rd $, then the claim follows from Theorem 26 in \cite{o.isde}. 
The proof of the case $ \sS \ne \Rd $ is the same as for the case $ \Rd $ because of \As{A4}, and is therefore omitted.
\PFEND

Note that the solution $ \X $ in \lref{l:27} 
is associated with the upper Dirichlet form $ (\E , \dom ) $ in the sense that 
the $ \Lm $-semi-group $ T_t $ given by $ (\E , \dom ) $ on $ \Lm $ 
satisfies $ T_t f (\sss ) = \EE _{\sss } [ f ( \ulab (\X _t)) ]$ 
for any $ f \in \Lm $, where 
$ \EE _{\sss } $ is the expectation with respect to $ \PPs $. 

In \cite[Theorem 3.1]{k-o-t.udf}, it was proved that the lower Dirichlet form $ (\lE ,\ld ) $ 
is also associated with a weak solution $ \X $ of \eqref{:27a} under mild additional constraints. 
Specifically, in \cite{k-o-t.udf}, $ \sigma $ was assumed to be the identity matrix. The generalization of the result in \cite[Theorem 3.1]{k-o-t.udf} to the present case is easy. 

We make the following assumption. 

\smallskip 

\noindent \As{A6} For each $ \rsp \in\N $ with $ r < s $, 
there exists $ \dmu _{ \rsp } \in C_b (\sS \times\sSS ) $ such that 
\begin{align} &\notag 
\limi{r}\limi{s}\limi{ \p }\sup_{R\ge r+s+1} \parallel \dmu _{ \rsp } - \dmu 
\parallel_{L_{ \mathrm{loc}}^1(\sS \times \sSS ,\, \muRsone )}
 = 0 \quad \text{ for $ \mu $-a.s.\! $ \sss $}
.\end{align}
Here, $ \muRsone $ is the reduced one-Campbell measure of 
$ \muRs = \mu (\cdot | \piRc (\sss ) ) $.

Condition \As{A6} is not difficult to check in practice; see \cite[Lemma 6.1]{k-o-t.udf} for a sufficient condition. 
The roles of the parameters $ (\rsp ) $ in \As{A6} are discussed in Section 6 of \cite{k-o-t.udf}. 
Recall that $ (\lE ,\ld ) $ on $ \Lm $ is a Dirichlet form. 
Hence, we have the associated Markovian $ \Lm $-semi-group $ \underline{T}_t $. 
We quote a result from \cite{k-o-t.udf}. 
\begin{lemma}[{\cite[Theorem 3.1]{k-o-t.udf}}] \label{l:28}
Assume that \As{A1} and \As{A3}--\As{A6} hold. 
Then, there exist a continuous $ \sSS $-valued process $ \underline{ \XX } $ and 
a family of probability measures $ \{ \underline{ \PP }_{ \sss } \}$ associated with $ (\lE ,\ld ) $ on $ \Lm $. 
Furthermore, there exists an $ \SN $-valued continuous process $ \underline{ \X } $ 
such that $ \underline{ \X } $ under $ \underline{ \PP }_{ \sss } $ is a solution 
of \eqref{:27a} with $ \underline{ \X }(0) = \lab (\sss ) $ for $ \mu $-a.s.\,$ \sss $. 
\end{lemma}

A key point of \lref{l:28} is that $ (\lE ,\ld ) $ on $ \Lm $ is not necessarily a quasi-regular Dirichlet form. 
If $ (\lE ,\ld ) $ on $ \Lm $ is quasi-regular, 
then the conclusion in \lref{l:28} follows immediately from the result in \cite{o.isde}. 
We do not know how to prove the quasi-regularity directly. 

The quasi-regularity can be proved from the identity $ (\lE ,\ld ) = (\E ,\dom ) $ 
in \lref{l:29}. The proof of \lref{l:29} follows from \lref{l:28} 
and the uniqueness in law of solutions of \eqref{:27a} (see \cite{k-o-t.udf}). 
To state the uniqueness result, we introduce two conditions for solutions $ \X $ of 
ISDE \eqref{:27a} defined on $ (\Omega , \mathcal{F}, \pP ,\{ \mathcal{F}_t \} ) $. 

\smallskip 

\noindent \As{$ \mu $-AC} 
$ \pP \circ \ulab (\X _t) ^{-1} $ is absolutely continuous with respect to 
$ \mu $ for each $ t > 0$. 

\noindent 
\As{NBJ} \ 
$ \pP ( \mrXX < \infty ) = 1 $ for each $ \rR , T \in \mathbb{N} $.

\smallskip 
\noindent 
Here, $ \mr $ 
 is such that, for $ \mathbf{w}=(w^n) $, 
\begin{align}\label{:29r}
\mr (\mathbf{w}) = \inf\{ m \in \N \, ;\, & 
\min_{t\in[0,T]} |w^n (t)|> \rR 
\\ & \notag 
\text{ for all } n \text{ such that } n > m \} 
.\end{align}
Condition \As{NBJ} plays an important role in controlling the labeled dynamics 
$ \X = \lpath (\XX ) $ using the unlabeled dynamics $ \XX $ (see \cite[Section 5]{o-t.tail}). 
Note that \As{NBJ} obviously holds if the cardinality of the particles is finite. 
Even if the number of particle is infinite, \As{NBJ} is easy to check. 
Lemma 10.3.\,in \cite{o-t.tail} states a sufficient condition for \As{NBJ}. 
We make the following assumption for ISDE \eqref{:27a}. 

\smallskip 
\noindent \As{A7} 
The uniqueness in law of weak solutions of ISDE \eqref{:27a} 
with the initial distribution $ \mu \circ \lab ^{-1}$ holds under constraints \As{$ \mu $-AC} and \As{NBJ}. 

\ms

\noindent 
References \cite{o-t.tail,k-o-t.ifc} provide sufficient conditions for \As{A7}. 

We quote a result from \cite{k-o-t.udf}. 
\begin{lemma}[{\cite[Theorem 3.2]{k-o-t.udf}}] \label{l:29}
Assume that \As{A1} and \As{A3}--\As{A7} hold. Then, \As{A2} holds. 
\end{lemma}

\subsection{A sufficient condition for \As{A1}: Quasi-Gibbs measures } \label{s:25}

We introduce a Hamiltonian on a bounded Borel set $ \SR $. 
For Borel-measurable functions $ \map{ \Phi }{ \sS }{ \R \cup \{ \infty \}} $ and 
$ \map{ \Psi }{ \sS \ts \sS }{ \R \cup \{ \infty \}} $ with $ \Psi (x,y) = \Psi (y,x) $, let
\begin{align} &\label{:2Xp}
 \HR ^{ \Phi , \Psi } 
(\mathfrak{x}) = 
\sum_{x_i \in \SR } \Phi ( x_i ) + 
\sum_{x_i , x_j \in \SR , i < j } \Psi ( x_i , x_j ) 
,\quad \text{ where } 
\mathfrak{x} = \sum _i \delta _{x_i } 
.\end{align}
The functions $ \Phi $ and $ \Psi $ are called free and interaction potentials, respectively. 

For a symmetric open set $ \ORm \subset \SRm $, we set 
\begin{align}&\notag 
 \OORm = \piR ^{-1} ( \ulab (\ORm ) ) , \quad 
 \OOR = \bigcup_{m = 1 }^{\infty} \OORm , \quad 
\OO = \bigcap_{ \rR = 1 }^{\infty} \OOR 
.\end{align}
Here, $ \ulab $ is the unlabeling map defined after \lref{l:25}. 
We regard $ \OORm $, $ \OOR $, and $ \OO $ as subsets of $ \sSS $. We note that 
$ \OORm $ and $ \OOR $ are also subsets of $ \SSR := \cup_{m=0}^{\infty} \SSRm $. 
We write $ \mu _1 \le \mu _2 $ if two measures $ \mu _1$ and $ \mu _2 $ on a measurable space $ (\Omega , \mathcal{B}) $ satisfy $ \mu _1(A)\le \mu _2(A) $ for all $ A\in\mathcal{B}$. 
We introduce the concept of a quasi-Gibbs measure. 
\begin{definition}\label{d:22}
A random point field $ \mu $ is said to be a $ (\Phi , \Psi ) $-quasi-Gibbs measure 
with $ \{ \ORm \}_{ \Rm \in \N } $ if \thetag{1}, \thetag{2}, and \thetag{3} hold. \\
\thetag{1} 
 $ \{ \ORm \}_{ \Rm \in \N } $ is a sequence of symmetric open sets such that 
 $ \ORm \subset \SR ^m $ for each $ \Rm \in \N $. 
\\\thetag{2} 
There exists a sequence of measures $ \{ \muk  \} $ on $  \sSS $ such that 
\begin{align}&\notag 
 \muk  \le \mukk \le \mu \quad \text{ for each $ k \in \N $}
,\\ \notag &
\limi{ k } \mukR  = \muR \quad \text{ weakly for each $ \rR \in \N $}
.\end{align}
Here, we set  measures $ \mukR = \muk \circ \piR ^{-1}$ and 
$ \muR = \mu \circ \piR ^{-1} $ on $  \sSS $ as in \eqref{:27v}. 

\noindent 
\thetag{3} The regular conditional probability measures 
\begin{align} & \notag %
\mukRs ^m = 
 \muk  (\piR (\mathfrak{x}) \in \cdot | \, \mathfrak{x}(\SR )= m ,\, 
 \piRc (\mathfrak{x}) = \piRc (\sss ) ) 
\end{align}
satisfy, for all $ k , \rR , m \in \N $ and $ \muk $-a.e.\,$ \sss \in  \sSS $, 
\begin{align}&\label{:2Xx}
\cref{;2y}^{-1} 
 1_{\OORm } e^{-\HR ^{ \Phi ,\Psi } }
 \LambdaRm (d\mathfrak{x}) 
 \le 
 1_{\OORm } 
\mukRs ^m (d\mathfrak{x}) 
\le \cref{;2y} 1_{\OORm } 
 e^{-\HR ^{ \Phi ,\Psi } }
 \LambdaRm (d\mathfrak{x}) 
.\end{align}
Here, $ \piRc (\sss ) = \sss (\cdot \cap \SR ^c )$. We set 
$ \Ct \label{;2y} = \cref{;2y}( k , \rR , m ,\piRc (\sss )) $ to be a positive constant and 
$ \LambdaRm = \LambdaR ( \cdot \cap \SSRm )$, where 
 $ \LambdaR $ is the Poisson random point field with intensity $ 1_{ \SR }dx $. 
\end{definition}

\begin{definition}[{\cite{o.rm}}]\label{d:23}
We call $ \mu $ a $ (\Phi , \Psi )$-quasi-Gibbs measure 
if we take $ \ORm = \SRm $ in \dref{d:22}. 
In this case, $ \OORm = \SSRm $ and $ \OO = \sSS $. 
\end{definition}

We say that a random point field $ \mu $ satisfies 
\As{QG} with $ \{ \ORm \}_{ \Rm \in \N } $ if the following hold: 
\\
\As{QG1} $ \mu $ is a $ (\Phi , \Psi ) $-quasi-Gibbs measure with 
$ \{ \ORm \}_{ \Rm \in \N } $. 
\\
\As{QG2} There exists a potential $ ( \Phi _0 , \Psi _0 ) $ such that 
$ ( \Phi _0 , \Psi _0 ) $ is locally bounded from below, 
 upper semi-continuous, and satisfies
\begin{align}\label{:2Xa}&
\cref{;2X} ^{-1} \Phi _0 ( x ) \le \Phi ( x ) \le \cref{;2X} \Phi _0 ( x ) , \ x \in \SR 
,\\ \notag &
\cref{;2X} ^{-1} \Psi _0 ( x - y ) \le \Psi ( x , y ) \le \cref{;2X} \Psi _0 ( x - y ), 
 \ x , y \in \SR , 
\end{align}
for each $ \rR $ with a constant $ \Ct := \cref{;2X}( \rR ) \label{;2X} > 0 $ 
depending on $ \rR $. 

\ms

Note that $ (\Phi , \Psi ) $ is unbounded in general, and that 
$ e^{-\HR ^{ \Phi ,\Psi } (\ulab (\x ))}$ is bounded and lower semi-continuous on 
$ \SRm $ for each $ \Rm \in \N $. 
Replacing $ \mu $ by $ \mu ( \cdot \cap \OO ) $ in the definition of $ (\ER , \dimu ) $, 
we set $ (\ER ^{\OO }, \di ^{\OO })$. 
The following result was essentially obtained in \cite{o.rm}. 
\begin{lemma} \label{l:2X}
Assume that $ \mu $ satisfies \As{QG} with $ \{ \ORm \}_{ \Rm \in \N } $. 
Then, $ (\ER ^{\OO }, \di ^{\OO })$ is closable on 
$ L^2 (\mu (\cdot \cap \OO )) $ and $ \Lm $ for each $ \rR \in \N $. 
\end{lemma}
\PF 
In the same manner as for Lemmas 3.4 and 3.5 in \cite{o.rm}, 
we can prove that $ (\ER ^{\OO }, \di ^{\OO })$ is closable on $ L^2 (\mu (\cdot \cap \OO )) $, which implies the first claim. 
Note that, if a sequence $ \{ f_n \} $ strongly  converges to zero in $ L^2 (\mu ) $, 
then $ \{ f_n \} $ strongly  converges to zero in $ L^2 (\mu (\cdot \cap \OO )) $. 
Using this, we deduce the second claim from the first. 
%
%
\PFEND

\section{Main results }\label{s:3}

In this section, we set up the problem and state the main results (\tref{l:31}--\tref{l:39}). 

\subsection{Universality of unlabeled dynamics }\label{s:31} 
Let $ \sS $, $ \sSS $, and $ \di $ be as in \sref{s:21}. 
Let $ \{ \muN \}_{ \nN \in \N }$ be a sequence of random point fields on $ \sS $. 
Let $ \DDDRm $ be the carr\'{e} du champ operator in \eqref{:21k}. 
We set 
\begin{align}&\notag 
\ER ^{ \nN , m } (f,g)=\int _\sSS \DDDRm [f,g] d\muN 
,\\ &\notag 
\dimuNm = \{ f \in \di \cap \LmN \, ;\, \ERNm (f,f) < \infty \} 
.\end{align}
The bilinear forms $ ( \EN , \di ^{ \nN } ) $ and $ ( \ERN , \di ^{ \nN } ) $ are given by 
\begin{align}& \label{:31h} 
\EN (f , g ) = \int_{ \sSS } \DDDa [ f , g ] d \muN 
, \quad 
\ERN (f,g)=\int _\sSS \DDDR [f,g] d \muN 	
, \\ \notag &%
\dimuN = \{ f \in \di \cap \LmN \, ;\, \EN (f,f) < \infty \} 
.\end{align}
For the existence of $ \muN $-reversible diffusion, we assume the following. 

\smallskip 
\noindent
\As{B1} \thetag{1} 
$ (\ER ^{ \nN , m } ,\dimuNm ) $ is closable on $ \LmN $ for each $ \nN , \Rm \in \N $. 
\\\thetag{2} 
 $ \muN $ has an $ m $-density function on $ \SR $ for each $ \nN , \Rm \in \N $, and 
 $ \muN $ satisfies 
\begin{align}& \notag 
\sum_{m=1}^{ \infty}m \muN (\SSRm )<\infty \quad \text{ for each $ \nN , \rR \in \N $}
.\end{align}
To take $ \dimu $ as $ \mathcal{C} $ in \dref{d:41}, we assume the following: 

\smallskip \noindent 
\As{B2} $ \dimu \subset \bigcap_{ \nN \in \N }\dimuN $. 

\smallskip \noindent 
Note that \As{B2} is a mild assumption. 
Indeed, if $ \muN (\sss (\sSS ) \le C_{ \nN } ) = 1 $ for some $ C _{ \nN } \in \N $ for each $ \nN \in \N $, then $ \dimuN = \dimu $. Hence, \As{B2} holds. 

To state the main theorems, we introduce cut-off Dirichlet forms and the associated unlabeled diffusions. 
From \As{B1}, we see that $ (\ERN ,\dimuN ) $ is closable on $ \LmN $. 
Then, we denote the closure by $ (\lERN , \ldRN ) $. 
Let $ \SOR = \{ s \in \sS ;\, |s| \le \rR \} $. We set 
\begin{align}& \label{:31i}
 \ldRNW = \{ \f \in \ldRN \, ;\, \text{$ \f $ is $ \sigma [\pi _{\SOR } ]$-measurable}\} 
.\end{align}
Then, $ (\lERN , \ldRNW ) $ is a Dirichlet form on $ \LmN $. 
Clearly, $ \ldRN \supset \ldRNW $. Hence, we have 
\begin{align}\label{:31j}&
(\lERN , \ldRN ) \le (\lERN , \ldRNW ) 
.\end{align}

Let $ \muRN = \muN \circ \piR ^{-1}$. 
Note that $ \muN (\{ \sss \in \sSS ;\, \sss (\partial \SR ) > 0 \} ) = 0 $. 
Then, we can regard $ \muRN $ as a probability measure on $ \CSRbar $, where 
$ \Cf ( A )$ denotes the configuration space over $ A $ for a topological space $ A $. 
We regard $ (\lERN , \ldRNW ) $ as a closed form on 
\begin{align*}&
 L^2 (\muRN ) := L^2 ( \CSRbar , \muRN )
.\end{align*}
It is easy to show that $ (\lERN , \ldRNW ) $ 
is a quasi-regular Dirichlet form on $ L^2 (\muRN ) $. 
Hence, we have a diffusion $ \XXRN $ with the state space $ \CSRbar $ 
associated with the Dirichlet form $ (\lERN , \ldRNW ) $ on $ L^2 (\muRN ) $. 
The associated labeled process $ \lpath (\XXRN )$ is described by SDE \eqref{:16c}. 
Note that the concept of quasi-regularity depends on the topology equipped on the measurable space. Hence, we take $ \CSRbar $ as the state space of the diffusion. 

Let $ \{ \rN \}_{ \nN \in \N } $ be a non-decreasing sequence in $ \N \cup \{ \infty \} $. 
Let $ (\XX , \PP ) $ be the diffusion in \lref{l:25}. We make the following assumptions. 

\smallskip 

\noindent 
\As{B3} The distributions of $ \XXrNN (0) $ and $ \XX (0) $ have densities 
$ \xi _{ \rN }^{ \nN }\in \LmrNN $ and $ \xi \in \Lm $, respectively. 
The functions $ \xi _{ \rN }^{ \nN } $ and $ \xi $ satisfy $ \limi{N}\xi _{ \rN }^{ \nN }= \xi $ 
strongly in the sense of \dref{d:42}. 

\smallskip 

\noindent 
\As{B4} 
For each $ \Rm \in \N $,
\begin{align}\label{:31k}&
\limi{ \nN } \Big\| \frac{ \sigma _{ \rR }^{ \Nm } }{ \sigma _{ \rR }^{m} } -1 \Big\|_{ \SRm } = 0 
.\end{align}
Here, $ \sigma _{ \rR }^{m}$ and $ \sigma _{ \rR }^{\nN ,m}$ are the $ m$-density functions of 
$ \mu $ and $ \muN $ on $ \SR $, respectively. Additionally, 
$ \| \cdot \|_{ \SRm } $ is the $ L^{ \infty}(\SRm , d\x ) $-norm. 
\begin{theorem}\label{l:31} 
Assume that \As{A1}--\As{A3}, and \As{B1}--\As{B4} hold. 
Then, there exists a non-decreasing sequence $ \{ \rN \}_{ \nN \in \N }$ 
in $ \N \cup \{ \infty \} $ satisfying the following:
\begin{align}\label{:31n}& 
 \limi{ \nN } \rN = \infty 
,\\\label{:31r}&
\limi{N} \XXrNN = \XX \text{ in finite-dimensional distributions} 
.\end{align}
Here, using the natural inclusion $ \CSrNbar \subset \sSS $, 
we regard $ \XXrNN $ as a $ \sSS $-valued process. 
\end{theorem}

\begin{remark}\label{r:31} 
We shall give a concrete value of $ \rN $ in \eqref{:52w} and the subsequent sentence. 
With this choice of $ \rN $, we have $ \rN = \infty $ if and only if 
$ \nN \ge \sup \{ \nN _{\nnn } ; \nnn \in \N \} $. 
Furthermore, $ \rN = \infty $ implies $ \muN = \mu $. Thus, if $ \rN = \infty $, then \eqref{:31r} yields a trivial result. 

We present a sufficient condition such that both 
$ \rN = \infty $ for all $ \nN \in \N $ and 
\eqref{:31r} hold for $ \muN $ such that $ \muN \ne \mu $. 
Let $ \Ct \label{;31}(\nN )$ be such that 
\begin{align}\label{:31a}& 
\cref{;31} ( \nN ) = 
\sup \Big\{ \Big\| 
\frac{ \sigma _{ \rR  }^{ \nN , m } }{ \sigma _{ \rR  }^{ m } } -1 
\Big\|_{ \sS _{ \rR  }^m } \, ; \, 1 \le m < \infty ,\, 1 \le \rR  < \infty \Big\}
.\end{align}
We introduce the condition stronger than \As{B4} as follows. 
\begin{align}\label{:31c}&
 \limi{\nN }  \cref{;31} ( \nN ) = 0 
.\end{align}
Let $ \cref{;52} $ be as in  \eqref{:52p}. 
By definition, we deduce  for any $ \kRN \in \N $ 
\begin{align}&\label{:31d}
\text{$ \cref{;52} ( \kRN ) \le  \cref{;31} ( \nN ) $}
.\end{align}
Let $ \cref{;53}$ be as in \eqref{:52a}. 
Then $  \cref{;53} ( \nN ) \le  \cref{;31} ( \nN ) $ from \eqref{:31d}.  
Replacing $  \cref{;53} ( \nN ) $ by $  \cref{;31} ( \nN ) $, we obtain \eqref{:31r} with $ \rN = \infty $ for all $  \nN \in \N $ with a slight modification of the proof in \sref{s:5} using 
$  \cref{;53} ( \nN ) \le  \cref{;31} ( \nN ) $. 
%
%
\end{remark}

Combining \tref{l:31} and \lref{l:29}, we obtain the following. 
\begin{corollary}\label{l:31a}
Assume that \As{A1}, \As{A3}--\As{A7}, and \As{B1}--\As{B4} hold. 
Then, we have the same conclusion as in \tref{l:31}. 
\end{corollary}
Next, we present variants of assumption \As{B4}. 

In general, the density $ \sigmaRm $ in \eqref{:31k} may vanish. 
To control the set $ \{ \sigmaRm = 0 \} $, we use the concept of capacity. 
See \cite[66p]{fot.2} for the definition of capacity. 

Note that $ \SR = \{ s \in \sS ;\, |s| < \rR \} $ is an open set in $ \Rd $. 
Let $ \CSR $ be the configuration space over $ \SR $. 
We equip $ \CSR $ with the vague topology. 
Note that the topology of $ \CSR $ is different from the relative topology 
as a topological subspace of $ \sSS $. 

Let $ (\ER ,\dR ) $ be the Dirichlet form on $ \Lm $ defined before \lref{l:22}. 
Recall that $ (\ER ,\dR ) $ is the closure of 
$ (\E , \dimu \cap \Br ) $ on $ \Lm $. 
Then, each $ f \in \dR $ is $ \sigma [\piR ]$-measurable. 
Hence, we regard $ f $ as a function on $ \CSR $. 

We regard $ \muR $ as a probability measure on $ \CSR $ instead of $ \sSS $, and 
 $ (\ER ,\dR ) $ as a Dirichlet form on $ L^2 ( \CSR , \muR )$. 
Then, we can prove that $ (\ER ,\dR ) $ is a quasi-regular Dirichlet form on $ \LmRC $. 
The proof is similar to that given in \cite{o.dfa} for $ (\E , \dom )$ on $ \Lm $, and is therefore omitted. 

Let $ \mathrm{Cap}_{ \rR }$ be the capacity given by the Dirichlet form 
$ (\ER ,\dR ) $ on $ \LmRC $. 
Let $ \CSRm = \{ \sss \in \CSR \, ;\, \sss (\SR ) = m \} $. We make the following assumption. 
\smallskip 

\noindent \As{ZC} 
For each $ \Rm \in \N $, the capacity $ \mathrm{Cap}_{ \rR }$ satisfies 
\begin{align}\label{:32u}&
\mathrm{Cap}_{ \rR } \big( \{ \sss \in \CSRm \st \sigmaRm (\sss ) = 0 \} \big) = 0 
.\end{align}
Here, $ \sigmaRm $ is regarded as a function on $ \CSRm $ in an obvious manner, that is, 
we set $ \sigmaRm (\sss ) = \sigmaRm (s_1 ,\ldots,s_m ) $ for $ \sss = \sum_{i=1}^m \delta_{\si }$. 

\smallskip 
For $ \Rnu , m \in \N $, let $ \ORnum $ be a symmetric open set in $ \SRm $. We set 
\begin{align} \label{:61z} & 
 \OORnum = \piR ^{-1} ( \ulab ( \ORnum ) ) , \quad 
\OORnu = \bigcup_{m=1}^\infty \OORnum , \quad 
\OOnu = \bigcap_{ \rR =1}^{\infty} \OORnu 
.\end{align}
We introduce the following conditions for $ \{ \ORnum \}_{ \Rnu , m \in \N } $. 
For each $ \Rnu , m \in \N $, 
\begin{align} \label{:63q}&
 \OORnu = \piR^{-1} (\piR (\OORRnu )) 
,\\\label{:63r}&
\OORnu \subset \OORnunu 
,\\ \label{:63s}&
 \mu \Big(\bigcup_{ \nu =1}^{ \infty} \OORnu \Big) = 1 
,\\ \label{:62c}&
\piR (\sSS \setminus \bigcup_{\nu = 1}^{\infty}\OORnu )
= \bigcup_{m=1}^{\infty} \big\{ \sss \in \SSRm \st \sigmaRm ( \sss ) = 0 \big\}
, \ \text{ $ \mu \circ \piR ^{-1}$-a.s}
,\\ \label{:63p}&
0 < \inf \{ \sigmaRm ( \sss ) \st \sss \in \OORnum \} ,\quad 
\sup \{ \sigmaRm ( \sss ) \st \sss \in \OORnum \} < \infty 
.\end{align}
We make the following assumptions regarding the density and correlation functions. 

\ms
\noindent 
\As{B4$'$} 
The density function $ \sigmaRm $ is continuous for each $ \Rm \in \N $ and satisfies 
\begin{align}\label{:32y}&
\limi{ \nN } \Big\| { \sigma _{ \rR }^{ \Nm } }- { \sigmaRm } \Big\|_{ \SRm } = 0 
\quad \text{ for each $ \Rm \in \N $}
.\end{align}
\As{B4$''$} 
The correlation functions 
$ \rNm $ and $ \rho ^m $ are continuous and satisfy 
\begin{align}\label{:33s}&
\limi{ \nN } \big\| { \rho _{}^{ \Nm } }- { \rho _{}^m } \big\|_{ \SRm } = 0 
\quad \text{ for each $ \Rm \in \N $}
, \\ \label{:33t}&
\sup_{ \nN \in \N ,\, \xm \in\SRm } \rNm (\xm ) \le 
 \cref{;21a}^m m^{ \cref{;21b}m} 
\quad \text{ for all } m \in \N 
.\end{align}
Here, 
$ \Ct = \cref{;21a} ( \rR ) \label{;21a} $ and $ \Ct = \cref{;21b} ( \rR ) \label{;21b} $ 
are constants satisfying $0 < \cref{;21a} < \infty $ and $ 0 < \cref{;21b}<1$, and 
 $ \rNm $ and $ \rho ^m $ are the $ m $-point correlation functions of $ \muN $ and 
$ \mu $, respectively. 

\begin{theorem}\label{l:32} 
Assume that \As{A1}--\As{A3} and \As{B1}--\As{B3} hold, and that 
\As{ZC} is satisfied. 
Assume that $ \mu $ and $ \muN $ satisfy \QGO. 
Finally, assume that either \As{B4$'$} or \As{B4$''$} holds. 
Then, the same conclusion as for \tref{l:31} holds. 
\end{theorem}

 \begin{remark}\label{r:32}
\thetag{1} 
Let $ ( \overline{\PP }_{ \rR } , \overline{\XX }_{ \rR })$ be the diffusion associated with 
the Dirichlet form $ (\ER ,\dR ) $ on $ \LmRC $. 
From the general theory of Dirichlet forms in \cite{fot.2}, \eqref{:32u} implies that
\begin{align} &\notag 
 \overline{\PP }_{ \rR }
 ( \overline{\XX }_{ \rR } (t) \in \{ \sss \in \CSRm \st \sigmaRm (\sss ) = 0 \} 
 \ \text{ for some } t ) = 0 
.\end{align}
\thetag{2} 
Let $ \mathrm{Cap}$ be the capacity associated with $ ( \E , \dom )$ on $ \Lm $. 
Note that \eqref{:32u} implies 
$ \mathrm{Cap} \big( \{ \sss \in \SSRm \st \sigmaRm (\sss ) = 0 \} \big) = 0 $ 
 for each $ \Rm \in \N $. 
 This follows from the relation $ ( \E , \dom ) \le ( \ER , \dR )$ by \lref{l:22}, 
 the variational formula of capacity, and the natural identification between 
 $ \CSRm $ and $ \SSRm $. 
\end{remark}

We set $ ( \E _{ \rN }^{ \nN }, \dom _{ \rN }^{ \nN })$ similarly to $ (\ER , \dR ) $ in \lref{l:22}, replacing $ \mu $ and $ \rR $ by $ \murNN $ and $ \rN $, respectively. 
Note that 
\begin{align*}&
 (\lE _{ \rN }^{ \nN } , \widetilde{\ld }_{ \rN}^{ \nN } ) \le 
 (\E _{ \rN }^{ \nN } , \dom _{ \rN}^{ \nN } ) 
.\end{align*}
Let $ T_t $ be the Markovian semi-group associated with $ (\E ,\dom ) $ on $ \Lm $.
Let $ ( \XX , \PP ) $ be the diffusion associated with $ (\E ,\dom ) $ on $ \Lm $. 
The following demonstrates the universality of boundary conditions on $ \partial \sS _{ \rN } = \{| x | = \rN \} $. 
\begin{theorem}	\label{l:34}
Under the same assumptions as for either \tref{l:31} or \tref{l:32}, let $ \4 $ be a (not necessarily quasi-regular) Dirichlet form on $ \LmrNN $ such that 
\begin{align}\label{:34a}&
 (\lE _{ \rN }^{ \nN } , \widetilde{\ld }_{ \rN}^{ \nN } ) \le \4 \le 
 (\E _{ \rN }^{ \nN } , \dom _{ \rN}^{ \nN } ) 
.\end{align}
Let $ \widehat{T}_{ \rN , t}^{ \nN } $ be the Markovian semi-group associated with 
$ \4 $ on $ \LmrNN $. 
Then, $ \widehat{T}_{ \rN , t}^{ \nN } $ converges to $ T_t $ for each $ t $ 
strongly in the sense of \dref{d:42}. 
In particular, if there exists a Markov process $ \widehat{\XX } _{ \rN }^{ \nN } $ 
satisfying \As{B3} associated with $ \widehat{T}_{ \rN , t}^{ \nN } $, then 
$ \widehat{\XX } _{ \rN }^{ \nN } $ converges to $ \XX $ in finite-dimensional distributions. 
\end{theorem}

\subsection{Universality of labeled dynamics and SDEs of finite-particle systems} \label{s:32}
In this section, we state several results on the convergence of labeled dynamics. 
We describe the limit dynamics in terms of solutions to an ISDE and 
strengthen the meaning of convergence at the path-space level. 

Let $ \labN $ be a label.
Let $ \lpathN $ be the label path map given by $ \labN $ as in \eqref{:26l} and \eqref{:26m}. 
For the sequence $ \{\labN \}_{\nN \in \N }$ and $ \lab $, we define the discontinuity set by 
\begin{align*}
\mathrm{Disc}[\lab ] = \{ \sss \in \sSS \, ;\,& \text{ $ \sss \notin \SSs $ or 
there exists $ \{\sss _{\nN } \}_{\nN \in \N }$ in $ \SSs $ }
 \\ \notag & \text{such that }
 \limi{\nN } \sss _{\nN } = \sss \text{ and } \limi{\nN } \labN ( \sss _{\nN }) \ne \lab ( \sss ) 
\} 
.\end{align*}
We make the following assumptions. 

\ms\noindent 
\As{C1} $ \mu (\mathrm{Disc}[\lab ]) = 0 $. 

\smallskip
\noindent 
From \As{B3} and \As{C1}, we have that, for each $ m \in \N $, 
\begin{align}\label{:35a}
\limi{N}\murNN \circ (\lab ^{N,1},\ldots,\lab ^{ \Nm })^{-1} = 
 \mu \circ (\lab ^{1},\ldots,\lab ^{m})^{-1} \quad \text{ weakly}
.\end{align}
Without loss of generality, we can take a non-decreasing label as follows:
\begin{align} &\notag 
\text{$ | \lab ^{ \nN , i } | \le | \lab ^{ \nN , i + 1 } | $ for all $ i \in \N $}
.\end{align}

Let $\XB $ be the diffusion given by the Dirichlet form 
$ (\lE _{ \rN }^{ \nN } , \widetilde{\ld }_{ \rN}^{ \nN } ) $ with the initial distribution $ \murNN $. 
Let $ \WSsNE $ be as in \eqref{:26i}. 

\ms
\noindent 
\As{C2} $ \PPrNN ( \XXrNN \in \WSsNE ) = 1 $ for each $ \nN \in \N $. 

\ms
\noindent 
We set the labeled process $ \XrNN = ( X _{\rN }^{\nN , i } )_{i=1}^{\infty }$ to 
\begin{align} & \notag 
\XrNN = ( \lpathN ( \XXrNN ) , o , o , ,\ldots ) \in C([0,\infty);\SN ) 
,\end{align}
where $ o $ denotes a constant path whose value is denoted by the same symbol $ o $. 
Here, $ o \in \sS $ is a dummy point that has no importance. 
From \As{C2}, we see that $ \XrNN $ under $ \PPrNN $ is well-defined. 
Let $ \mr $ be as in \eqref{:29r}. 

\ms
\noindent \As{C3} For each $ \rR , T \in \N $ and $ \epsilon > 0 $, there exists $ l \in \N $ 
such that 
\begin{align} &\notag 
\sup_{\nN \in \N } \PPrNN (\{ \mr ( \lpathN ( \XXrNN ) ) > l \} ) 
< 1 - \epsilon 
.\end{align}
To prove \tref{l:35}, we need the tightness of $ \{ \XXrNN \}_{\nN \in \N }$ 
in $ \WSS $. We deduce this from the tightness of 
$ \{ \XrNN \}_{\nN \in \N }$ in $ C([0,\infty);\SN ) $. 
We use \As{C3} for this. 

\begin{theorem}\label{l:35}
Assume that $ \As{A2}$--$ \As{A4}$ and $ \As{B1}$--$ \As{B3}$ hold. 
Assume that either \As{B4}, \As{B4$'$}, or \As{B4$''$} is satisfied. 
Assume that $ \mu $ and $ \muN $ satisfy \QGO. 
Finally, assume that \As{ZC} and \As{C1}--\As{C3} hold. 
Then, $ \XB $ satisfy 
\begin{align} &\notag 
\limi{N} \XrNN = \X \quad \text{ in law in $ C([0,\infty); \SN ) $}
\end{align}
and, in particular, for each $ m \in \N $,
\begin{align}\label{:35b}&
\limi{N}( X _{ \rN }^{ \nN ,1},\ldots, X _{ \rN }^{ \nN ,m})=( X ^{1},\ldots, X ^{m}) 
\quad \text{ in law in $ C([0,\infty); \Sm ) $}
.\end{align}
Here, $ \X = (X^i)_{i\in\N }$ is a solution of \eqref{:27a} with the initial distribution 
$ ( \xi d\mu ) \circ \lab ^{-1}$. 
\end{theorem}

The most critical assumption in Theorems \ref{l:31}--\ref{l:35} is \As{A2}, which asserts 
the uniqueness of Dirichlet forms. Combining \tref{l:35} and \lref{l:29}, we immediately obtain the following. 

\begin{corollary}\label{l:35A}
Under the same assumptions as for \tref{l:35} and replacing \As{A2} by \As{A5}--\As{A7}, the same conclusion as in \tref{l:35} holds. 
\end{corollary}

We present an SDE describing finite particle systems in \tref{l:35}. 

Assume that $ \muN $ has the logarithmic derivative $ \dmuN _{ \rR }$ on $ \SR $ 
in the sense of \dref{d:21}. 
Let $ \nnnN $ be the number of particles in $ \SO _{\rR }$ at time zero. 
Note that $ \nnnN $ is unchanged as $ t $ increases. 
Hence, we have that $ \nnnN = \XXN _{\rN } (t) (\SO _{\rN })$. 
Let 
\begin{align}\notag &
 \XXNtid = \sum_{ j \ne i }^{ \nnnN } \delta _{ X _{ \rN }^{ \nN , j } (t)} 
.\end{align}
Then, we introduce the SDE of 
$ ( X _{ \rN }^{ \nN , i })_{ i = 1 }^{ \nnnN }$ as follows. 
For each $ 1 \le i \le \nnnN $, 
\begin{align}\label{:36b}
X_{ \rN }^{ \nN , i } (t) - X_{ \rN }^{ \nN , i } (0) 
= & 
\int_0^t \sigma ( X_{ \rN }^{ \nN , i } (u) , \XXNuid ) dB^{ \nN , i }(u) 
\\ \notag & + 
 \half \int_0^t \{ \nabla a + 
 a \dmuN _{ \rN } \} (X_{ \rN }^{ \nN , i } (u) , \XXNuid ) du 
\\ \notag & + 
 \half \int_0^t 
 a (X_{ \rN }^{ \nN , i } (u) , \XXNuid ) \mathbf{n}^{ \rN }( X_{ \rN }^{ \nN , i } ( u ) ) 
L_{ \rN }^{ \nN , i } (d u ) 
.\end{align} 
Here, $ \mathbf{n}^{ \rN } ( x ) $ is the inward normal unit vector at $ x \in \partial \sS _{ \rN }$ 
and $ \dmuN _{ \rN } $ is the logarithmic derivative of $ \muN $ on $ \SrN $. 
Furthermore, 
$ L_{ \rN }^{ \nN , i }$ is the local time of $ X_{ \rN }^{ \nN , i } $ on the boundary 
$ \partial \sS _{ \rN }$. 
That is, $ L_{ \rN }^{ \nN , i } $ is a continuous non-decreasing process such that 
\begin{align} &\notag 
L_{ \rN }^{ \nN , i } ( t ) = \int_0^t 
\mathbf{1}_{ \partial \sS _{ \rN }} ( X_{ \rN }^{ \nN , i } ( u ) ) L_{ \rN }^{ \nN , i } (d u ) 
.\end{align}

If $ \rN = \infty $, then we delete the term $ L_{ \rN }^{ \nN , i } $ from \eqref{:36b}. 
It is easy to show that if $ \dmuN $ is given by a pair of potentials $ ( \Phi ^{ \nN } , \Psi ^{ \nN } ) $ 
such that 
\begin{align*}&
 \dmuN \xs = - \nabla \Phi ^{ \nN } ( x ) - \sum _i \nabla \Psi ^{ \nN } ( x - \si )
,\end{align*}
then the logarithmic derivative 
$ \dmuN _{ \rN } \xs $ of $ \muN $ on $ \sS _{ \rN } $ is given by 
\begin{align}\label{:36d}
 \dmuN _{ \rN } \xs = - \nabla \Phi ^{ \nN } ( x ) & - 
\sum_{\si \in \sS _{ \rN } } \nabla \Psi ^{ \nN } ( x - \si ) 
\\ \notag &
- 
 \int_{ \sS \backslash \sS _{ \rN } } \nabla \Psi ^{ \nN } ( x - y ) \rNone ( y ) dy 
,\end{align}
where $ \rNone $ is the one-point correlation function of $ \muN $, as before. 
Thus, under the reduced one-Campbell measure $ \murNNone $ of 
$ \muN _{\rN }= \muN \circ \pi _{\rN }^{-1} $, we have that
\begin{align*}&
 \dmuN _{ \rN } \xs = \dmuN \xs - 
 \int_{ \sS \backslash \sS _{ \rN } } \nabla \Psi ^{ \nN } ( x - y ) \rNone ( y ) dy 
.\end{align*}
From \eqref{:36b} and \eqref{:36d}, we obtain \eqref{:16c}. 

From \As{A4} and \As{B1}, we find that $ \XrNN $ is a solution of \eqref{:36b}. 
Taking this into account, we make the following assumption. 

\smallskip 

\noindent \As{C4} 
The uniqueness in law of weak solutions of SDE \eqref{:36b} 
with the initial distribution $ \murNN \circ (\labN )^{-1} $ 
holds under constraints \As{$ \murNN $-AC} and \As{NBJ} for each $ \nN \in \N $. 

\smallskip 

\noindent 
Here, conditions \As{$ \murNN $-AC} and \As{NBJ} were given in \sref{s:23}. 
Clearly, we do not need condition \As{NBJ} in \As{C4} if $ \nnnN < \infty $. 

\smallskip 
From \tref{l:35}, we obtain the convergence of solutions of SDEs. %

\begin{corollary}\label{l:35B}
Consider the same assumptions as for \tref{l:35}. In addition, we assume that \As{C4} holds. 
Let $ \XrNN $ be the solution of \eqref{:36b} in \As{C4}. 
Then, the first $ m $ components of $ \XrNN $ satisfy \eqref{:35b} for each $ m \in \N $. 
\end{corollary}
\PF 
Applying \lref{l:27} to the Dirichlet form 
$ (\lE _{ \rN }^{ \nN } , \widetilde{\ld }_{ \rN}^{ \nN } ) $ on $ \LmrNN $, 
we see that $ \XrNN = \lpath ^{ \nN } ( \XXrNN ) $ is a solution of \eqref{:36b}. 
By \As{C4}, all weak solutions are equivalent in law to $ \XrNN $ if they have common initial distributions. 
Hence, Corollary \ref{l:35B} follows from \tref{l:35}. 
\PFEND

\smallskip

Finally, we strengthen the statement of \tref{l:31}. 
\begin{theorem} \label{l:39}
Under the same assumptions as for \tref{l:35}, there exists a non-decreasing sequence $ \{ \rN \}_{ \nN \in \N }$ 
in $ \N \cup \{ \infty \} $ satisfying \eqref{:31n} and 
\begin{align} &\notag 
\limi{N} \XXrNN = \XX \text{ weakly in } \WSS 
.\end{align}
\end{theorem}

\section{Generalized Mosco convergence}\label{s:4} 
We now introduce the concept of generalized Mosco convergence in the sense of Kuwae-Shioya \cite{k-s}. 
Using this concept, we shall prove the main results in \sref{s:3}. 

\begin{definition}\label{d:41}
Let $H_{ \nN } $ ($ \nN \in \N $) and $ H $ be real Hilbert spaces. 
We say that $ \{H_{ \nN }\}_{ \nN \in \N }$ converges to $ H $ 
 if there exists a dense subspace $ \mathcal{C}\subset H$ and a sequence of operators 
\begin{align} &\notag 
\map{ \Phi _{ \nN }}{ \mathcal{C}}{H_{ \nN }}
\end{align}
such that, for any $u \in \mathcal{C}$,
\begin{align} &\notag 
\limi{N}||\Phi _{ \nN }u||_{H_{ \nN }}=||u||_{H}
.\end{align}
\end{definition}

\begin{definition}\label{d:42}
\noindent
(1) We say that a sequence $ \{u_{ \nN }\}$ with $u_{ \nN }\in H_{ \nN }$ strongly converges to $u\in H$ if there exists $ \{ \tilde{u}_{M}\}\subset \mathcal{C} $ such that
\begin{align}\label{:41c}&
\limi{M}||\tilde{u}_M -u ||_{H}=0
,\\\label{:41d}& 
\limi{M}\limsupi{N}||\Phi_{ \nN } \tilde{u}_M -u_{ \nN } ||_{H_{ \nN }}=0
.\end{align}
\noindent
(2) We say that $ \{u_{ \nN }\}$ with $u_{ \nN }\in H_{ \nN }$ weakly converges to $u\in H$ if
\begin{align} &\notag 
\limi{N}(u_{ \nN } ,v_{ \nN })_{H_{ \nN }} =(u,v)_{H}
\end{align}
for any sequence $ \{v_{ \nN }\}$ with $v_{ \nN }\in H_{ \nN }$ that strongly converges to $v\in H$.
\end{definition}

\begin{definition} \label{:43}
Let $ L (H) $ denote the set consisting of linear operators on $ H $. 
We say that a sequence of bounded operators $ \{B_{ \nN }\} $ with $B_{ \nN }\in L(H_{ \nN }) $ strongly converges to an operator $B\in L(H) $ if, for any sequence $ \{u_{ \nN }\}$ with $u_{ \nN }\in H_{ \nN }$ that strongly converges to $u\in H $, $ \{B_{ \nN } u_{ \nN }\}$ strongly converges to $Bu$. 
\end{definition}

Let $ (\E ,\dom ) $ be a non-negative, symmetric bilinear form 
$ \map{ \E}{ \dom \ts \dom }{ \R }$, where $ \dom $ is a subspace of the Hilbert space $ H $. 
We identify the bilinear form $ \E (\cdot , *) $ with the function $ \E \lA \cdot \rA $ 
on $ H $ such that 
\begin{align}\label{:41f}&
\E \lA u \rA =
\begin{cases}
\E (u,u), \quad& u\in\dom ,
\\
\infty , \quad& u\notin \dom 
.\end{cases}
\end{align}
We say that $ \E $ is a bilinear form on $ H $ if the domain of 
$ \E $ is a subset of $ H $. 

\begin{definition}\label{d:44}
We say that a sequence of bilinear forms $ ( \E ^{ \nN } , \dom ^{\nN } ) $ 
on $ H_{ \nN } $, $ \nN \in \N$, 
is Mosco convergent to a bilinear form $ ( \E , \dom ) $ on $ H $ 
if the following conditions hold.

\noindent
(1)
If a sequence $ \{u_{ \nN }\}$ with $u_{ \nN } \in H_{ \nN }$ weakly converges to $u\in H$, then 
\begin{align}\label{:41g}&
\E \lA u \rA \le \liminfi{N} \E ^{ \nN }\lA u_{ \nN } \rA 
.\end{align}
\noindent
(2)
For any $u\in H$, there exists a strongly convergent sequence $ \limi{N}u_{ \nN }=u$ with $u_{ \nN }\in H_{ \nN }$ such that 
\begin{align}\label{:41h}&
\E \lA u \rA =\limi{N} \EN \lA u_{ \nN } \rA 
.\end{align}
\end{definition}
\begin{lemma}[{\cite{kol.2006}}] \label{l:41}
Let $ ( \EN , \dom ^{ \nN }) $ and $ ( \E , \dom ) $ be Dirichlet forms on $ \LmN $ and $ \Lm $, respectively. 
Let $ T_t^{ \nN }$ and $ T_t $ be the associated semi-groups on $ \LmN $ and $ \Lm $, respectively. 
Then, the following are equivalent.

\noindent
(1)
$ \limi{N}\EN =\E $ in the sense of Mosco convergence. 

\noindent
(2)
$ \limi{N}T_t^{ \nN }=T_t $ strongly for all $t>0$.
\end{lemma}
Thus, we see that the Mosco convergence of Dirichlet forms is equivalent to the strong convergence of the associated semi-groups, which implies the convergence of finite-dimensional distributions.

\section{Proof of \tref{l:31}}\label{s:5}
In this section, we prove \tref{l:31}. 
Throughout this section, we assume that \As{A1}, \As{A3}, and \As{B1}--\As{B4} hold. 
We utilize the concept of Mosco convergence in \dref{d:44}. 
We take $H_{ \nN } = \LmrNN $, $H=\Lm $, $ \mathcal{C} = \dimu $, and 
$ \Phi _{ \nN } = \7 \cdot $ in \dref{d:41}, where 
$ \rN $ is given by \eqref{:52w} and $ \7 $ is given before \eqref{:56b}. 
Then, we have that $ \limi{N}H_{ \nN }=H$ in the sense of \dref{d:41}. 

We take the Dirichlet forms $ \EN $ and $ \E $ in \sref{s:3} 
as the bilinear forms $ \EN $ and $ \E $ in \sref{s:4}. 
Hereafter, $ \E (f) $ and $ \DDDa [f]$ denote $ \E (f,f) $ and $ \DDDa [f,f] $, respectively. 
We shall use the same convention for bilinear objects when they appear. 
The difference between $ \E (f) $ and $ \E \lA f \rA $ should be clearly distinguished. 
Note that $ \E ( f )$ is only defined for $ f \in \dom $, while $ \E \lA f \rA $ is defined 
for all $ f \in \Lm $ and 
\begin{align*}&
\E \lA f \rA = \begin{cases} \infty & \text{for } f \notin \dom \\
\E ( f , f ) & \text{for } f \in \dom 
.\end{cases}
\end{align*} 

Let $ \{ \akR \}_{ \kR \in \N } $ be a sequence of natural numbers satisfying 
\begin{align}\label{:51o}&
 \akR < \akRR ,\ \quad \akR < \akkR \quad \text{ for all $ \kR \in \N $}
.\end{align}
For $ \{ \akR \}_{ \kR \in \N } $, we set 
\begin{align}\label{:51x}&
 \KKkR = \{ \sss \in \sSS \, ;\, \sss ( \SR ) \le \akR \} , \quad 
 \KKk = \cap_{ \rR = 1}^{\infty} \KKkR 
.\end{align}

\begin{lemma} \label{l:51} 
There exists a sequence $ \{ \akR \}_{ \kR \in \N } $ satisfying \eqref{:51o} and 
\begin{align}\label{:51a}&
\inf_{ \nN \in \N } \muN ( \KKk ) \ge 1 - \frac{1}{ \kappa } 
,\quad 
\mu ( \KKk ) \ge 1 - \frac{1}{ \kappa } 
\quad \text{ for each $ \kappa \in \N $}
.\end{align}
\end{lemma}
\PF
Recall that a subset $ \mF{A} $ in $ \sSS $ is relatively compact if and only if there exists an increasing sequence of natural numbers $ M_{ \rR } $ such that 
$ \mF{A} \subset \{ \sss ; \sss (\SR ) \le M_{ \rR } \text{ for all }\rR \} $. 
Because $ \muN $ converges weakly to $ \mu $, 
$ \{ \muN \} $ is tight. 
Hence, there exists a sequence $ \{ \akR \}_{ \kR \in \N } $ 
satisfying \eqref{:51o} and \eqref{:51a}. 
\PFEND

For $ \kRN \in \N $, we set $ \Ct \label{;52}( \kRN ) $ to
\begin{align}\label{:52p}& 
\cref{;52} ( \kRN ) = 
\max \Big\{ \Big\| 
\frac{ \sigma _{ r }^{ \nN , m } }{ \sigma _{ r }^{ m } } -1 
\Big\|_{ \sS _{ r }^m } \, ; \, 1 \le m \le \akR ,\, 1 \le r \le \rR ' \Big\}
.\end{align}

\begin{lemma} \label{l:50} For each $ \nnn \in \N $,
\begin{align}\label{:52s}&
\limi{\nN } \cref{;52} ( \nN , \nnn , \nnn , \nnn ) = 0 
.\end{align}
\end{lemma}
\PF 
From \As{B4}, we have that $ \limi{\nN }\cref{;52} ( \kRN ) = 0 $ for each $ \kR , \rR ' \in \N $. 
Then, taking $ \kappa = \rR = \rR ' = \nnn $, we have \eqref{:52s}. 
\PFEND

For each $ \nnn \in \N $, let 
\begin{align}\label{:52t}&
 \nN _{ \nnn } = \min \{ \nN \,;\, \cref{;52} ( \nN ' , \nnn , \nnn , \nnn ) \le 2^{- \nnn } 
 \text{ for all $ \nN ' \ge \nN $}\} 
.\end{align}
From \eqref{:52s} and \eqref{:52t}, we find that $ \nN _{ \nnn } < \infty $ for each $ \nnn \in \N $. 
Furthermore, it is easy to see that $ \{ \nN _{\nnn } \}_{ \nnn \in \N } $ is a non-decreasing sequence. Hence, we denote the inverse function of $ \nN _{\nnn }$ on $ \N $ as $ \rN $ . 
Indeed, we take $ r _{1} = 1 $ for $ \nN < \nN _{1}$ and 
\begin{align}\label{:52w}&
\rN = 1 + \sup \{ \nnn \in \N \,;\,  \nN _{\nnn } \le \nN \le  \nN _{\nnn +1 } \} 
.\end{align}
If $ \{ \nN _{\nnn } \}_{ \nnn \in \N } $ is bounded, then 
we set $ \rN = \infty $ for $ \nN \ge \sup \{ \nN _{\nnn } ; \nnn \in \N \} $. 
We set 
\begin{align}\label{:52a}&
\cref{;53} ( \nN ) =
\begin{cases} 
 \cref{;52} ( \nN , \rN , \rN , \rN ) & \text{ if } \rN < \infty 
,\\
0 & \text{ if } \rN = \infty 
.\end{cases}
\end{align}

\begin{lemma} \label{l:52}
The sequences $ \rN $ and $ \Ct \label{;53} ( \nN ) $, $ \nN \in \N $, satisfy 
\begin{align} \label{:52b}&
\limi{ \nN } \rN = \infty 
,\\ \label{:52c}& 
 \limi{\nN } \cref{;53} ( \nN ) = 0 
.\end{align}
\end{lemma}
\PF 
It is sufficient to assume that $ \rN < \infty $ for all $ \nN \in \N $. 
From \eqref{:52s}, we can deduce that $ \nN _{ \nnn } < \infty $. 
From this and \eqref{:52w}, we have \eqref{:52b}. From \eqref{:52a}, we deduce 
\begin{align}\label{:52f}& 
\cref{;53} ( \nN ) = 
 \cref{;52} ( \nN , \rN , \rN , \rN ) =
 \cref{;52} ( \nN , \nnn , \nnn , \nnn ) \le 2 ^{ - \nnn }= 2 ^{-\rN } 
.\end{align} 
From \eqref{:52b} and \eqref{:52f}, we obtain \eqref{:52c}. 
\PFEND

\subsection{Lower schemes of Dirichlet forms}\label{s:51}
In \sref{s:51}, we check \dref{d:44} (1). 

Let $ \KKkR $ and $ \rN $ be as in \eqref{:51x} and \eqref{:52w}, respectively. We set 
\begin{align}& \label{:53a}
 \fH = 1_{ \KKrN } \f 
.\end{align}
\begin{lemma} \label{l:53} 
Assume that $ \cref{;53} ( \nN ) \le 1/2 $ and $ \f \in \ldW _{ \rN }^{ \nN } $. Then, 
$ \fH \in \ldW _{ \rN }^{ \nN }$ and 
\begin{align}\label{:53b}& 
 | \lErN ( \fH ) - \lErNN (\fH ) | \le \cref{;53} ( \nN ) \, \lErN ( \fH ) 
.\end{align}
\end{lemma}
\PF 
Because $ \f \in \ldW _{ \rN }^{ \nN } $, we have that $ \f $ is $ \sigma [\pirN ]$-measurable. 
Hence, $ \fH $ is also $ \sigma [\pirN ]$-measurable. 
Using this and \eqref{:53a}, we have 
\begin{align}\label{:53e}&
\lErNN ( \fH ) \le \lErNN ( \f ) < \infty 
.\end{align}
This implies that $ \fH \in \ldW _{ \rN }^{ \nN }$. 
Note that $ \ldW _{ \rN }^{ \nN } \cap \di $ is dense in $ \ldW _{ \rN }^{ \nN }$ 
 with respect to 
 \begin{align*}&
 \lErNNone := \lErNN + (\cdot , * )_{\LmrNN }
. \end{align*}
Then, for each $ \f \in \ldW _{ \rN }^{ \nN } $, we have an $ \lErNNone $-Cauchy sequence 
$ \{ \f _{p}\} $ in $ \ldW _{ \rN }^{ \nN } \cap \di $ such that 
\begin{align}& \label{:53i}
\limi{ p } \lErNNone ( \f _{p} - \f ) = 0 
.\end{align}
Note that $ \fH _p -\fH = \widehat{ \f _{p} - \f } $. 
Then from \eqref{:53a}, \eqref{:53e}, and \eqref{:53i}, 
we deduce that $ \{ \fH _{p}\} $ is an $ \lErNNone $-Cauchy sequence satisfying 
\begin{align}& \label{:53j}
\limi{ p } \lErNNone ( \fH _{p} - \fH ) = 
\limi{ p } \lErNNone ( \widehat{ \f _{p} - \f } ) 
\le \limi{ p } \lErNNone ( \f _{p} - \f ) = 0 
.\end{align}
Hence, we have, from \eqref{:53e} and \eqref{:53j}, 
\begin{align}& \label{:53k}
\limi{ p } \lErNNone ( \fH _{p}) = \lErNNone ( \fH ) < \infty 
.\end{align}

Next, we assume that $ \f \in \ldW _{ \rN }^{ \nN } \cap \di $. 
Then, from \eqref{:53a}, we have that 
\begin{align} \notag 
\lErNN (\fH ) 
= &
 \sum_{m=1}^{\infty} \frac{1}{ m! } 
 \int_{ \SrNm } \DDDa [\fH ] \sigma _{ \rN }^{ \nN , m } d\xm 
 \\ \notag 
= &
 \sum_{m=1}^{ \akrN } \frac{1}{ m! } 
 \int_{ \SrNm } \DDDa [\f ] \sigma _{ \rN }^{ \nN , m } d\xm 
\\ \notag =& 
 \sum_{m=1}^{ \akrN } \frac{1}{ m! } 
 \int_{ \SrNm } \DDDa [\f ] 
 \sigma _{ \rN }^{ m } 
 \Big\{ \frac{ \sigma _{ \rN }^{ \nN , m } }{ \sigma _{ \rN }^{ m } } -1 \Big\} 
d\xm + \lErN ( \fH ) 
.\end{align}
Together with \eqref{:52p} and \eqref{:52a}, this implies that
\begin{align} \label{:53g}
 | \lErN ( \fH ) - \lErNN (\fH ) | & = \Big|
\sum_{m=1}^{ \akrN } \frac{1}{ m! } 
 \int_{ \SrNm } \DDDa [\f ] 
 \sigma _{ \rN }^{ m } 
 \Big\{ \frac{ \sigma _{ \rN }^{ \nN , m } }{ \sigma _{ \rN }^{ m } } -1 \Big\} 
d\xm \Big|
\\ \notag & \le 
 \cref{;53} ( \nN ) \, \lErN ( \fH ) 
.\end{align}

Applying \eqref{:53g} to $ \f _p - \f _q $ and noting that
$ \widehat{ \f _p - \f _q } = \fH _p -\fH _q $, 
we obtain 
\begin{align}& \notag 
 | \lErN ( \fH _p -\fH _q ) - \lErNN ( \fH _p -\fH _q ) | \le 
 \cref{;53} ( \nN ) \, \lErN ( \fH _p -\fH _q ) 
.\end{align}
Using this, $ \cref{;53} ( \nN ) \le 1/2 $, the fact that 
$ \{ \fH _{p}\} $ is an $ \lErNNone $-Cauchy sequence satisfying \eqref{:53j}, and \As{B4}, 
we deduce that $ \{ \fH _{p}\} $ is an $ \lErNone $-Cauchy sequence satisfying 
\begin{align} \notag & 
\limi{ p } \lErNone ( \fH _{p} - \fH ) = 0 
.\end{align}
Hence, we have 
\begin{align}\label{:53m}&
\limi{ p } \lErNone ( \fH _{p} ) = \lErNone ( \fH ) < \infty 
.\end{align}
Applying \eqref{:53g} to $ \{ \f _{p}\} $, we have 
\begin{align}\label{:53l}&
 | \lErN ( \fH _{p} ) - \lErNN (\fH _{p} ) | \le \cref{;53} ( \nN ) \, \lErN ( \fH _{p} ) 
.\end{align}
Using \eqref{:53k} and \eqref{:53m} in \eqref{:53l}, 
we obtain \eqref{:53b} for $ \fH $. 
\PFEND

\begin{proposition} \label{l:54} 
\XXX 
Then, $ (\lE , \ld )$ and 
$ \{ ( \lE _{ \rN }^{ \nN } , \ldW _{ \rN }^{ \nN } ) \}_{ \nN \in \N } $ satisfy 
\begin{align} \label{:54a}& 
 \lE \lA \f \rA \le \liminfi{ \nN } \lE _{ \rN }^{ \nN } \lA \fN \rA 
.\end{align}
\end{proposition}
\PF 
If $ \liminfi{ \nN } \lE _{ \rN }^{ \nN } \lA \fN \rA = \infty $, then \eqref{:54a} is obvious. 
Hence, we assume that
\begin{align}\notag & 
 \liminfi{ \nN } \lE _{ \rN }^{ \nN } \lA \fN \rA < \infty 
.\end{align}
From this and \eqref{:41f}, we have 
\begin{align}\label{:54f}& 
 \liminfi{ \nN } \lE _{ \rN }^{ \nN } ( \fN ) < \infty 
.\end{align}

Assume that $ \cref{;53} ( \nN ) \le 1/2 $ and $ \fN \in \ldW _{ \rN }^{ \nN }$. 
Using \lref{l:53}, we have from \eqref{:53b} that
\begin{align} \label{:54g}&
( 1 - \cref{;53} ( \nN ) ) \lErN ( \fNW ) 
\le 
 \lErNN (\fNW ) 
\le 
( 1 + \cref{;53} ( \nN ) ) \lErN ( \fNW ) 
.\end{align}
Note that $ \lE _{ \rN }^{ \nN } ( \fN ) < \infty $ implies $ \fN \in \ldW _{ \rN }^{ \nN }$. 
Hence, using \eqref{:54f}, 
we deduce that $ \fN \in \ldW _{ \rN }^{ \nN }$ for infinitely many $ \nN $. 
From this, we see that \eqref{:54g} holds for infinitely many $ \nN $. 
Combining this with \eqref{:52c} and \eqref{:53e}, we deduce from \eqref{:54g} that
\begin{align} & \label{:54h}
 \liminfi{ \nN } \lErN ( \fNW ) = \liminfi{ \nN } \lErNN (\fNW ) 
 \le \liminfi{ \nN } \lErNN (\fN ) 
 < \infty 
.\end{align}
Because $ \fNW $ is $ \sigma [\pirN ]$-measurable, we see that 
$ \lE (\fNW ) = \lErN ( \fNW ) $. Hence, \eqref{:54h} yields
\begin{align}\label{:54i}& 
 \liminfi{ \nN } \lE (\fNW ) < \infty 
.\end{align}
Furthermore, from \eqref{:52p} and \eqref{:52a}, we obtain
\begin{align} & \notag 
( 1 - \cref{;53} ( \nN ) ) \| \fNW \|_{ L^2 ( \mu ) }
\le 
 \| \fNW \|_{ \LmrNN }
\le 
( 1 + \cref{;53} ( \nN ) ) \| \fNW \|_{ L^2 ( \mu ) } 
.\end{align}
Combining this with \eqref{:52c}, we deduce 
\begin{align} & \label{:54l} 
\limsupi{ \nN } \| \fNW \|_{ L^2 ( \mu ) } = \limsupi{ \nN } \| \fNW \|_{ \LmrNN } 
.\end{align}
Because $ \{ \fN \} $ with $ \fN \in \LmN $ weakly converges to $ \f \in \Lm $, we have 
\begin{align}\label{:54j}&
\sup_{ \nN \in \N } \| \fN \|_{ \LmrNN } < \infty 
.\end{align}
Clearly, $ \| \fNW \|_{ \LmrNN } \le \| \fN \|_{ \LmrNN } $. 
Hence, we have from \eqref{:54l} and \eqref{:54j} that
\begin{align}\label{:54m}&
 \limsupi{ \nN } \| \fNW \|_{ L^2 ( \mu ) }
 \le 
\sup_{ \nN \in \N } \| \fN \|_{ \LmrNN } < \infty 
.\end{align}

We set 
 $ \lEp (\cdot, \cdot) = \lE (\cdot, \cdot) +\p ^{-1} (\cdot, \cdot)_{ \Lm }$ for $ \pN $. 
Then, $ \ld $ is a Hilbert space with the inner product $ \lEp $ for each $ \pN $. 
From \eqref{:54i} and \eqref{:54m}, we deduce 
\begin{align} & \notag 
 \liminfi{ \nN } \lEp ( \fNW ) < \infty 
.\end{align}
Hence, we can choose an $ \lEp $-bounded subsequence from $ \{ \fNW \} $. 

We can choose an $ \lEp $-weak convergent subsequence with the limit $ \f $ 
from an arbitrary, $ \lEp $-bounded subsequence $ \{ \fW _{ \nN '} \} $ of $ \{ \fNW \} $. 
Furthermore, we obtain $ \limi{ \nN } 1_{ \KKrN } (\sss )= 1 $ for $ \mu $-a.s.\,$ \sss $ 
from \eqref{:51a}. Hence, for $ \mu $-a.s.\,$ \sss $, we have that 
\begin{align}\notag &
\f (\sss ) = \limi{ \nN } \fNW (\sss ) , \quad 
\DDDa [ \limi{ \nN } \fNW ] (\sss ) = \limi{ \nN } \DDDa [ \fNW ] (\sss ) 
.\end{align}
Collecting these results together, we can deduce for all $ \p \in \N$ that 
\begin{align}\label{:54w}&
 \lEp (\f ) = \lEp ( \limi{ \nN } \fNW ) \le \liminfi{ \nN } \lEp ( \fNW ) 
.\end{align}
Then, from \eqref{:54m} and \eqref{:54w}, we obtain
\begin{align}\label{:54'}
 \lE (\f ) \le \lEp (\f ) &\le \liminfi{ \nN }\Big\{ \lE (\fNW ) + \frac{1}{ p }
 \| \fNW \|_{ L^2 ( \mu ) }^2 \Big\} 
 \\ \notag &
 \le \liminfi{ \nN } \lE (\fNW ) + \frac{1}{ p } \limsupi{ \nN } \| \fNW \|_{ L^2 ( \mu ) }^2 
.\end{align}
Note that \eqref{:54'} holds for all $ \p \in \N $. 
From \eqref{:54m} and \eqref{:54'}, we then obtain 
\begin{align}\label{:54"}
 \lE (\f ) & \le \liminfi{ \nN } \lE (\fNW ) 
.\end{align}

Putting these together, we obtain
\begin{align}\label{:54(}
 \lE (\f ) & \le \liminfi{ \nN } \lE (\fNW ) \quad \quad \text{by \eqref{:54"}}
\\ \notag &
= \liminfi{ \nN } \lErN ( \fNW ) 
\quad \text{by the $ \sigma [\pirN ]$-measurability of $ \fNW $}
 \\ \notag & \le 
 \liminfi{ \nN } \lErNN ( \fN ) 
 \quad \text{by \eqref{:54h}}
.\end{align}
Finally, we have \eqref{:54a} from \eqref{:54(}. 
\PFEND 

\subsection{Upper schemes of Dirichlet forms and proof of \tref{l:31}}\label{s:52}
In \sref{s:52}, we check \dref{d:44} \thetag{2}. 
Let $ \akR $ be as in \eqref{:51o}. 
For $ \akR $ and a label $ \lab = (\labj )_j $
satisfying $|\lab _j(\sss )|\le |\lab _{j+1}(\sss )|$ for all $j$, we set 
\begin{align*} &
 J_{ \kR } (\sss ) =\{j\st j> \akR ,\lab _j(\sss )\in \SR \} 
,\\ \notag &
\da = \Big\{
 \sum _{j\in J_{ \kR } (\sss ) }(\rR -|\lab _j (\sss ) |)^2 
 \Big\}^\half 
.\end{align*}
Let $ \rho \in C^\infty(\R) $ be a function satisfying $ \rho (t) \in [0,1]$ for any $t\in\R$, 
$ \rho (t)=1 $ for $t\le 0$, $ \rho (t)=0 $ for $t\ge 1$, and 
$ \rho '(t) \le \sqrt{2} $ for any $ t \in \R$. We set 
\begin{align} &\notag 
\chika (\sss )=\rho \circ \da 
.\end{align}
Let $ \KKkR = \{ \sss \in \sSS \, ;\, \sss ( \SR ) \le \akR \} $ be as in \eqref{:51x}. 
We set 
\begin{align} &\notag 
\SSaa = \{ \sss \in \sSS \, ;\, \sss (\sS _{ \rR - 1 } ) \le \akR \} 
.\end{align}
Then, $ \SSa \subset \SSaa $ by construction. 
The next lemma shows that $ \chika $ is a cut-off function of $ \SSa $. 
Let $ \Eone  (f,g) = \E (f,g) + (f,g)_{ \Lm }$. 
We can prove \lref{l:55} in a similar manner to \cite[Lemma 2.5]{o.dfa}. 
\begin{lemma} \label{l:55} 
Let $ \cref{;21B} $ be as in \eqref{:21c}. For each $ \kR \in \N $, the following hold: 

\noindent \thetag{1} 
$ \chika \f \in \dimu $ for each $ \f \in \dimu $. 

\noindent \thetag{2} 
$ \chika =1 $ on $ \SSa $ and $ \chika =0 $ on $ \sSS \backslash \SSaa $. 

\noindent \thetag{3} 
$ \DDDa [\chika ] =0 $ on $ (\SSaa \backslash \SSa )^c $ and 
$0 \le \DDDa [\chika ] \le \cref{;21B} $ on $ \SSaa \backslash \SSa $. 

\noindent \thetag{4} 
$ \Eone  (\chika \f ) \le 2 \E ( \f ) + 3 \cref{;21B}\| \f \|_{\Lm }^2 $ for each $ \f \in \dimu $. 
\end{lemma}
\PF A straightforward calculation shows that
\begin{align}\label{:55f}
\DDDa [ \chika ] \le & \cref{;21B} \DDD [ \chika ]
\\ \notag = & \frac{ \cref{;21B} }{ 2 } 
\Big\{ \frac{ \rho ' ( \da ) }{ \da } \Big\} ^2 
\sum _{j\in J_{ \kR } (\sss ) }(\rR -|\lab _j (\sss ) |)^2 
\\ \notag = & 
 \frac{ \cref{;21B} }{ 2 } \rho ' ( \da ) ^2 \le \cref{;21B} 
.\end{align}
From \eqref{:55f}, we have 
\begin{align}\label{:55g}
\DDDa [ \chika \f ] &= 
 \chika ^2 \DDDa [ \f ] + \f ^2 \DDDa [ \chika ] + 2 \chika \f \, \DDDa [ \chika , \f ] 
 \\ \notag & \le 
 2\{ \chika ^2 \DDDa [ \f ] + \f ^2 \DDDa [ \chika ] \} 
 \\ \notag & \le 
 2 \{ \DDDa [ \f ] + \cref{;21B} \f ^2 \} 
.\end{align}
From \eqref{:55g}, we see that $ \E ( \chika \f ) < \infty $, which implies \thetag{1}. 
\thetag{2} is clear by construction. 
\thetag{3} follows from \thetag{2} and \eqref{:55f}. \thetag{4} follows from \eqref{:55g}. 
\PFEND

Let $ ( \EN , \di ^{ \nN } ) $ be as in \eqref{:31h}. 
Then, using \lref{l:23} \thetag{1}, we see 
from \As{B1} \thetag{1} that $ ( \EN , \di ^{ \nN } ) $ is closable on $ \LmN $. 
We denote the closure of $ ( \EN , \di ^{ \nN } ) $ on $ \LmN $ as $ (\EN , \dN ) $. 
\begin{proposition}	\label{l:56} 
Let $ \rN $ be as in \eqref{:52w}. For each $ \f \in \Lm $, there exists 
a sequence $ \{ \gN \}_{N \in \N }$ satisfying the following: 
\begin{align} \label{:56x}& \ 
\text{$ \gN $ is $ \sigma [ \pi _{\rN }]$-measurable, 
$ \limsupi{\nN} \| \gN \|_{ \LmrNN } < \infty $} 
,\\ &\label{:56y}
\limi{ \nN }\gN =f \text{ strongly in the sense of \dref{d:42}}
,\\&\label{:56z}
\limi{ \nN }\EN \lA \gN \rA = \E \lA \f \rA 
.\end{align}
\end{proposition}
\PF 
If $ f \notin \dom $, then $ \E \lA f \rA = \infty $ and we can easily take $ \gN $ so as to satisfy 
\eqref{:56x}--\eqref{:56z}. 

We next suppose that $ f \in \dom $. Then, $ \E \lA f \rA = \E ( f ) < \infty $. 
Furthermore, there exists a sequence $ \{ \fn \}_{ \nnn \in \N }$ in $ \dimu $ such that 
\begin{align}\label{:56a}&
\limi{ \nnn } \Eone  ( \fn - \f ) = 0 
.\end{align}
Because each $ \fn \in \dimu $ is a local function,
 we can assume that $ \fn $ is $ \sigma [ \pi _{ \nnn }]$-measurable. 

Let $ \{ \qN \} $ be a non-decreasing sequence of natural numbers such that 
$ \qN + 1 \le \rN $ for $ \rN \ge 2 $ and $ \limi{ \nN } \qN = \infty $. 
We set $ \7 = \chi _{ \qN + 1 , \qN } $. Let 
\begin{align}\label{:56b}& 
 \gN = \7 \fqN 
.\end{align}
Because $ \7 $ and $ \fqN $ are 
 $ \sigma [ \pi _{ \qN }]$-measurable and $ \qN + 1 \le \rN $, we find that 
$ \gN $ is $ \sigma [ \pi _{\rN }]$-measurable. 
From $ \fqN \in \Lm $, we have that $ \gN \in \Lm $. 
Combining these results with \As{B4}, we find that $ \limsupi{\nN} \| \gN \|_{ \LmrNN } < \infty $. 
Thus, we obtain \eqref{:56x}.

We next check \eqref{:56y}. 
Recall that we take $ \mathcal{C} = \dimu $. 
For $ \gN = \7 \fqN $, we take $ \widetilde{\g }_{ \mM } = \fqM $. 
Then, $\widetilde{\g }_{ \mM } \in \mathcal{C} $. 
Using \eqref{:56a} and $ \limi{ \mM } \qM = \infty $, we have 
\begin{align}\label{:56c}&
\limi{ \mM } \| \widetilde{\g }_{ \mM } - \f \|_{ \Lm } ^2 = 
\limi{ \mM } \| \fqM - \f \|_{ \Lm } ^2 = 0 
.\end{align}
Let $ \Phi _{ \nN } ( \cdot ) = \7 \cdot $. Then, 
$ \Phi _{ \nN } ( \widetilde{\g }_{ \mM } ) = \Phi _{ \nN } ( \fqM ) = \7 \fqM $ 
by construction and \eqref{:56b}. 
Hence, from this, \eqref{:52p}, and \lref{l:55}, we have that
\begin{align}\label{:56e}&
\limi{ \mM } \limsupi{ \nN } \|\Phi _{ \nN } ( \widetilde{\g }_{ \mM } ) - \gN \|_{ \LmrNN } 
\\ \notag = &
 \limi{ \mM } \limsupi{ \nN } \| \7 \fqM - \7 \fqN \|_{ \LmrNN } 
 \quad \text{ by definition}
 \\ \notag \le & 
 \limi{ \mM } 
 \limsupi{ \nN }
 \cref{;52} ( \nN , \qN + 1, \qN + 1 , \qN ) 
 \| \7 \fqM - \7 \fqN \|_{ \Lm } 
 \text{ by \eqref{:52p}}
 \\ \notag 
 \le & \limi{ \mM } \limsupi{ \nN } 
\cREFa 
 \| \7 \fqM - \7 \fqN \|_{ \Lm } 
=0 
.\end{align}
Here, we have used $ \cref{;52} ( \nN , \qN + 1, \qN + 1 , \qN ) \le \cREFa  $, which 
follows from $ \qN + 1 \le \rN $. 
The last equality follows from \lref{l:55} and \eqref{:56a}. 

Thus, we see that $ \{ \gN \} $ and $ \{ \widetilde{\g }_{ \mM } \} $ satisfy \eqref{:41c} and \eqref{:41d} 
from \eqref{:56c} and \eqref{:56e}, respectively. 
We have already checked $\widetilde{\g }_{ \mM } \in \mathcal{C} $. 
Hence, we obtain \eqref{:56y}.

By \lref{l:55}, \eqref{:56a}, and \eqref{:56b}, we have that
\begin{align} \label{:56g}& 
\sup_{ \nN \in \N } \Eone ( \gN ) = \sup_{ \nN \in \N } \Eone (\7 \fqN ) < \infty 
.\end{align}
Obviously, we have 
\begin{align}\label{:56h} &
 | \EN (\gN ) - \E ( \f ) | \le 
 |\EN (\gN ) - \E (\gN ) | 
+ 
|\E (\gN ) - \E ( \f ) |
.\end{align}
From \lref{l:52}, \eqref{:52p} and \eqref{:56g}, we see that 
\begin{align}\label{:56i}
|\EN (\gN ) - \E (\gN ) |
 \le & \cref{;52} ( \nN , \qN + 1, \qN + 1 , \qN ) 
 \E ( \gN ) 
\\ \notag \le &
\cREFa 
 \E ( \gN ) 
\xrightarrow[\nN \to \infty]{}
 0 
.\end{align}
By a straightforward calculation, \lref{l:55}, \eqref{:56c}, and \eqref{:56g} imply that
 \begin{align}\label{:56j} 
 |\E (\gN ) - \E (\f ) |
 &=
 \Big | 
\int_{ \sSS } 
\DDDa [ \7 \fqN ] - \DDDa [\f ] d\mu 
\Big| 
\\ \notag 
&=
 \Big | 
\int_{ \sSS } 
\DDDa [ \7 ] \fqN ^2 +
2 \DDDa [ \7 , \fqN ] \7 \fqN 
+ \7 ^2 \DDDa [ \fqN ] - \DDDa [\f ] d\mu 
\Big| 
\\ \notag 
&\le 
 \Big | 
\int_{ \sSS } 
\DDDa [ \7 ] \fqN ^2 
 + 2 \DDDa [ \7 , \fqN ] 
 d\mu \Big | 
+ 
 \Big | 
\int_{ \sSS }
\7 ^2 \DDDa [ \fqN ]- \DDDa [\f ] d\mu 
\Big| 
\\ \notag &
\xrightarrow[\nN \to \infty]{}
0 .\end{align}
 From \eqref{:56h}, \eqref{:56i}, and \eqref{:56j}, we obtain \eqref{:56z}. 
This completes the proof. 
\PFEND 

\noindent
{\em Proof of \tref{l:31}. } 
\XXX 
From \pref{l:54}, we see that
$ (\E , \ld )$ and $ \{ ( \lE _{ \rN }^{ \nN } , \ldW _{ \rN }^{ \nN } ) \}_{ \nN \in \N } $ 
satisfy 
\begin{align}\label{:59c}& 
 \E \lA \f \rA \le \liminfi{ \nN } \lE _{ \rN }^{ \nN } \lA \fN \rA 
.\end{align}
By \As{A2}, we have $ (\E , \dom ) = (\lE , \ld ) $. 
From this and \eqref{:59c}, we see that 
$ (\E , \dom )$ and $ \{ ( \lE _{ \rN }^{ \nN } , \ldW _{ \rN }^{ \nN } ) \}_{ \nN \in \N } $ 
satisfy \eqref{:41g}. Thus, we obtain \dref{d:44} \thetag{1}. 

Let $ \gN $ be as in \pref{l:56}. From \eqref{:56z}, we see that
\begin{align}\label{:59e}&
 \E \lA \f \rA = \limi{\nN } \EN \lA \gN \rA 
.\end{align}
Because $ \gN $ is $ \sigma [ \pi _{\rN }]$-measurable, we find that
$ ( \lE _{ \rN }^{ \nN } , \ldW _{ \rN }^{ \nN } ) $ satisfies 
\begin{align}& \label{:59f}
 \EN \lA \gN \rA = \lE _{ \rN }^{ \nN } \lA \gN \rA 
.\end{align}
Combining \eqref{:59e} and \eqref{:59f}, we obtain 
\begin{align*}&
 \E \lA \f \rA = \limi{\nN } \EN \lA \gN \rA = \limi{\nN } \lE _{ \rN }^{ \nN } \lA \gN \rA 
.\end{align*} 
Hence, we obtain \dref{d:44} \thetag{2} for 
$ (\E , \dom )$ and $ \{ ( \lE _{ \rN }^{ \nN } , \ldW _{ \rN }^{ \nN } ) \}_{ \nN \in \N } $. 

Thus, the Mosco convergence in \dref{d:44} holds for 
$ (\E , \dom )$ and $ \{ ( \lE _{ \rN }^{ \nN } , \ldW _{ \rN }^{ \nN } ) \}_{ \nN \in \N } $. 
The Mosco convergence of Dirichlet forms 
implies the strong convergence of the associated $ L^2$-semi-groups, 
which yields the convergence of finite-dimensional distributions of $ \XXrNN $ to $ \XX $ (see \cite[Section 7]{k-s}). 
\PFEnd

\section{Cut-off Dirichlet forms in infinite volumes} \label{s:6} 
In this section, we construct schemes of cut-off Dirichlet forms in infinite volumes. 
We shall use these schemes in \sref{s:7} to prove Theorems \ref{l:32}--\ref{l:39}. 
 
Let $ \OORnu $ be as in \eqref{:61z}. 
We set 
\begin{align}\label{:61a}&
 \ERnu (f,g) = \int _{ \OORnu } \DDDaR [f,g] \, d\mu 
.\end{align}
\begin{lemma} \label{l:61} 
Assume that $ \mu $ satisfies \QGg. Assume that \eqref{:63q} holds. 
Then, $ (\ERnu ,\dimu ) $ is closable on $ \Lm $. 
\end{lemma}
\PF 
Using \lref{l:2X}, we see that $ ( \ER ^{ \OOnu } , \di ^{ \OOnu } ) $ is closable on $ \Lm $. 
From \eqref{:63q} and \eqref{:61a}, we have $ \ER ^{ \OOnu } = \ERnu $. 
Hence, we find that $ (\ERnu , \di ^{ \OOnu } ) $ is closable on $ \Lm $. 
Combining this with $ \dimu \subset \di ^{ \OOnu }$, we conclude that 
$ (\ERnu ,\dimu ) $ is closable on $ \Lm $. 
\PFEND

\subsection{Lower schemes of Dirichlet forms in infinite volumes}\label{s:61}

We assume that \eqref{:63q}--\eqref{:63s} hold. 
From \lref{l:61}, we see that $ (\ERnu , \dimu ) $ is closable on $ \Lm $. 
Hence, we denote the closure of $ (\ERnu , \dimu ) $ on $ \Lm $ as $ (\lERnu , \ldRnuW ) $. 
Then, we deduce from \eqref{:63q} that, for each $ \nuN $, 
\begin{align}\label{:62f}&
 (\lERnu , \ldRnuW ) \le 
 ( \lE _{ \rR ', \nu } , \ld _{ \rR ', \nu } ) 
\quad \text{ for } \rR \le \rR ' 
.\end{align}
From \eqref{:63r}, it is easy to show that, for each $ \rR \in \N $,
\begin{align}\label{:62g}&
 (\lERnu ,\ldRnuW ) \le ( \lE _{ \Rnu '} , \ld _{ \Rnu ' }) 
\quad \text{ for } \nu \le \nu ' 
. \end{align}
From \eqref{:62f} and \eqref{:62g}, we have that, 
for $ f \in \bigcap_{ \Rnu \in \N }\ldRnuW $, 
\begin{align}\label{:62i}& 
\limi{\rR } \{ \limi{\nu } \lERnu (f,f) \} = 
 \limi{\nu }\{ \limi{\rR } \lERnu (f,f) \} 
.\end{align}

For $ f , g \in \dimu $, we set
\begin{align}&\notag 
 \Enu (f,g) = \int _{ \OOnu } \DDDa [f,g] (\sss ) \, d\mu 
.\end{align}
Because $ \Enu ( f , f ) = \ERnu ( f , f ) $ for $ f \in \dimu \cap \Br $, 
we naturally extend the domain of $ \Enu $ to $ \ldRnuW $. 
By \eqref{:62f}, we see that $ \{ (\lERnu , \ldRnuW ) \}_{\rR \in \N } $ 
is increasing for each $ \nu \in \N $. 
Hence, we set the closed form $ ( \lEInu , \ldWnu ) $ on $ \Lm $ such that 
\begin{align} \label{:62z}&
 \lEInu ( f , f ) = \limi{ \rR } \lERnu ( f , f ) \quad \text{ for } f \in \ldWnu 
 ,\\ & \notag 
 \ldWnu = \{ f \in \bigcap_{ \rR \in \N }\ldRnuW \, ;\, \limi{ \rR } \lERnu ( f , f ) < \infty \} 
.\end{align}
From \eqref{:62f}, \eqref{:62g}, and \eqref{:62z}, we see that 
the sequence of the forms $ ( \lEInu , \ldWnu ) $ is increasing in $ \nuN $. 
Hence, we define the closed form $ (\E _{\infty , \infty} , \ld _{\infty , \infty } )$ as 
\begin{align}\label{:62a}&
\E _{\infty , \infty} (f,f) = \limi{\nu} \lEInu ( f , f ) 
\quad \text{ for } f \in \bigcap_{\nu = 1 }^{\infty} \ldWnu 
, \\ \notag &
\ld _{\infty , \infty } = 
\{ f \in \bigcap_{\nu = 1 }^{\infty} \ldWnu \,;\, \limi{\nu } \lEInu (f,f) < \infty \} 
.\end{align}
\begin{lemma} \label{l:62} 
Assume that $ \mu $ satisfies \QGO. 
Assume that \As{ZC} holds. Then,
\begin{align}\label{:62b} & 
 (\E _{\infty , \infty} , \ld _{\infty , \infty } ) = (\lE , \ld ) 
.\end{align}
\end{lemma}
\PF 
From \eqref{:62g}, we set the closed form $ (\E _{\RI }, \ld _{\RI })$ by 
\begin{align}\label{:62j}&
\E _{\RI } (f,f) = \limi{\nu } \E _{\Rnu } (f,f), 
\\ \notag &
\ld _{\RI } = 
\{ f \in \bigcap_{\nu = 1} ^{\infty} \ld _{\Rnu }\, ;\, \limi{\nu } \E _{\Rnu } (f,f) < \infty \} 
.\end{align}
Obviously, we have 
\begin{align}\label{:62k}&
 (\E _{\RI }, \ld _{\RI }) \le (\lER , \ldR ) \le (\ER , \dom _{\rR })
.\end{align}
It is easy to show that $ (\lER , \ldR ) $ is a quasi-regular Dirichlet form. 
From \eqref{:62k}, we have that $ \underline{\mathrm{Cap}}_{\rR } \le \mathrm{Cap}_{ \rR } $, 
where $ \underline{\mathrm{Cap}}_{\rR }$ is the capacity given by the Dirichlet form 
$ (\lER , \ldR ) $ on $ \Lm $. 
Using \As{ZC}, \eqref{:62c}, and the inequality of the capacities as above, we obtain 
\begin{align}\label{:62l}&
 \underline{\mathrm{Cap}}_{\rR }
 \Big( \sSS \setminus \bigcup_{\nu = 1}^{\infty}\OORnu \Big) = 0 
.\end{align}
Using \eqref{:62j} and \eqref{:62l} and the definition of $ (\lER , \ldR ) $, we deduce 
\begin{align}\label{:62m}&
 (\E _{\RI }, \ld _{\RI }) = (\lER , \ldR ) 
.\end{align}

Using \eqref{:23b}, \eqref{:62i}--\eqref{:62a}, \eqref{:62j}, and \eqref{:62m}, we obtain 
\begin{align} \notag 
(\E _{\infty , \infty} , \ld _{\infty , \infty } ) 
&
= \limi{\nu } ( \lEInu , \ldWnu ) && \text{ by \eqref{:62a}}
\\ \notag &
= \limi{\nu }\limi{\rR } (\E _{\Rnu }, \ld _{\Rnu }) && \text{ by \eqref{:62z}}
\\ \notag &
= \limi{\rR }\limi{\nu } (\E _{\Rnu }, \ld _{\Rnu }) && \text{ by \eqref{:62i}}
\\ \notag &
= \limi{\rR } (\E _{\RI }, \ld _{\RI }) && \text{ by \eqref{:62j}}
\\ \notag &
= \limi{\rR } (\lER , \ldR ) && \text{ by \eqref{:62m}}
\\ \notag &
 = (\lE , \ld ) && \text{ by \eqref{:23b}}
.\end{align}
This completes the proof of \eqref{:62b}. 
\PFEND

\subsection{Upper schemes of Dirichlet forms in infinite volumes}\label{s:62}

We now proceed with the upper scheme. 
We define the domain $ \Brnu $ such that 
\begin{align}\label{:63w} 
&\Brnu = \{f \in\Br ;\, f \text{ is constant on each connected component of }
 \sSS \backslash \OORnu \} 
.\end{align}
From \eqref{:63q} and \eqref{:63w}, we deduce the following for each $ \nuN $: 
\begin{align}\label{:63z}&
\Brnu \subset \mathcal{B}_{ \RRnu } 
 \quad \text{ for all $ \rR \in \N $}
.\end{align}
Using \lref{l:2X}, we see that
$ (\ER ^{\OOnu }, \di ^{\OOnu })$ is closable on 
$ \Lm $ for each $ \Rnu \in \N $. 
Because $ \dimu \cap \Brnu \subset \di ^{\OOnu }$ and 
\begin{align*}&
\text{$ \E (f,f) = \ER ^{\OOnu } (f,f) $\quad for $ f \in \dimu \cap \Brnu $}
,\end{align*}
we have that $ (\E ,\dimu \cap \Brnu ) $ is closable on $ \Lm $. 
Hence, we define $ (\ERnu ,\dRnu ) $ as the closure of $ (\E ,\dimu \cap \Brnu ) $ on $ \Lm $. 
Note that, by construction,
\begin{align}\label{:63a}&
 (\lERnu ,\ldRnuW ) \le (\ERnu ,\dRnu ) 
.\end{align}
\begin{lemma}\label{l:63} 
Assume that $ \mu $ satisfies \QGO. 
Then, $ \{(\ERnu ,\dRnu )\}_{ \rR \in \N }$ is decreasing in $ \rR $ for each $ \nuN $. 
The strong resolvent limit $ (\Enu ,\domK ) $ of $ \{(\ERnu ,\dRnu )\}_{ \rR \in \N }$
 is the closure of $ (\Enu ,\cup_{ \rR } \dRnu ) $. 
\end{lemma}
\PF 
From \eqref{:63z}, we have the following for each $ \nuN $:
\begin{align}& \label{:63f}
\dimu \cap \Brnun \subset \dimu \cap \mathcal{B}_{ \RRnu }
\quad \text{ for all } \rR \in \N 
.\end{align}
Because $ (\ERnu ,\dRnu ) $ is the closure of $ (\ERnu ,\dimu \cap \Brnun ) $, we see that 
$ \{(\ERnu ,\dRnu )\}_{ \rR \in \N }$ is decreasing in $ \rR $ for each $ \nuN $ 
from \eqref{:63f}. This implies the first claim. 

The strong resolvent limit $ (\Enu ,\domK ) $ 
is the closure of the largest closable part of $ (\Enu ,\cup_{ \rR } \dRnu ) $. 
Hence, it only remains to prove that $ (\Enu ,\cup_{ \rR } \dRnu ) $ is closable on $ \Lm $. 

We set $ \1 = \cup_{\rR } (\dimu \cap \Brnu ) $. 
Then, $ \1 \subset \dimu $. Recall that $ (\E ,\dimu ) $ is closable on $ \Lm $. 
Hence, $ (\E , \1 ) $ is closable on $ \Lm $. Note that $ \E (f,f) = \Enu (f,f)$ for $ f \in \1 $. 
Hence, the closability of $ (\Enu ,\1 ) $ on $ \Lm $ follows from that of $ (\E , \1 ) $ on $ \Lm $. 

Note that $ \1 \subset \cup_{ \rR } \dRnu $ and that $ \1 $ is dense in $ \cup_{ \rR } \dRnu $ 
with respect to the inner product 
$ \Enu (\cdot , *) + (\cdot , *)_{\Lm }$. 
Combining these facts with the closability of $ (\Enu ,\1 ) $ on $ \Lm $, 
we deduce that $ (\Enu ,\cup_{ \rR } \dRnu ) $ is closable on $ \Lm $. 
This completes the proof of the second claim. 
\PFEND

\begin{lemma}\label{l:64} 
Consider the same assumptions as for \lref{l:63}. \\
\thetag{1} 
$ \{(\Enu ,\domK )\}_{ \nuN } $ is decreasing in $ \nuN $. For each $ \nuN $, we have that
\begin{align}\label{:64s}&
 (\Enu ,\domK ) \ge (\E ,\dom ) ,\quad 
 \Enu (f , f ) = \E (f , f ) \quad \text{ for all } f \in \domK 
.\end{align}
Furthermore, $ (\E , \cup_{ \nu =1}^{\infty} \domK ) $ is closable on $ \Lm $. 
\\\thetag{2} 
The closure $ (\Edi ,\domI ) $ of $ (\E , \cup_{ \nu =1}^{\infty} \domK ) $ satisfies 
\begin{align} \label{:64w}&
\text{$ (\Edi ,\domI ) = \limi{ \nu } (\Enu ,\domK ) $ in the strong resolvent sense}
,\\ \label{:64x}&
 (\Edi ,\domI ) \ge (\E ,\dom ) 
,\quad 
 \Edi (f,f) = \E (f,f) \quad \text{ for all } f \in \domI 
.\end{align}
\end{lemma}
\PF 
From \eqref{:63r}, we have $ \OORnu \subset \OO _{\Rnu + 1} $. 
Hence, $ \Brnun \subset \Brnunn $ by \eqref{:63w}. 
Then, 
\begin{align}\label{:64y}&
(\E ,\dimu \cap \Brnun ) \ge (\E ,\dimu \cap \Brnunn ) \quad \text{ for each $ \Rnu \in \N $}
.\end{align}
Taking the closures of both sides of \eqref{:64y} and using $ \dom _{ \rR ,\nu + 1 } \subset \dR $, 
we obtain
\begin{align}\label{:64a}& 
 (\ERnu ,\dRnu ) \ge (\ERkK , \dom _{ \rR ,\nu + 1 } ) \ge (\ER ,\dR ) 
\quad \text{ for each $ \Rnu \in \N $}
.\end{align}
By \lref{l:63}, $ (\Enu ,\domK ) $ is the strong resolvent limit of 
$ \{(\ERnu ,\dRnu )\}_{ \rR \in \N }$ 
for each $ \nuN $. From this and the first inequality in \eqref{:64a}, 
$ \{(\Enu ,\domK )\}_{ \nuN } $ is decreasing. 

We see that $ (\E ,\dom ) = \limi{ \rR } (\ER ,\dR ) $ in the strong resolvent sense 
from \lref{l:24} \thetag{1}. 
We can deduce $ (\Enu ,\domK ) =\limi{ \rR }(\ERnu ,\dRnu ) $ in the strong resolvent sense from 
\lref{l:63}. 
Hence, taking $ \rR \to\infty $ in \eqref{:64a}, we obtain the inequality in \eqref{:64s}. 

From $ \E ( f , f ) = \ERnu ( f , f ) $ for $ f \in \dRnu $, we have the equality in \eqref{:64s}. 

From \eqref{:64s}, we have $ (\E ,\dom ) \le (\E , \cup_{ \nu =1}^{\infty} \domK ) $. 
Recall that $ (\E ,\dom ) $ is a closed form on $ \Lm $. 
Hence, $ (\E , \cup_{ \nu =1}^{\infty} \domK ) $ is closable $ \Lm $. 
Thus, we obtain \thetag{1}. 

By \thetag{1}, $ \{ (\Enu ,\domK ) \}_{ \nuN } $ is decreasing. 
Therefore, $ (\Enu ,\domK ) $ converges to 
the closure of the largest closable part of $ (\Edi , \cup_{ \nu =1}^{\infty} \domK ) $ in the strong resolvent sense. 
By \thetag{1}, we see that $ (\Edi , \cup_{ \nu =1}^{\infty} \domK ) $ is closable. 
Hence, we obtain \eqref{:64w}. 
Equation \eqref{:64x} is clear because of \eqref{:64s} and the fact that 
$ (\Edi ,\domI ) $ is the closure of $ (\E , \cup_{ \nu =1}^{\infty} \domK ) $. 
\PFEND 

\begin{proposition}\label{l:65}
Consider the same assumptions as for \lref{l:63}, and furthermore assume that \As{ZC} holds. 
Let $ (\Edi ,\domI ) $ be as in \lref{l:64}. 
Then, $ \domI $ is dense in $ \dom $ with respect to the inner product 
$ \Eone  := \E (\cdot , *) + (\cdot , * )_{ \Lm }$. 
Furthermore, 
\begin{align}\label{:65y}& 
 (\E ,\dom ) = (\Edi ,\domI ) 
,\\\label{:65z}&
\text{$ (\E ,\dom ) = \limi{ \nu } (\Enu ,\domK ) $ in the strong resolvent sense}
.\end{align}
\end{proposition}
\PF 
We regard $ \OORnu ^c $ as a subset of $ \CSR $. Then, 
\begin{align}\label{:65x}&
 \bigcap_{ \3 = 1 }^{\infty} \OORnu ^c = 
 \bigcup_{ m =1}^{ \infty} \{ \sss \in \CSRm \st \sigmaRm (\sss ) = 0 \} 
.\end{align}
From \As{ZC}, we see that
\begin{align}\label{:65w}&
\mathrm{Cap}_{ \rR } ( \bigcup_{ m =1}^{ \infty} \{ \sss \in \CSRm \st \sigmaRm (\sss ) = 0 \} ) 
= 0 
.\end{align}
Using \eqref{:63r}, \eqref{:65x}, and \eqref{:65w}, we obtain
\begin{align} &\notag 
\limi{ \3 } \mathrm{Cap}_{ \rR } ( \OORnu ^c ) = 0 
\quad \text{ for each $ \rR \in \N $}
.\end{align}
Hence, for each $ \rR \in \N $, 
there exist a decreasing sequence of open sets $ \{ \MRnu \}_{\nuN } $ and 
a sequence of functions $ \{ \varphiR \}_{\nuN } $ satisfying 
\begin{align}\label{:65a}&
 \OORnu ^c \subset \MRnu ,\quad 
\limi{ \3 } \mathrm{Cap}_{ \rR } ( \MRnu ) = 0 
,\\\label{:65b} &
\varphiR = 1 \ \text{ on }\MRnu ,\quad 
0 \le \varphiR (\sss ) \le 1 \quad \text{ for all }\sss \in \sSS 
,\\\label{:65d}&
\varphiR \in \dR , \quad \limz{ \3 } \ERone  (\varphiR ) = 0 
.\end{align}
Here, $ \ERone  := \ER (\cdot , *) + (\cdot , * )_{ \LmR }$. 
For $ \qQ \le \rR $, the set $ \OO _{ \Qnu }$ can be regarded as an open set in $ \CSR $, and 
the function $ \varphiQ $ can be considered as an element in $\dR $. 
Recall the definition of $ \dRnu $ given before \eqref{:63a}. 
Then, because of \eqref{:65a}--\eqref{:65d}, we have the following for each $ \rR \in \N $:
\begin{align}\label{:65e}& 
\psiR := \prod_{ \qQ = 1}^{ \rR} ( 1-\varphi _{ \Qnu }) 
 \in \dRnu \cap L^{ \infty }(\mu ) 
.\end{align}
From \eqref{:65b} and \eqref{:65e}, we have $ \psiR = 0 $ on $ \MRnu $. 
Hence, from \eqref{:65b}--\eqref{:65e} with a straightforward calculation, we have that
\begin{align}\label{:65g}&
\limi{ \nu } \Eone  ( 1 - \psiR ) = \limi{ \nu } \ERone  ( 1 - \psiR ) = 0 
\end{align}
and a subsequential limit such that 
\begin{align}& \label{:65h} 
 \limz{ \3 } \psiR (\sss ) = \prod_{\qQ = 1}^{\rR } \limi{ \3 }( 1-\varphi _{ \Qnu } (\sss ) ) = 1 
\quad \text{ for $ \mu $-a.s.\,$ \sss $}
.\end{align}
Combining \eqref{:65g}, \eqref{:65h}, 
$ f \in \dom \cap L^{ \infty }(\mu ) $, and $ | 1 - \psiR | \le 1 $, we obtain 
\begin{align}\label{:65i} & 
\Eone  (f- f\psiR ) =\Eone  (f ( 1 - \psiR )) 
\\ \notag 
 = &\int_{ \sSS } \DDDa [f]| 1 - \psiR |^2 + |f|^2 \DDDa [ 1 - \psiR ] 
+ 2 \DDDa [ f , 1 - \psiR ] + |f|^2 | 1 - \psiR |^2 d\mu 
\\ \notag 
\le& \int_{ \sSS } 2 \left\{ \DDDa [f]| 1 - \psiR |^2 + |f|^2 \DDDa [ 1 - \psiR ] \right\} 
+ |f|^2 | 1 - \psiR |^2 d\mu 
\\ \notag 
\to &\, 0 ,\quad \nu \to \infty 
.\end{align}
For $ f \in \dom \cap L^{ \infty }(\mu ) $, we have that $ f \psiR = 0 $ on $ \MRnu $. 
For $ f \in \dom \cap L^{ \infty }(\mu ) $, this yields 
\begin{align}\label{:65j}& 
 f \psiR \in \dRnu \cap L^{ \infty }(\mu ) \subset \domK \cap L^{ \infty }(\mu ) 
.\end{align}
From \eqref{:65i}, \eqref{:65j}, and 
$ \cup_{ \3 = 1}^{ \infty} \domK \cap L^{ \infty }(\mu ) = \domI \cap L^{ \infty }(\mu ) $, 
we deduce that $ \domI \cap L^{ \infty }(\mu ) $ is dense in $ \dom \cap L^{ \infty }(\mu ) $ with respect to $ \Eone  $. 
Furthermore, it is not difficult to see that $ \dom \cap L^{ \infty }(\mu ) $ is dense in $ \dom $ with respect to $ \Eone  $. 
Collecting these results, we see that $ \domI $ is dense in $ \dom $ with respect to $ \Eone  $, which completes the proof of the first claim. 

The second claim \eqref{:65y} follows immediately from the first claim and \eqref{:64x}. 
The third claim \eqref{:65z} follows from \eqref{:64w} and \eqref{:65y}. 
\PFEND

\subsection{Construction of $ \{ \ORnum \}_{ \Rnu , m \in \N } $}\label{s:63}
In this subsection, we construct $ \{ \ORnum \}_{ \Rnu , m \in \N } $ 
satisfying \eqref{:63q}--\eqref{:63p}. 
\begin{lemma} \label{l:66}
Assume that $ \mu $ satisfies \As{QG} with $ \{ \SRm \}_{ \Rm \in \N } $ and \As{ZC}. 
Assume that $ \sigma _{\rR }^m $ is uniformly continuous on $ \SRm $ for each $ \Rm \in \N $. 
Then, we have a sequence of symmetric open sets $ \{ \ORnum \}_{ \Rnu , m \in \N } $ such that $ \mu $ satisfies \QGO. 
\end{lemma}

\PF 
Let $ (\Phi _0 , \Psi _0 )$ be as in \eqref{:2Xa}. 
Let $  \Rnu , m \in \N $ be fixed. 
Let $ T _1 = \sS _1 $ and $ \TQ = \SQ \setminus \SQQ $ for $ 2 \le \qQ \le \rR $. 
For $ \x = (x_i) \in \SRm $, we set 
\begin{align} &\notag 
\HQd (\x ) = \sum_{ x_i \in \TQ } \Phi _0 ( x_i ) + 
\sum_{x_i ,\, x_j \in \TQ \atop i < j } \Psi _0 ( x_i , x_j )
 + 
\sum_{\qQ '=\qQ +1}^{\rR }
\, \sum_{x_i  \in \TQ ,\,  x_j \in \TQQ } \Psi _0 ( x_i , x_j )
.\end{align}
Let $  \TQl = \{ \x = (x_i) \in  \SRm \, ;\, \sharp \{ i \,;\, x_i \in \TQ \} = l \} $. 
We set $ \AQ ^{0} = \emptyset $. Furthermore, we set for $ \qQ   \in \N $ and $ l \in \N $ 
\begin{align}\label{:66h}& 
\AQ ^{ l } = \{ \x \in \TQl \, ;\, \HQd (\x ) = \infty \}
.\end{align}
Let  $ \mathbf{L}(m) = \{ 
(l_{\qQ }) \in \{ 0 ,\ldots, m \}^{ \rR } \, ;\, 
 \sum_{\qQ =1}^{\rR }l_{\qQ } = m \} $. We set 
\begin{align}\label{:}&
\BRm = \sum _{(l_{\qQ }) \in \mathbf{L}(m) } \sum_{\qQ = 1}^{\rR } \AQl 
.\end{align}
Let $ \HR ^{ \Phi , \Psi } $ be as in \eqref{:2Xp}. We set 
\begin{align} \notag &
 \HR (\x )= \HR ^{ \Phi _0 , \Psi _0 } ( \ulab (\x ) ) 
.\end{align}
Recall that $ ( \Phi _0 , \Psi _0 ) $ is bounded from below because $ \mu $ satisfies 
\As{QG} with $ \{ \SRm \}_{ \Rm \in \N } $. 
Let $ \cref{;2y}(\rR , k , m ,\piRc (\sss )) $ be as in \eqref{:2Xx}. We set 
\begin{align*}&
\Omega (n) = \{ \sss \in \sSS \, ;\, n ^{-1} \le \cref{;2y}(\rR , k , m ,\piRc (\sss )) \le n \}
.\end{align*}
Let $ \sigmaRnm $ be the density function of $ \mu ( \cdot \cap \Omega (n) ) $ on $ \SRm $. 
Suppose that $ \mu ( \Omega (n) ) > 0 $. 
Then, from \eqref{:2Xx} and \eqref{:2Xa}, 
we have a positive constant $ \Ct \label{;69}$ such that
\begin{align}\label{:66i}&
\cref{;69}^{-1} \exp\{ - \cref{;2X} \HR (\x ) \} \le \sigmaRnm (\x ) \le 
\cref{;69} \exp\{ - \cref{;2X}^{-1} \HR (\x ) \}
\quad \text{ for } \x \in \SRm 
.\end{align}
Thus, for $ \x \in \SR ^m $, 
we deduce that $ \sigmaRnm (\x ) = 0 $ if and only if $ \HR (\x ) = \infty $. 
It is clear that $ \mu ( \cup_n \Omega (n) ) = 1 $. 
Hence, $ \{ \x \in \SRm ;\, \sigmaRm (\x ) = 0 \} $ if and only if $ \x \in \BRm $. 

Let $ \enu $, $ \nuN $, be functions defined on $ \N ^2 $ such that,
 for each $ \Ql \in \N $,  
\begin{align}\label{:66d}&
\enu (\Ql ) > \enunu (\Ql ) > 0 , \ \limi{\nu } \enu (\Ql ) = 0 
.\end{align}
We set 
\begin{align} & \notag 
 \AQe  = \{ \x \in \TQl \, ;\, 
\inf _{\y \in  \AQ ^{ l } } | \x - \y | > \enu (\Ql  )
\} 
, \\ \notag &
\ORnum =
\sum _{(l_{\qQ }) \in \mathbf{L}(m) }
\sum_{\qQ = 1}^{\rR } 
\AQee  
.\end{align}
Then, it is easy to show that $ \ORnum $ is a symmetric open set satisfying \eqref{:63q} and \eqref{:63r}. 

Recall that $ \mu ( \piR ^{-1}( \ulab ( \{ \x ; \, \sigmaRm (\x ) = 0 \}) ) ) = 0 $ from \As{ZC}. 
Then, $ \mu ( \piR ^{-1}( \BRm ) )= 0 $ from \eqref{:66h} and \eqref{:66i}. 
Hence, using  \eqref{:66d}, we obtain \eqref{:63s}. 
Equation \eqref{:62c} follows from \eqref{:66h} and \eqref{:66i}. 
Equation \eqref{:63p} is clear because $ \sigmaRm $ is uniformly continuous on 
$ \SRm $ and $\BRm $ is relatively compact in $ \SRm $. 

To prove that $ \mu $ satisfies \QG, we begin by checking \dref{d:22}. 
\dref{d:22} \thetag{1} is obvious. Using \As{ZC}, we have 
\begin{align*}&
\muR (\{ \bigcup_{\nu \in \N }\OORnu \}^c ) = 0 
.\end{align*}
Hence, taking $ \mu _{\Rnu } = \mu (\cdot \cap \OORnu ) $, we obtain \dref{d:22} \thetag{2} 
(replacing $ \mu _{\rR , k}$ by $ \mu _{\Rnu }$). 
By assumption, $ \mu $ satisfies \As{QG} with $ \{ \SRm \}_{ \Rm \in \N } $. 
Hence, applying \eqref{:2Xx} to $ \SRm $ with a simple calculation, we find that 
\begin{align}&\label{:66p}
\cref{;2y}^{-1} e^{-\HR ^{ \Phi ,\Psi } } \LambdaRm (d\mathfrak{x}) \le 
\mukRs ^m (d\mathfrak{x}) \le 
 \cref{;2y}e^{-\HR ^{ \Phi ,\Psi } } \LambdaRm (d\mathfrak{x}) 
.\end{align}
Using \eqref{:66p}, $ \ORnum \subset \SRm $, and retaking $ \ORnum $ if necessary, 
we can easily find that 
 $ \mu $ is a $ (\Phi , \Psi )$-quasi-Gibbs measure with $ \{ \ORnum \}_{ \Rm \in \N } $. 
 Thus, \dref{d:22} \thetag{3} holds. Collecting these results, we obtain \dref{d:22}, 
 which implies \As{QG1}. 
 By assumption, \As{QG2} is satisfied. Hence, we find that $ \mu $ satisfies \QG. 
This completes the proof of \lref{l:66}. 
\PFEND

\section{Proof of \tref{l:32}--\tref{l:34} } \label{s:7} 
In this section, we shall prove \tref{l:32}--\tref{l:34}. 
%

For simplicity, we make the following assumption.\\
\As{B4$^*$} Either \As{B4}, \As{B4$ '$}, or \As{B4$ ''$} holds. 
\smallskip 

\noindent 
We introduce condition \As{B4$_{\nu } $}, which is similar to \As{B4}. 

\smallskip 
\noindent 
\As{B4$_{\nu } $} For each $ \Rnu , m \in \N $, 
\begin{align}\label{:71d}&
\limi{ \nN } \Big\| 
\frac{ \sigma _{ \rR }^{ \Nm } }{ \sigma _{ \rR }^{m} } -1 
 \Big\|_{ \ORnum } 
 = 0 
.\end{align}
\begin{lemma} \label{l:71}
Assume that \As{B4$^*$} is satisfied. Then, \As{B4$_{ \nu } $} holds. 
\end{lemma}
\PF
If \As{B4} holds, then we obviously have \As{B4$_{ \nu }$}. 
From \As{B4$'$}, we have that
\begin{align}\label{:71c}&
\limi{ \nN } \Big\| { \sigma _{ \rR }^{ \Nm } }- { \sigmaRm } \Big\|_{ \ORnum } = 0 
\quad \text{ for each $ \Rnu , m \in \N $}
.\end{align}
From \eqref{:63p}, we see that $ \sigma _{ \rR }^{m}$ is uniformly positive on $ \ORnum $. 
Hence, \eqref{:71c} yields \eqref{:71d}. This implies \As{B4$_{\nu } $}. 
Because \As{B4$''$} implies \As{B4$'$}, \As{B4$''$} yields \As{B4$_{\nu } $}. 
\PFEND

\subsection{Lower schemes of cut-off Dirichlet forms}\label{s:7a}
The main result of this subsection is \pref{l:74}, which presents an inequality for the lower scheme. 

Let $ \OORnum $ and $ \OOnu $ be as in \eqref{:61z}. 
Note that $ \OORnu ^m \cap \OORnu ^{m'} = \emptyset $ for $ m \ne m' $ and that 
 $ \OOnu \subset \OORnu = \cup_{ m = 1 }^{\infty} \OORnu ^m $. 
Replacing $ \mu $ by $ 1_{ \OORnum } \muN $ and $ 1_{ \OOnu } \muN $, 
we introduce the cut-off bilinear forms $ \ERnu ^{ \Nm } $ and $ \ERnu ^{ \nN } $ such that 
\begin{align} \notag 
 \ERnu ^{ \Nm } (f,g) =& \int _{ \OORnum } \DDDaR [f,g] \, d\muN 
,\\ \notag 
\ERnu ^{ \nN } (f,g) = & \sum_{m=1}^{ \infty} \ERnu ^{ \Nm } (f,g) 
 = \int _{ \OOnu } \DDDaR [f,g] d\muN 
.\end{align}

\begin{lemma} \label{l:72X}
Assume that $ \muN $ satisfies \QG. Then, $ (\ERnu ^{ \nN },\dimuN ) $ is closable on $ \LmN $. 
\end{lemma}
\PF 
We obtain \lref{l:72X} from \lref{l:61} by replacing $ \mu $ with $ \muN $. 
\PFEND 
We denote the closure of $ (\ERnu ^{ \nN },\dimuN ) $ on $ \LmN $ as $ (\lERnuN ,\ldRnuW ^{ \nN }) $ . 
 \begin{lemma} \label{l:72} 
Assume that $ \muN $ satisfies \QG. 
Let $ \ldRNW $ be as in \eqref{:31i}. 
Then, for all $ \Rnu , \nN \in \N $, $ (\lERnuN ,\ldRnuW ^{ \nN }) $ satisfies
\begin{align}\label{:72a}& 
 (\lERnuN ,\ldRnuW ^{ \nN }) \le (\lERN , \ldRNW ) 
.\end{align}
\end{lemma}
\PF 
 We have $ (\ERnu ^{ \nN },\dimuN \cap \Br ) \le (\ER ^{ \nN },\dimuN \cap \Br ) $. 
 Hence, $(\lERnuN ,\ldRnuW ^{ \nN }) \le (\lERN , \ldRN ) $. 
Combining this with \eqref{:31j}, we obtain \eqref{:72a}. 
\PFEND 

 Let $ \rN $ be as in \eqref{:52w}. We set 
 \begin{align*}&
\ldW _{ \rN , \nu }^{ \nN } = \{ f \in \ld _{ \rN , \nu }^{ \nN } 
\, ;\, \text{ $ f $ is $ \sigma [\pi _{\rN } ] $-measurable}\} 
. \end{align*}
\begin{lemma}\label{l:73}
Assume that $ \mu $ and $ \muN $ satisfy \QG\ such that \eqref{:63q}--\eqref{:62c} hold. 
Assume that \As{B4$^*$} holds. 
	\XXX 
Then, for each $ \nuN $, 
\begin{align}\label{:73a}
 \lEInu \lA \f \rA \le \liminfi{ \nN }\lErNnuN \lA \fN \rA 
\end{align}
holds for $ ( \lEInu , \ldWnu )$ and 
$ \{ ( \lE _{ \rN , \nu }^{ \nN } , \ldW _{ \rN , \nu }^{ \nN } ) \}_{\nN \in \N }$. 
\end{lemma}
\PF Using \lref{l:71}, we deduce \As{B4$_{\nu } $} from \As{B4$^*$}. 
Then, we obtain \eqref{:73a} from \As{B4$_{\nu } $} in the same fashion as in \pref{l:54}. 
Indeed, we replace $ ( \lE _{ \rN }^{ \nN } , \ldW _{ \rN }^{ \nN } ) $ and 
$ (\E , \dom )$ in \pref{l:54} by $ ( \lE _{ \rN , \nu }^{ \nN } , \ldW _{ \rN , \nu }^{ \nN } ) $ 
and $ ( \lEInu , \ldWnu )$, respectively. 
The remainder of the proof is the same as that of \pref{l:54}, and so we omit the details. 
\PFEND 

\begin{proposition}\label{l:74}
Consider the same assumptions as for \lref{l:66} regarding $ \mu $ and $ \muN $. 
Assume that \As{B4$^*$} holds. 
\XXX 
Then, $ ( \lE , \ld )$ and 
$ \{ ( \lE _{ \rN }^{ \nN } , \ldW _{ \rN }^{ \nN } ) \}_{\nN \in \N } $ satisfy 
\begin{align}\label{:74a}&
 \lE \lA \f \rA \le \liminfi{ \nN } \lErNN \lA \fN \rA 
.\end{align}
\end{proposition}
\PF 
From the definition of $ \E \lA \cdot \rA $, \eqref{:62b}, and \eqref{:62a}, 
we have the following for $ f \in \Lm $
\begin{align}\label{:74d}&
\lE \lA \f \rA = \E _{\infty , \infty} \lA f \rA = \limi{ \nu } \lEInu \lA \f \rA 
.\end{align}
From \eqref{:73a} and \eqref{:72a}, for $ f \in \Lm $, we have that
\begin{align}\label{:74c}& 
 \lEInu \lA \f \rA \le \liminfi{ \nN }\lErNnuN \lA \fN \rA \le \liminfi{ \nN } \lErNN \lA \fN \rA 
.\end{align}
Combining \eqref{:74d} and \eqref{:74c}, we obtain \eqref{:74a}. 
\PFEND 

\subsection{Upper schemes of cut-off Dirichlet forms } \label{s:7b} 
The main result of this section is \pref{l:76}, which presents 
the convergence of the upper scheme of the cut-off Dirichlet forms. 
The argument is similar to that in \sref{s:52}, with the replacement of \eqref{:31k} by \As{B4$_{\nu } $} . 

Let $ \OOnu $ be as in \eqref{:61z}. We set 
\begin{align} &\notag 
\ERnu ^{ \nN } (f,g) = \int _{ \OOnu } \DDDaR [f,g] d\muN 
.\end{align}
\begin{lemma} \label{l:75X}
Assume that $ \muN $ satisfies \QG. 
Then, $ (\ERnu ^{ \nN },\dimu \cap \Br ) $ is closable on $ \LmN $.
\end{lemma}
\PF 
Using \As{B2}, we have $ \dimu \subset \cap_{ \nN \in \N }\dimuN $. 
Then, 
\begin{align}\label{:75Xa}&
\dimu \cap \Br \subset \dimu \subset \cap_{ \nN \in \N }\dimuN \subset \dimuN 
.\end{align}
From \lref{l:72X}, we see that $ (\ERnu ^{ \nN },\dimuN ) $ is closable on $ \LmN $. 
Combining this with \eqref{:75Xa} completes the proof. 
\PFEND

We denote the closure of $ (\ERnu ^{ \nN },\dimu \cap \Br ) $ on $ \LmN $ as $ ( \ERnu ^{ \nN }, \dRnu ^{ \nN } ) $ . 
We set 
\begin{align} &\notag 
\Enu ^{ \nN } (f,g) = \int _{ \OOnu } \DDDa [f,g] d\muN 
.\end{align}
Similar to \lref{l:23}, we see that 
$ (\Enu ^{ \nN }, \cup_{\rR = 1}^{\infty} \dimu \cap \Br )$ is closable on 
$ \LmN $. 
Let $ (\Enu ^{ \nN } , \domKN )$ be the closure of 
$ (\Enu ^{ \nN }, \cup_{\rR = 1}^{\infty} \dimu \cap \Br )$ on $ \LmN $. 
Then, similar to \lref{l:24}, we see that $ \{ ( \ERnu ^{ \nN }, \dRnu ^{ \nN } ) \} $ converges to 
$ (\Enu ^{ \nN } , \domKN )$ on $ \LmN $ as $ \rR \to \infty $ in the strong resolvent sense. 
Let $ (\Enu ,\domK ) $ be as in \lref{l:63}. 
\begin{lemma} \label{l:75} 
Assume that \As{B4$^*$} holds. 
Then, $ (\Enu ^{ \nN } , \domKN )$ and $ (\Enu ,\domK ) $ satisfy the following. 
For each $ \f \in \Lm $, there exists a sequence $ \{ \gNnu \}_{\nN \in \N } $ satisfying 
\begin{align}
 \label{:75a} & 
\text{$ \gNnu $ is $ \sigma [ \pi _{\rN }]$-measurable, 
$ \gNnu \in \LmrNN $
}
,\\
&\label{:75b}
\limi{ \nN }\gNnu = \f \text{ strongly in the sense of \dref{d:42}}
,\\&\label{:75c}
\limi{N}\Enu ^{ \nN }\lA \gNnu \rA = \Enu \lA \f \rA 
.\end{align}
\end{lemma}
\PF 
From \lref{l:71}, we deduce \As{B4$_{\nu } $} from \As{B4$^*$}. 
We then obtain \eqref{:75a}--\eqref{:75c} 
from \As{B4$_{\nu } $} in the same fashion as in \pref{l:56}. 
We omit the details of the proof. 
\PFEND

Next, we check \dref{d:44} (2).
\begin{proposition}\label{l:76}
Consider the same assumptions as for \lref{l:66}. Assume that \As{B4$^*$} holds. 
Then, $ (\EN , \dN )$ and $ (\E , \dom )$ satisfy the following. 
For any $ \f \in \Lm $, there exists a sequence $ \{ \gN \}_{ \nN \in \N } $ such that 
\begin{align} \label{:76x} & 
\text{$ \gN $ is $ \sigma [ \pi _{\rN }]$-measurable, 
$ \gN \in \LmrNN $}
,\\\label{:76y} &
\limi{ \nN }\gN = \f \text{ strongly in the sense of \dref{d:42}}
,\\\label{:76z} & 
\limi{ \nN } \E ^{ \nN }\lA \gN \rA = \E \lA \f \rA 
.\end{align}
\end{proposition}

\PF 
From \lref{l:66}, we see that the assumptions of \pref{l:65} are fulfilled. 
Using \pref{l:65}, we find that $ \cup_{ \3 = 1}^{ \infty} \domKN $ is dense in $ \dom $. 
Hence, without loss of generality, we can assume that $ \f \in \cup_{ \3 = 1}^{ \infty} \domKN $. 
From \pref{l:65}, we obtain 
\begin{align}\label{:76d}
\E ( \f ) &=\limi{ \nu }\Enu ( \f ) 
\end{align}
Let $ \gNnu $ denote the sequence in \lref{l:75}. Then, $ \gNnu $ satisfies \eqref{:76x}. 
Combining \eqref{:75a}--\eqref{:75c} and \eqref{:76d}, 
we can take $ \gN $ satisfying \eqref{:76y} and \eqref{:76z} 
by choosing a subsequence of $ \{ \gNnu \}_{ \nN , \nu \in \N } $. 
\PFEND

\subsection{Proof of \tref{l:32}--\tref{l:34} } \label{s:7c}
\noindent {\em Proof of \tref{l:32}. } 
We first assume that \As{B4$''$} holds. 
\XXX 
Let $ (\lErNN , \ldW _{ \rN }^{ \nN } ) $ be as in \eqref{:31j} with $ \rR = \rN $. 

From \pref{l:74}, we see that
$ ( \lE , \ld ) $ and $ \{ ( \lE _{ \rN }^{ \nN } , \ldW _{ \rN }^{ \nN } ) \}_{\nN \in \N } $ 
 satisfy 
\begin{align}\label{:7Xz}&
 \lE \lA \f \rA \le \liminfi{ \nN } \lErNN \lA \fN \rA 
.\end{align}
From \As{A2}, we have $ (\E ,\dom ) = (\lE ,\ld )$. Combining these results, we see that
$ ( \E , \dom ) $ and $ \{ ( \lE _{ \rN }^{ \nN } , \ldW _{ \rN }^{ \nN } ) \}_{\nN \in \N } $ 
satisfy \dref{d:44} \thetag{1}. 

Let $ \gN $ be as in \pref{l:76}. Applying \eqref{:7Xz} to $ \{ \gN \} $, we have that
\begin{align}&\label{:7Xb}
 \E \lA \f \rA \le \liminfi{ \nN } \E _{ \rN }^{ \nN } \lA \gN \rA 
.\end{align}
Clearly, $ \E _{ \rN }^{ \nN } \lA \gN \rA \le \E ^{ \nN } \lA \gN \rA $. 
Hence, from \pref{l:76}, we obtain 
\begin{align}\label{:7Xc}&
 \limsupi{ \nN } \E _{ \rN }^{ \nN } \lA \gN \rA 
 \le \limsupi{ \nN } \EN \lA \gN \rA = \E \lA \f \rA 
.\end{align}
Combining \eqref{:7Xb} and \eqref{:7Xc}, we have 
\begin{align}\label{:7Xd}&
 \E \lA \f \rA = \limi{ \nN } \0 \lA \gN \rA 
.\end{align}
Hence, we find that
$ ( \E , \dom ) $ and $ \{ ( \lE _{ \rN }^{ \nN } , \ldW _{ \rN }^{ \nN } ) \}_{\nN \in \N } $ 
 satisfy \dref{d:44} \thetag{2}. 

From \eqref{:7Xz} and \eqref{:7Xd}, we see that 
the Mosco convergence of $ ( \0 , \ldW _{ \rN }^{ \nN } ) $ on 
$ \LmrNN $ to $ ( \E , \dom ) $ on $ \Lm $ holds. 
Combining this with \lref{l:41}, we conclude that \eqref{:31r} holds. 

Next, we suppose that \As{B4$''$} holds. 
Combining \eqref{:33s} and \eqref{:33t}, we see that
\begin{align}\label{:7Ya}&
\sup_{ \xm \in \SRm }\rho ^m (\xm )\le \cref{;21a}^m m^{ \cref{;21b}m}
.\end{align}
Then, from \eqref{:33t} and \eqref{:7Ya}, we obtain 
\begin{align}\label{:7Yb} 
\sigma _{ \rR }^{ \nN ,m} (\xm ) & = \sum_{ n =0}^{ \infty} \frac{(-1)^{ n }}{ n !} 
\int_{ \SR ^{ n } } \rho ^{ \nN , m + n }(\xm ,\yn ) d\yn 
,\\
\label{:7Yc} 
\sigmaRm (\xm ) &= \sum_{ n =0}^{ \infty} \frac{(-1)^{ n }}{ n !} 
\int_{ \SR ^{ n } } \rho ^{ m + n }(\xm ,\yn ) d\yn 
.\end{align}
Combining \eqref{:33s}, \eqref{:7Yb}, and \eqref{:7Yc} and 
using the Lebesgue convergence theorem, we obtain \eqref{:32y}, which implies that 
 \As{B4$'$} holds. 
Hence, \eqref{:31r} follows from the first part of the proof. 
\PFEnd 

\noindent {\em Proof of \tref{l:34}. }
We first check the Mosco convergence of $ \4 $ on $ \LmrNN $. 
\XXX 
From \eqref{:7Xz} and the first inequality in \eqref{:34a}, 
$ \4 $ satisfies 
\begin{align} \label{:7Za}&
 \lE \lA \f \rA \le \liminfi{ \nN } \lE _{ \rN }^{ \nN } \lA \fN \rA 
.\end{align}
Combining this with \As{A2} implies \dref{d:44} \thetag{1}. 

Let $ \gN $ be as in \pref{l:76}. 
Then, from \eqref{:7Xc} and the second inequality in \eqref{:34a}, $ \4 $ satisfies 
\begin{align}\label{:7Zb}&
 \limsupi{ \nN } \E _{ \rN }^{ \nN } \lA \gN \rA \le \E \lA \f \rA 
.\end{align}
Taking $ \fN = \gN $ in \eqref{:7Za} and combining this with \eqref{:7Zb} and \As{A2}, 
we obtain \dref{d:44} \thetag{2}. 

Thus, the Mosco convergence of $ \4 $ on $ \LmrNN $ to $ ( \E , \dom ) $ on $ \Lm $ holds. 
Using this and \lref{l:41} completes the proof of \tref{l:34}. 
\PFEnd

\section{Proof of \tref{l:35}--\tref{l:39}} \label{s:8}
Let $\XB $ be the unlabeled diffusion given by the Dirichlet form 
$ (\lE _{ \rN }^{ \nN } , \widetilde{\ld }_{ \rN}^{ \nN } ) $ as in \tref{l:31}. 
Originally, $ \XXrNN $ was a $ \CSrNbar $-valued process. 
We regard $ \XXrNN $ as a $ \WSs $-valued process in an obvious manner. 

Let $ (\XX , \PP ) $, $ \PP = \{ \PP _{\sss } \}_{\sss \in \sSS } $, 
 be the $ \sSS $-valued, $ \mu $-reversible diffusion 
associated with $ (\E ,\dom ) $ on $ \Lm $ (see \lref{l:25}). 
Let $ \xi $ be as in \As{B3}. 
Let $ \PP _{ \xi d\mu } = \int _{\sSS } \PP _{\sss } \xi (\sss ) d\mu $. 

Recall that the diffusion processes $ \XXrNN $ and $ \XX $ are given by 
$ \XXrNN (t) = \ww (t)$ and $ \XX (t) = \ww (t)$ as functions defined on $ \WSs $, 
where $ \SSs $ is defined by \eqref{:26g}. 
Thus, we write 
$ \lpathN (\XXrNN (\ww ) ) = \lpathN (\ww ) $ and $ \lpath (\XX (\ww )) = \lpath (\ww )$. 

We consider the discontinuity set of the sequence 
$ \{\lpathN (\ww ) (t)\}_{\nN \in \N }$ converging to $ \lpath (\ww ) (t) $ such that 
\begin{align*}
\mathrm{Disc}[ \lpath (\ww )& ( t ) ]= 
\{ \ww\in \WSS \, ;\, \text{ $ \ww \notin \WSs $ or 
there exists $ \{\ww _{\nN } \}_{\nN \in \N }$
 in $ \WSs $ }
 \\ \notag & \text{such that }
 \limi{\nN } \ww _{\nN } = \ww \text{ and }
 \limi{\nN } \lpathN (\ww ) (t) \ne \lpath (\ww ) (t) 
\} 
.\end{align*}
We set $ \mathrm{Disc}[ \lpath (\ww )( u ) - \lpath (\ww ) ( t ) ]$ similarly. 
\begin{lemma} \label{l:81}
For each $ 0 \le t < u < \infty $, 
\begin{align} \label{:81c}&
 \PP _{ \xi d\mu } ( \mathrm{Disc}[ \lpath (\ww) ( t ) ] ) = 0 
,\\ \label{:81a}& 
 \PP _{ \xi d\mu } ( \mathrm{Disc}[ \lpath (\ww )( u ) - \lpath (\ww ) ( t ) ] ) = 0 
.\end{align}
\end{lemma}
\PF 

Let $ \WSsiNE $ and $ \WSsNE $ be as in \eqref{:26i}. 
By definition, $ \WSsiNE \subset \WSsNE \subset \WSs $. 
Recall that $ \Pmu ( \WSsiNE ) = 1 $ by \As{A4}. 

Let $ t $ be fixed. 
It is easy to show that $ \lpath (\ww) ( t ) $ restricted on $ \WSsNE $ 
is a continuous function in $ ( \lab (\ww (0)) , \ww ) $ in the sense that if $ (\lab (\ww _N (0)) , \ww _N )$ converge to $ ( \lab (\ww (0)) , \ww ) $, then 
$ \lpathN (\ww _N ) ( t ) $ converge to $ \lpath (\ww) ( t ) $. 
Hence, from \As{C1}, we have 
\begin{align}& \label{:81b}
\PP _{ \mu } ( \mathrm{Disc}[ \lpath (\ww ) ( t ) ] ) = 0 
.\end{align}
Because $ \PP _{ \xi d\mu }$ is absolutely continuous with respect to 
$ \PP _{ \mu }$, we deduce \eqref{:81c} from \eqref{:81b}. 

From \eqref{:81c} with a simple calculation, we obtain \eqref{:81a}. 
\PFEND

\noindent{\em Proof of \tref{l:35}. }
From \As{A4}, \As{C2}, and \eqref{:35a}, we can construct the labeled processes 
$ \XrNN = ( \lpath ( \XXrNN ), o , o ,\ldots ) $ and $ \X = \lpath (\XX ) $. 
Note that the initial distribution of $ \XXrNN $ has a density in $ \LmrNN $ from \As{B3}. 
Hence, it is sufficient for the tightness to prove the case in which $ \XXrNN $ start from the stationary distribution $ \murNN $. Let 
\begin{align*}& 
\mu _{\rN }^{ \nN , m } = \muN \circ \pi _{\SO _{\rN }}^{-1}
 (\, \cdot \, | \sss (\SO _{\rN }) \ge m ) 
.\end{align*}
By construction, 
$ \mu _{\rN }^{ \nN , m } (\{ \sss ; \, m \le \sss (\SO _{\rN }) < \infty \} ) = 1 $. 

We assume that $ \XXrNN (0) \elaw \mu _{\rN }^{ \nN , m } $ in the rest of the proof. 
We write $ \XrNN = (X_{\rN }^{\nN , i })_{i=1}^{\infty}$. 
To apply the Lyons-Zheng decomposition to $X_{\rN }^{\nN , i }$, $ 1 \le i \le m $, 
we use the $ m $-labeled process such that 
\begin{align*}&
\XrN ^{\nN ,[ m ]} = ( (X _{\rN }^{\nN , i } )_{i=1}^{ m }, 
\sum_{i = m + 1 }^{\nnnN } \delta_{ X _{\rN }^{\nN , i } } ) 
.\end{align*}
Then, $ \XrN ^{\nN ,[ m ]} $ is the diffusion process associated with the Dirichlet form 
$(\E _{\rN }^{\nN , [ m ]}, \ldrNNm ) $ on $ \2 $. 
Here, $ \murNNm $ is the $ m $-Campbell measure of 
$ \mu _{\rN }^{ \nN , m } $ and $ \E _{\rN }^{\nN , [ m ]} $ is the Dirichlet form such that 
\begin{align} \notag & 
\E _{\rN }^{\nN , [ m ]} (f,g) = 
\int_{\SrNSS } \DDD ^{a, [ m ]} [f,g] d \murNNm 
.\end{align}
Furthermore, $ \DDD ^{a, [ m ]}$ is the carr\'{e} du champ on $ \SrNSS $ such that 
\begin{align*}&
 \DDD ^{a, [ m ]} [f,g] (\mathbf{x},\sss ) = \half \sum_{i=1}^m 
(a (x_i, \{\sum _{j \not=i }^{ m } \delta_{x_j} \} + \sss ) 
 \nabla_{x_i} f (\mathbf{x},\sss ) , \nabla _{x_i} g (\mathbf{x},\sss ) )_{ \Rd } 
 + \DDDa [ f , g ] (\mathbf{x},\sss )
,\end{align*}
where $ \mathbf{x} = (x_1,\ldots, x_m) \in \Sm $, and we regard $ \DDDa $ as 
 the carr\'{e} du champ on $ \SrNSS $ in an obvious fashion. 
The domain $ \ldrNNm $ is taken to be the closure of 
\begin{align*}&
\{ f \in C_0^{\infty}(\Sm ) \otimes \di ;\, \E _{\rN }^{\nN , [ m ]} ( f , f ) < \infty , f \in 
 \2 \} 
.\end{align*}

We see that $ X _{ \rN }^{ \nN ,i} $, $ 1\le i \le m $, 
 is an additive functional of the $ \murNNm $-symmetric, conservative diffusion 
 $ \XrN ^{\nN ,[ m ]} $. 
Moreover, $ X _{ \rN }^{ \nN ,i} $ is a Dirichlet process of $ \XrN ^{\nN ,[ m ]} $. 
 Here, a Dirichlet process is an additive functional of a Markov process associated with a Dirichlet form given by the composition of the Markov process with a function belonging to the domain of the Dirichlet form locally. 
Thus, we can apply the Lyons-Zheng decomposition to $ X _{ \rN }^{ \nN ,i} $. 
Note that we cannot apply the Lyons-Zheng decomposition to $ X _{ \rN }^{ \nN ,i} $ 
 as an additive functional of the unlabeled diffusion $ \XXrNN $ directly, because 
$ X _{ \rN }^{ \nN ,i} $ is not a Dirichlet process of $ \XXrNN $. 
See Section 9 in \cite{k-o-t.ifc} for the proof of the Lyons-Zheng decomposition.

Let $M^{[x_i ]} $ be a continuous martingale additive functional of $\XrNNm $ such that 
\begin{align}\label{:82b}&
M_t^{[x_i ]} \big(\XrNNm \big) = 
 \int_0^t \sigma (X _{ \rN }^{ \nN , i }(u), \XX _{ \rN }^{ \nN , i\dia }(u)) 
 dB^{ \nN , i }(u)
.\end{align}
For $ T > 0 $, we set $\mathcal{\rR }_T (\mathbf{w}) (t) := \mathbf{w} (T-t) $. 
For each $ 0 \le t \le T $ and $ 1 \le i \le m $, we set 
\begin{align}\label{:82bb}&
M_t^i = M_t^{[x_i ]} \big(\XrNNm \big) 
, \\ \notag &
M_t^{*i} = M_{T-t}^{[x_i ]}
 \big(\mathcal{\rR }_T (\XrNNm ) \big) - 
 M_{T}^{[x_i ]} \big(\mathcal{\rR }_T (\XrNNm ) \big) 
.\end{align}
Using the Lyons-Zheng decomposition for solutions of SDE \eqref{:36b} with the function $ x_i $, we have that, for each $ 0 \le t \le T $ and $ 1 \le i \le m $, that 
\begin{align}\label{:82a}&
 X _{ \rN } ^{ \nN , i} (t) - X _{ \rN }^{ \nN ,i} (0) = 
 \half \big \{ M_t^i + M_t^{*i} 
\big\}
\end{align}
and that under $ \PPxs $ for $ \murNNm $-a.e.\,$ \xs $, $ M^i $ and $ M^{*i}$ are continuous martingales such that $ M_0^i = M_0^{*i}=0$. 
Here, $ \PPxs $ is the distribution of the diffusion process 
 $ \XrN ^{\nN ,[ m ]} $ associated with the Dirichlet form 
$(\E _{\rN }^{\nN , [ m ]}, \ldrNNm ) $ on 
$ \2 $ starting at $ (x,\mathsf{s})$.

From \eqref{:21c} and \eqref{:82b}--\eqref{:82a}, 
there exists a constant $ \Ct \label{;82}$ independent of $ i $ such that 
\begin{align}\label{:82B}&
E[ | X _{ \rN } ^{ \nN , i} (t) - X _{ \rN }^{ \nN ,i} (u) |^4 ] 
\le \cref{;82} | t - u |^2 \ \text{ for all } 0 \le t,u \le T 
.\end{align}

Using \eqref{:35a}, we see that $ \{ X _{ \rN } ^{ \nN , i} (0) \}_{ \nN \in \N } $ is tight in $ \sS $ 
 for each $ 1 \le i \le m $. 
Combining this with \eqref{:82B}, we easily obtain the tightness of 
 $ \{ X _{ \rN } ^{ \nN , i} \}_{ \nN \in \N } $ in $ C([0,T]; \sS ) $ for each $ 1 \le i \le m $. 
Because $ T $ is arbitrary, this implies the tightness of 
$ \{ X _{ \rN } ^{ \nN , i} \}_{ \nN \in \N } $ in $ C([0,\infty) ; \sS ) $ for each $ 1 \le i \le m $. 

 Taking an arbitrary $ m \in \N $, we obtain the tightness of 
 $ \{ X _{ \rN } ^{ \nN , i} \}_{ \nN \in \N } $ in $ C([0,\infty) ; \sS ) $ for all $ i \in \N $. 
From this, we deduce the tightness of 
$ \XrNN = (X_{\rN }^{\nN , i })_{i=1}^{\infty}$ in $ C([0,\infty) ; \SN ) $. 
Here, we endow $ C([0,\infty) ; \SN ) = \prod_{i\in\N} C([0,\infty) ; \sS ) $ with the product topology. 
We use the fact that, in general, the tightness of random variables with the value of a countable product of Polish spaces follows from that of each component-wise random variable. 

Recall that $ \XrNN = ( \lpath (\XXrNN ) , o , o ,\ldots )$ and $ \X = \lpath (\XX )$. 
From the tightness of $ \XXrNN (0)$, \As{C3}, and \eqref{:82B}, it easily follows that 
$ \XXrNN $ is tight in $ \WSS $. 

From \tref{l:31}, \lref{l:81}, \As{LIN}, and the tightness of $ \XXrNN $ in $ \WSS $, 
we see that the random variables $ \XrNN ( t ) $ and $ \XrNN ( u ) - \XrNN ( t ) $ 
converge weakly to $ \X ( t ) $ and $ \X ( u ) - \X ( t ) $, respectively. 
Thus, we have the convergence of the finite-dimensional distributions of 
$ \XrNN $ to those of $ \X $. 
From this, we deduce the convergence of the finite-dimensional distributions of 
$ \XrN ^{ \nN ,m} $ to $ \X ^{ m } $. 
Collecting these results, we obtain \tref{l:35}. 
\PFEnd

\noindent{\em Proof of \tref{l:39}. } 
From the proof of \tref{l:35}, we see that $ \{ \XXrNN \} $ is tight in $ \WSS $. 
Hence, \tref{l:39} follows from \tref{l:35}. 
\PFEnd

\section{A sufficient condition for \As{C3}} \label{s:X}

In this section, we present a sufficient condition for \As{C3} 
in terms of initial distributions. 
Let $ \lab ^{\nN }$ and $ \murNN $ be as in \sref{s:32}. 
Let $ \mathscr{R}(t) = \int_t^{\infty} (1/\sqrt{2\pi }) e^{-|x|^2/2} dx $ 
be a (scaled) complementary error function. 
We specify the following condition on the initial distributions 
$ \murNN \circ (\labN )^{-1} $, $ \nN \in \N $. 

\noindent \As{D} For any positive numbers $ \rR $ and $ T $, 
\begin{align}\label{:X1x} & \quad 
\limi{l} \Big\{
 \sup_{\nN \in \N } \int \sum_{i \ge l }\mathscr{R}( \frac{|\si |- \rR }{T}) 
 \murNN \circ (\labN )^{-1} (d\mathbf{s}) \Big\} 
= 0 
.\end{align}
\begin{lemma} \label{l:44} 
Assume that \As{D} holds. Then, $ \XB $ satisfies \As{C3}. 
\end{lemma}
\PF 
Let $ \XB $ be the diffusion defined in \sref{s:32}. 
Then, $ \XB $ is the diffusion associated with 
 $ (\lE _{ \rN }^{ \nN } , \widetilde{\ld }_{ \rN}^{ \nN } ) $ on $ \LmrNN $. 
We denote the distribution of the diffusion with the initial distribution $ \murNN $ using the same symbol $ \PPrNN $.

From \eqref{:29r}, we have that
\begin{align}\label{:X1y}
 \PPrNN \Big( \mr ( \lpathN ( \XXrNN ) ) > l \Big) = & 
 \Pmgone \Big( \bigcup_{i\ge l} \{ \inf_{t\in[0,T]}|\Xti | \le\rR \} \Big)
 \\ \notag 
 \le & \sum_{i \ge l} \Pmgone \Big( \inf_{t\in[0,T]}|\Xti | \le\rR \Big)
. \end{align} 

Let 
$(\E _{\rN }^{\nN , [ m ]}, \ldrNNm ) $ be the Dirichlet form on $ \2 $ 
given in the proof of \tref{l:35} in \sref{s:8}. 
Let $ \PPxs $ be the distribution of the associated diffusion process starting at $ \xs $, as before. 

From \eqref{:82a}, we see that, under $ \PPxs $ for $ \murNNm $-a.e.\,$ \xs $, 
the stochastic processes 
$ M^i $ and $ M^{*i}$ are continuous martingales with $ M_0^i = M_0^{*i}=0$ satisfying 
\begin{align}\label{:X1a}&
 X _{ \rN } ^{ \nN , i} (t) - X _{ \rN }^{ \nN ,i} (0) = 
 \half \big \{ M_t^i + M_t^{*i} 
\big\}
.\end{align}
We set $ X^i = X _{ \rN } ^{ \nN , i} $, $ \Xm = (X^1,\ldots,X^m)$, and 
$ \XXrNNms = \sum_{i > m }^{\nnnN } \delta_{X^i }$. 
 Here, $ \nnnN $ is the number of particles in $ \SO _{\rN } $ such that 
 $ \nnnN = \XXrNN (\SO _{\rN } )(0)$. 
By construction, $ \Xm = \labNm (\XXrNN ) $. 
Using \thetag{2.18} of Theorem 2.4 in \cite{o.tp}, we have the identity 
\begin{align}\label{:X1b}&
\PPxs = \PPrNN ( (\Xm , \XXrNNms ) \in \,\cdot\, | (\Xm , \XXrNNms ) (0) = \xs 
 )
.\end{align}
Using \eqref{:X1a} and \eqref{:X1b}, we see that
 $ X^i = X _{ \rN } ^{ \nN , i} $ satisfies \eqref{:X1a} under $ \PPrNN $. 
 Hence, we have that
\begin{align} \label{:X1d}
 & \Pmgone ( \inf_{t\in[0,T]}|\Xti | \le\rR ) \le \, 
\Pmgone ( \sup_{t\in[0,T]}|\Xti -\Xzi | \ge |\Xzi | -\rR ) 
\\ \notag \le \, &
\Pmgone ( \sup_{t\in[0,T]} | \Mti | \ge |\Xzi | -\rR ) + 
\Pmgone ( \sup_{t\in[0,T]} | M_t^{*i} | \ge |\Xzi | -\rR ) 
\quad \text{by \eqref{:X1a}, \eqref{:X1b}}
\\ \notag = \, & 2 \, 
 \Pmgone ( \sup_{t\in[0,T]} |\Mti | \ge |\Xzi | -\rR ) 
.\end{align}
Let $ \Emg $ denote the expectation with respect to $ \Pmg $. 
Then, using the martingale inequality, 
we see that there exists a positive constant $ \Ct \label{;44}$ such that 
\begin{align}\label{:X1e}
 \Pmgone ( \sup_{t\in[0,T]} |\Mti | \ge |\Xzi | -\rR ) 
\le & 
 \Emg ( 
\mathscr{R}( \frac{|\Xzi |- \rR }{{\cref{;44}}T} ) 
) 
\\ \notag 
= & 
\int \mathscr{R}( \frac{|\si |- \rR }{ \cref{;44} T }) \murNN \circ (\labN )^{-1} (d\mathbf{s}) 
.\end{align}

Combining \eqref{:X1x}, \eqref{:X1y}, \eqref{:X1d}, and \eqref{:X1e}, we can deduce that
\begin{align}& \notag 
\sup_{\nN \in \N } \PPrNN ( \mr ( \lpathN ( \XXrNN ) ) > l ) \le 
 2 \sup_{\nN \in \N } \sum_{i \ge l} 
\int_{\sSS }
 \mathscr{R}( \frac{|\si |- \rR }{{\cref{;44}}T} ) 
 \murNN \circ (\labN )^{-1} (d\mathbf{s})
.\end{align}
Hence, \As{D} implies \As{C3}. This completes the proof. 
\PFEND

\section{Examples of dynamical universality}\label{s:9}
In this section, we give some examples of dynamical universality. 
We consider the sine$ _\beta $ random point field, $ \beta = 1,2,4 $, 
and the Ginibre random point field. 
All examples satisfy the assumptions in \tref{l:31}--\tref{l:39}, 
and the main theorems are thus applicable to these examples. 

For these random point fields, \As{A1} is proved in \cite{o.rm}. 
We have \As{A2} from \cite{k-o-t.udf}. Assumption \As{A3} obviously holds. 
We obtain \As{A4} from \cite{o-t.tail} based on the result in \cite{o.col}. 
We have \As{A5} from \cite{o.isde}. 
We check \As{A6} in \cite{k-o-t.udf}. 
Both \As{ZC} and \As{C2} hold according to \cite{o.col}. In \cite{o.col}, these are only proved for $ \rR = \infty $. The current case can be proved in a similar fashion. 
The quasi-Gibbs property is checked in \cite{o.rm}. 
Assumptions \As{B1} and \As{B2} are clear because the random point fields in these conditions 
are supported on finite particle systems. 
Condition \eqref{:33t} in \As{B4$''$} holds because the correlation functions come from determinants of matrices given by kernels that are uniformly bounded on $ \SR \ts \SR $ for each $ \rR \in \N $ (see, for example, \cite[Lemma 10.1]{o.rm} for the Ginibre random point field). 
As for \As{C3}, we present a sufficient condition in \sref{s:X}, which is satisfied by all the examples. Condition \As{C4} obviously holds. 

\subsection{The sine$_\beta $ interacting Brownian motion with $ \beta = 1,2,4$}\label{s:91}
Let $ \muVtN $ be the random point field whose density is $ \mVtN $ given by \eqref{:11e}. 
Then, the logarithmic derivative $ \dmuN _{ \rN }$ of $ \muVtN $ on $ \SrN $ is given by 
\begin{align}\label{:91a}
 \dmuN _{ \rN } \xs = & - \frac{1}{\varrho_{V}(\theta )} 
\Vb ' \Big(\frac{x}{N\varrho_{V}(\theta )} + \theta \Big) 
\\ \notag &
+ \frac{\beta }{2} 
 \sum_{\si \in \sS _{ \rN } } \frac{1}{ x - \si } 
+ \frac{\beta }{2}
 \int_{ \sS \backslash \sS _{ \rN } } \frac{1}{ x - y } \rNone ( y ) dy 
.\end{align}
From \eqref{:36b}, \eqref{:36d}, and \eqref{:91a}, we see that the associated SDE is given by 
\begin{align}\label{:91b}
 X_{ \rN }^{ \nN , i } (t) & - X_{ \rN }^{ \nN , i } (0) = 
B ^{i} (t) - 
 \int_0^t \frac{1}{\varrho_{V}(\theta )} 
\Vb ' \Big(\frac{X_{ \rN }^{ \nN , i } (u) }{N\varrho_{V}(\theta )} + \theta \Big) du 
\\ \notag & + \frac{\beta }{2}
 \int_0^t \sumN \frac{1}{ X_{ \rN }^{ \nN , i } (u) - X_{ \rN }^{ \nN , j } (u) } du 
 + 
\frac{\beta }{2}
 \int_{ \sS \backslash \sS _{ \rN } } \frac{1}{ X_{ \rN }^{ \nN , i } (t) - y } \rNone ( y ) dy 
\\ \notag & + 
 \half \int_0^t \mathbf{n}^{ \rN }( X_{ \rN }^{ \nN , i } ( u ) ) 
L_{ \rN }^{ \nN , i } (d u ) 
.\end{align} 
For $ \beta = 1,2,4$ and $ \V $ as in \eqref{:11q}, condition 
\eqref{:33s} in \As{B4$ ''$} is proved in \cite{dkmvz}. 
For $ \beta = 2 $ with a real analytic function $ \V $ satisfying \eqref{:11Q}, condition \eqref{:33s} in \As{B4$ ''$} is proved in \cite{dkmvz.2}. 
 Hence, we apply the results in \tref{l:31}--\tref{l:39} to these models. 

\begin{proposition}	\label{l:91}
The following equation holds. In particular, \eqref{:13f} holds. 
\begin{align}\notag 
\limi{ \nN } \pP \Big( \max_{0 \le t \le T } 
 \Big\{
 X^{ \nN , i } (t) & - X^{ \nN , i }(0) - \frac{\beta }{2}
 \int_0^t \sum_{ j\neq i }^{ \nnnN} \frac{1}{X^{ \nN , i } (u) - X^{ \nN , j } (u) } du 
 \\ \notag &
 + \cref{;91} (\beta ) t 
 - \frac{\beta }{2}
 \int_{ \sS \backslash \sS _{ \rN } } \frac{1}{ X^{ \nN , i } (t) - y } \rNone ( y ) dy 
 \Big\} \ge a 
\Big)
\\ \notag &
= 2 \int_{ x \ge a }\frac{1}{\sqrt{2\pi T }} e^{- |x|/2T} dx 
.\end{align}
Here we set $ \Ct (\beta )= \cref{;13}\label{;91}$ for $ \beta = 1,2$ and 
$ \cref{;91} (4) = 2 \cref{;13}$. 
\end{proposition}
 \PF
 We write $ \muN = \muVtN $. We set $ \murNN = \muN \circ \pi _{ \rN }^{-1}$. 
 Recall that the unlabeled dynamics $ \XXrNN $ are $ \murNN $-reversible. 
 Let $ \XrNN = ( \lpath ( \XXrNN ), o , o ,\ldots ) $ be as in the proof of \tref{l:35}. 
 Let $ \widetilde{\pP }^{ \nN } = \pP ^{ \nN } \circ (\XrNN )^{-1}$. 
 We write $ \mathbf{w}=(w^i)_{i\in\N } $. 
Applying the Lyons-Zheng decomposition to $ w^i $, we have that
 \begin{align} \label{:91f}&
 w^i (t) - w^i(0) = \half \{ B^i (\mathbf{w}) (t) + B^i ( \mathcal{\rR }_T (\mathbf{w}))(t) \} 
, \end{align} 
where 
$ \map{\mathcal{\rR }_T }{C([0,T];\RNN )}{C([0,T];\RNN )}$ such that 
$ \mathcal{\rR }_T (\mathbf{w}) (t) = \mathbf{w} (T-t)$. 
The function $ B^i$ is a Brownian motion under $ \widetilde{\pP }^{ \nN } $. 

From \eqref{:11u} and \eqref{:13e}, we have that
\begin{align}\label{:91k}&
\limi{ \nN } 
\frac{1}{\varrho_{V}(\theta )} 
\Vb ' \Big(\frac{ x }{N\varrho_{V}(\theta )} + \theta \Big) = 
\cref{;91}
\end{align} 
uniformly in $ x \in \SR $ for each $ \rR \in \N $. 
Hence, from \eqref{:91f} and \eqref{:91k}, we have the following for each $ \epsilon > 0 $
\begin{align}\label{:91i}&
\limi{ \nN } \widetilde{\pP }^{ \nN } \Big( 
\max_{0 \le t \le T } 
\Big| 
 \int_0^t 
\frac{1}{\varrho_{V}(\theta )} 
\Vb ' \Big(\frac{ w ^{i} (u) }{N\varrho_{V}(\theta )} + \theta \Big) du - 
 \cref{;91} t \Big| \ge \epsilon 
\Big) = 0 
.\end{align}
Because $ \limi{\nN } \rN = \infty $, \eqref{:91f} implies that
\begin{align}\label{:91j}&
\limi{ \nN } \widetilde{\pP }^{ \nN } \Big( 
\max_{0 \le t \le T } 
\Big| 
 \half \int_0^t \mathbf{n}^{ \rN }( w ^{i} ( u ) ) 
L_{ \rN }^{i} (d u ) \Big| \ge \epsilon \Big) = 0 
.\end{align}
Under $ \widetilde{\pP }^{ \nN } $, we can rewrite \eqref{:91b} as 
\begin{align} \notag 
 w^{ i } (t) - w ^{i} (0) 
= 
B ^{i} (t) - & 
 \int_0^t 
\frac{1}{\varrho_{V}(\theta )} 
\Vb ' \Big(\frac{ w ^{i} (u) }{N\varrho_{V}(\theta )} + \theta \Big) du 
 + \frac{\beta }{2}
 \int_0^t \sumN \frac{1}{ w ^{i} ( u ) - w ^{j} (u) } du 
 \\ \notag 
 + & 
\frac{\beta }{2}
 \int_{ \sS \backslash \sS _{ \rN } } \frac{1}{ w ^{i} (t) - y } \rNone ( y ) dy 
 + \half \int_0^t \mathbf{n}^{ \rN }( w ^{i} ( u ) ) 
L_{ \rN }^{i} (d u ) 
.\end{align}
Then, we have that
\begin{align} \label{:91g}
& w^{ i } (t) - w ^{i} (0) + 
 \int_0^t 
\frac{1}{\varrho_{V}(\theta )} 
\Vb ' \Big(\frac{ w ^{i} (u) }{N\varrho_{V}(\theta )} + \theta \Big) du 
 \\ \notag &
 - \frac{\beta }{2}
 \int_0^t \sumN \frac{1}{ w ^{i} ( u ) - w ^{j} (u) } du 
 -
\frac{\beta }{2}
 \int_{ \sS \backslash \sS _{ \rN } } \frac{1}{ w ^{i} (t) - y } \rNone ( y ) dy 
 \\ \notag &
 - 
 \half \int_0^t \mathbf{n}^{ \rN }( w ^{i} ( u ) ) 
L_{ \rN }^{i} (d u ) 
= B ^{i} (t) 
.\end{align}
The boundary of the set 
$\{ w \in C([0,\infty); \R ) \,;\, \max_{0 \le t \le T } w (t) \ge a \} $ 
has Wiener measure zero. 
Combining this with \eqref{:91i}--\eqref{:91g}, and using \tref{l:35}, we obtain 
\begin{align}\notag 
 \limi{ \nN } \widetilde{\pP }^{ \nN } & \Big( \max_{0 \le t \le T } 
 \Big\{
 w^{i}(t) - w^{i}(0) + \cref{;91} t 
 - 
 \frac{\beta }{2} 
 \int_0^t \sum_{ j\neq i }^{ \nN} \frac{1}{w^{i}(u) - w^{j}(u) } du 
 \\ \notag &\quad \quad \quad \quad \quad \quad \quad \quad \quad 
 - 
\frac{\beta }{2}
 \int_{ \sS \backslash \sS _{ \rN } } \frac{1}{ w ^{i} (t) - y } \rNone ( y ) dy 
 \Big\} 
 \ge a \Big) 
 \\ \notag &
= 
P ( \max_{0 \le t \le T } B(t) \ge a)
= 
2 \int_{ x \ge a }\frac{1}{\sqrt{2\pi T }} e^{- |x|/2T} dx 
.\end{align}
Thus, the proof is complete. 
 \PFEND

\subsection{The Ginibre interacting Brownian motion } \label{s:92}

We apply our result to the random matrix model with 
strong non-Hermiticity introduced in \cite{acv}. 
Condition \eqref{:33s} in \As{B4$''$} follows from \pref{l:11}. 
 The SDE for the finite particle system is given by \eqref{:15g} 
with the addition of the following two terms on the right-hand side.:
\begin{align}
 &\notag 
 \int_{ \sS \backslash \sS _{ \rN } }
 \frac{ X_{ \rN }^{ \nN , i } (t) - y }{ | X_{ \rN }^{ \nN , i } (t) - y |^2} \rNone ( y ) dy 
+ 
 \half \int_0^t 
 \mathbf{n}^{ \rN }( X_{ \rN }^{ \nN , i } ( u ) ) 
L_{ \rN }^{ \nN , i } (d u ) 
. \end{align}
Thus, we obtain the results of \tref{l:31}--\tref{l:39} for the Ginibre random point field. 

\section{Acknowledgements}

We thank Stuart Jenkinson, PhD, from Edanz (https://jp.edanz.com/ac) for editing a draft of this manuscript. 
This work was supported by JSPS KAKENHI Grant Numbers JP16H06338, JP20K20885, JP21H04432, and JP21K13812.

\end{document}